\documentclass{amsart}   
\usepackage{amscd,amssymb,amsmath,amsfonts}
\usepackage[matrix, arrow]{xy}  
\usepackage{array} 

\def\intsi{\stackrel{\scriptscriptstyle{o}}{\sigma}} 
\newcommand\Id{\operatorname{Id}}
\newcommand\ad{\operatorname{Ad}}
\newcommand\C{\mathbb C}
\newcommand\K{\mathcal{K}}
\newcommand\E{\mathcal{E}}
\newcommand\A{\mathcal{A}}
\newcommand\de{\delta}
\newcommand\red{\text{red}}
\newcommand\N{\mathbb N}
\newcommand\R{\mathbb R}

\newcommand\mx{\text{max}}
\renewcommand\top{\text{top}}
\newcommand\Z{\mathbb Z}
\newcommand\ts{{\otimes}}
\newcommand\rtm{{\rtimes_{\mx}}}
\newcommand\Ga{{\Gamma}}
\newcommand\Si{{\Sigma}}
\newcommand\ga{{\gamma}}
\newcommand\lto{{\longrightarrow}}
\newcommand\eps{\varepsilon}
\newcommand\x{{X(\Gamma)}}
\newcommand\xin{{X^\infty(\Gamma)}}
\newcommand\ror{{\rho_{P_r(\Ga)}}}
\newcommand\rork{{\rho_{P_r(\Ga)/\Ga_k}}}
\newcommand\rori{{\rho_{P_r(\Ga)/\Ga_i}}}
\renewcommand\a{A_\Gamma}
\newcommand\ai{{A^\infty_\Gamma}}
\newcommand\bi{{B^\infty_\Gamma}}
\newcommand\bio{{B^\infty_{\Gamma,0}}}
\newcommand\cmx{C^*_{max}(X(\Gamma))}
\newcommand\cm{C^*_{max}}
\newcommand\du{\ds\coprod}
\renewcommand\l{\ell}

\newcommand\Ind{\operatorname{Ind}}

\newcommand\diag{\operatorname{diag}}

\newcommand\ds{\displaystyle}

\renewcommand\L{\mathcal{L}}
\newcommand\D{{\mathcal{D}}}

\newcommand\CC{{\mathcal{C}}}

\newcommand\li{{\ell^\infty(X(\Gamma),\mathcal{K}(H))}}

\newcommand\rt{{\rtimes}}

\newcommand\pr{{P_r(\Gamma)}}

\newcommand\I{\operatorname{I}}
\newcommand\res{\operatorname{Res}}
\theoremstyle{plain}
\newtheorem{theorem}{Theorem}[section] 
\newtheorem{proposition}[theorem]{Proposition}
\newtheorem{corollary}[theorem]{Corollary}

\newtheorem{lemma}[theorem]{Lemma}
\newtheorem{remark}[theorem]{Remark}
\newtheorem{example}[theorem]{Example}
\newtheorem{definition}[theorem]{Definition}
\title[$K$-theory for the maximal Roe algebra of certain expanders]
{$K$-theory for the maximal Roe algebra of certain  expanders}
\begin{document}
 \author[H. Oyono-Oyono]{Herv\'e Oyono-Oyono}
 \address{Universit\'e Blaise Pascal \& CNRS, Clermont-Ferrand, France
   and PIMS, University of Victoria, Canada}
 \email{oyono@math.cnrs.fr}
\author[G. Yu]{Guoliang Yu}
 \address{Vanderbilt University, Nashville, USA}
 \email{guoliang.yu@vanderbilt.edu}

\numberwithin{equation}{section} 
\renewcommand{\theenumi}{\roman{enumi}} 
\renewcommand{\labelenumi}{(\theenumi)} 

\maketitle
\begin{abstract}
We study in this paper the maximal version of the coarse Baum-Connes
assembly map for families  of expanding graphs arising from residually finite
groups.  Unlike for  the usual Roe algebra, we show that this assembly map is closely related to the
(maximal) Baum-Connes assembly map for the group and is an isomorphism for a class of expanders. We also introduce a
quantitative Baum-Connes  assembly map and discuss
its connections to  K-theory of (maximal) Roe algebras.
\end{abstract}
\begin{flushleft}{\it Keywords: Baum-Connes Conjecture, Coarse Geometry, Expanders,
  Novikov Conjecture, Operator Algebra K-theory, Roe Algebras}

\end{flushleft}
\tableofcontents
\section{Introduction}
In this paper, we study  K-theory of  (maximal) Roe algebras for a class of expanders. The Roe algebra was introduced by John Roe in his study of higher index theory of elliptic operators  on noncompact spaces \cite{r}.  The K-theory 
of Roe algebra is the receptacle for the higher indices of elliptic
operators. If a space is coarsely embeddable into Hilbert space, then
K-theory of Roe algebra and  higher indices of elliptic operators  are
computatble \cite{y}. Gromov discovered  that expanders do not admit
coarse embedding into Hilbert space \cite{g}.  The purpose of this
paper is to completely  or partially compute K-theory of the (maximal)
Roe algebras associated to certain expanders. In particular, we prove
the  maximal version of the coarse Baum-Connes conjecture for a
special class of expanders. The coarse Baum-Connes conjecture is a
geometric analogue of the Baum-Connes conjecture \cite{bc} and
provides an algorithm of computing  K-theory of Roe algebras and
higher indices of elliptic operators.  We also prove the (maximal)
coarse Novikov conjecture for a class of expanders. The coarse Novikov
conjecture gives  a partial computation of K-theory of  Roe algebras
and an algorithm to determine non-vanishing of higher indices for
elliptic operators. Our results on the coarse Novikov conjecture are
more general than  results obtained in 
\cite{ctwy,gwy,gty}. The question of computing K-theory of (maximal)
Roe algebras associated to general expanders   remains open. We show
that this question is closely related to  certain quantitative Novikov
conjecture and the quantitative Baum-Connes  conjecture for the  K-theory of (maximal) Roe algebras. We explore this connection to prove the quantitative Novikov conjecture and the quantitative 
Baum-Connes  conjecture in some cases.

\smallskip
 The class of expanders under examination in this paper is those
 associated to a finitely generated and  residually finite group $\Ga$
 with respect to a family
 $$\Ga_0\supset\Ga_1\supset\ldots\Ga_n\supset\ldots$$ of finite index
 normal  subgroups.  The behaviour of the Baum-Connes assembly map
 for $\Ga$ and of the coarse Baum-Connes assembly map for the metric
 space   $\x=\coprod_{i\in\N}\Ga/\Ga_i$ can differ quite  substantially: if
 $\Ga$ has the property $\tau$ with respect to the family
 $(\Ga_i)_{i\in\N}$, then $(\Ga/\Ga_i)_{i\in\N}$ is a family of
 expanders and the coarse
 assembly map for $\x$ fails to be an isomorphism,  althougth for example for
 $\Ga=SL_2(\Z)$, the assemply maps is an isomorphism. In \cite{gwy}
 was introduced the maximal Roe algebra of a coarse space and a
 maximal coarse  assembly map with value in this algebra was
 defined. As we shall see, the behaviour of this maximal coarse
 Baum-Connes  assembly map for $\x$, and of the maximal Baum-Connes
 assemply map for the group $\Ga$ with coefficients in
 $\ell^\infty(\x,\K(H))/C_0(\x,\K(H))$ turn out to be equivalent. In
 particular, as a consequence of \cite{tumoy}, if $\Ga$ sastifies
 the strong Baum-Connes conjecture, then the maximal coarse assembly
 map for $\x$ is an isomorphism. As a as a spin-off we also obtain the
 injectivity of the coarse assembly map when $\Ga$ coarsely embeds
 in a Hilbert Space.

This suggests that the properties of the maximal coarse assembly map
for $\x$ is closely related to some asymptotic properties of the 
maximal Baum-Connes assembly maps for $\Ga$ with  coefficients in the 
family $\{C(\Ga/\Ga_i\}_{i\in\N}$. For this purpose, we define 
quantitative assembly maps that take into account the propagation in the
crossed product $\{C(\Ga/\Ga_i)\rtm\Ga\}_{i\in\N}$. Notice that
althought $C(\Ga/\Ga_i)\rtm\Ga$ and $C_\mx^*(\Ga_i)$ are Morita
equivalent, the imprimitivy bimodule between these two algebras
introduces some distortion in the propagation and the relevant
propagation is the one coming from $\Ga$. In this setting, we show
that injectivity and bijectivity of the maximal coarse assembly map are
equivalent to some asymptotic statements for these quantative assembly
maps. For surjectivity, and up to a slight modification in the
sequence of normal subgroups, we  also obtain similar results.

The paper is organised as follows.
In section \ref{sec-roe}, we review results from \cite{gwy} and
\cite{hr}  concerning maximal Roe algebras and coarse assembly
maps. We also show the existence of a short exact sequence (see
section \ref{sub-roe-resfin})
\begin{equation}\label{equ-suite-ex}0\lto \K(\l^2(\x)\ts H)\lto\cmx
  \lto\a\rtm\Ga\lto 0.\end{equation}
In section \ref{sec-assembly-map}, we collect results about Baum-Connes
assembly maps that we use later on.
In  section \ref{sec-khom}, we state for the left hand side of the
Baum-Connes assembly map an  analogue of the exact
sequence of equation \ref{equ-suite-ex} . We show that assembly maps intertwines this exact
sequence with the one induced in K-theory by the exact sequence  \ref{equ-suite-ex}, and
obtain injectivity and bijectivity results for the 
maximal coarse assembly map for $\x$.
In section \ref{sec-asymp}  we set asymptotic  statements
concerning some quantitative assembly maps and we discuss examples of
groups that satisfy these  statements.

\section{K-theory for maximal Roe algebras}\label{sec-roe}
\subsection{Maximal Roe algebra of a locally compact metric space}
In this section, we collect from \cite{gwy}  results concerning the
maximal Roe algebra of a locally compact metric space  that we will
need in this paper.

\subsubsection{The case of a discrete space}
Let $\Sigma$ be a discrete space equipped with a proper distance $d$.
Let us denote by $C[\Sigma]$ the algebra of locally compact operators with
finite propagation of $\l^2(\Sigma)\ts H$, where $H$ is a separable
Hilbert space, i.e (bounded) operators $T$ of $\l^2(\Sigma)\ts H$  such
that when written as a family $(T_{x,y})_{(x,y)\in \Sigma^2}$ of operator
on $H$, then 
\begin{itemize}
\item $T_{x,y}$ is compact for all $x$ and $y$ in $\Sigma$;
\item there exists a real $r$ such that $d(x,y)>r$ implies that
  $T_{x,y}=0$ ($T$ is said to have propagation less than $r$).
\end{itemize}
For any real  $r$, we define $C_r[\Sigma]$ as  the set of elements of 
$C[\Sigma]$ with propagation less than $r$. It is straightforward to check
that $C[\Sigma]$ is a $*$-algebra. The (usual) Roe algebra $C^*(\Si)$ is the
closure of $C[\Sigma]$ viewed as a subalgebra of operator of
$\L(\l^2(\Sigma)\ts H)$.  
The next lemma, proved in \cite{gwy}, shows that if $\Si$ has bounded
geometry, then $C[\Sigma]$ admits an
envelopping algebra.
 
\begin{lemma}\label{lem-norm} Let $\Si$ be a discrete metric space
  with bounded geometry. For any positive number $r$, there
  exists a  real $c_r$
 such that 
for any $*$-representation $\phi$  of $C[\Sigma]$ on a Hilbert
space $H_\phi$ and any $T$ in $C_r[\Sigma]$, then  $\|\phi(T)\|_{\L(H_\phi)}\leq c_r\|T\|_{\l^2(\Sigma)\ts H}$.
\end{lemma}
This envelopping algebra is then

\begin{definition} \cite{gwy} The maximal Roe algebra of a discrete
  metric 
  space $\Sigma$ with bounded geometry,
  denoted by $C^*_\text{max}(\Sigma)$, is the
  completion of $C[\Sigma]$ with respect to the $*$-norm 
$$\|\phi(T)\|=\sup_{(\phi,H_\phi)}\|\phi(T)\|_{\L(H_\phi)},$$ when
$(\phi,H_\phi)$ runs through representations $\phi$ of $C[\Sigma]$  on a Hilbert space
$H_\phi$.
\end{definition}

\subsubsection{The general case}
Let $X$ be a locally compact space, equipped with a metric $d$.   A
$X$-module  is a Hilbert space $H_X$ together with a $*$-representation
$\rho_X$ of $C_0(X)$ in $H_X$. We shall often write $f$ instead of
$\rho_X(f)$ for the action of $f$ on $H_X$. If the representation is non-degenerate, the
$X$-module is said to be non-degenerate.
A $X$-module is called standard if no non-zero function of $C_0(X)$
acts as a compact operator on $H_X$. In the litterature, the
terminology $C_0(X)$-ample is also used for such a
representation \cite{hr,tzanev}.
\begin{definition}
Let $H_X$ be a standard non-degenerate $X$-module and let $T$ be a
bounded operator on $H_X$.
\begin{enumerate}
\item The support of $T$ is the complement of the open subset of
  $X\times X$
\begin{equation*}\begin{split}\{(x,y)\in X\times X\text{ s.t.  there exist }& f\text{ and }
g\text{ in }C_0(X)\text{  satisfying }\\
&f(x)\neq 0,\,g(y)\neq 0\text{
  and }f\cdot T\cdot g=0\}.\end{split}\end{equation*}
\item If there exists a real $r$ such that for any $x$ and $y$ in $X$
  such that $d(x,y)>r$, then $(x,y)$ is not in the support of $T$,
  then the operator $T$  is said to
  have finite propagation (in this case    propagation less than $r$).
\item The operator $T$ is said to be locally compact if $f\cdot T$ and
  $T\cdot f$ are compact for any $f$ in
  $C_0(X)$. We then define   $C[X]$ as the set of locally compact and finite
  propagation bounded operators of $H_X$.
\item The operator $T$ is said to pseudo-local if $[f,T]$ is compact
  for all $f$ in $C_0(X)$.
\end{enumerate}
\end{definition}
It is straightforward to check
that $C[X]$ is a $*$-algebra and that for a discrete space, this
definition coincides with the previous one. Moreover, up to
(non-canonical) isomorphism, $C[X]$ does not depend on the choice of
$H_X$.
The  Roe algebra $C^*(X)$ is then the norm closure of $C[X]$ in 
the algebra $\L(H_X)$ of bounded operators on $H_X$. Although $C^*(X)$
is not canonically defined, we shall see later on that up to canonical
isomorphism, its K-theory
does not depend on the choice a
non-degenerated standard $X$-module.

\begin{definition}\label{def-net} A net in a 
 locally compact space  $X$ is a countable  
subset $\Sigma$ such that there exists 
numbers $\eps$ and $r$ satisfying
\begin{itemize}
\item $d(y,y')\geq \eps$ for any distinct elements $y$ and $y'$ of $\Sigma$;
\item For any $x$ in $X$, there exists $y$ in 
$\Sigma$ such that $d(x,y)\leq r$.
\end{itemize}
\end{definition}
The following result was proved in 
\cite{gwy}
\begin{lemma}\label{lem-coarse-contain}
If a locally compact space $X$ contains a net with bounded geometry, then with notation of definition \ref{def-net},  there exists a unitary map $\Psi:H_\Si\to H_X$ that
fullfills the following conditions:
\begin{enumerate}
\item The homomorphism $\L(H_\Si)\to \L(H_X);T\mapsto \Psi^*\cdot T\cdot
  \Psi$
restricts to an algebra $*$-isomorphism 
$C[\Si]\to C[X]$;
\item There exists a number $r$ such that for every $x$ in $X$ and $y$
  in $\Si$ with  $d(x,y)<r$, then there exists $f$ in $C_0(X)$ and
  $g$ in $C_0(\Si)$ which satisfy $f(x)\neq 0$, $g(y)\neq 0$ and
  $f\cdot\Psi\cdot g=0$ (i.e $\Psi$ has propagation less than $r$).
\end{enumerate}
\end{lemma}
Then, if $\Psi:H_\Si\to H_X$ is a unitary map as in lemma
\ref{lem-coarse-contain},  the   $*$-isomorphism 
$C[\Si]\to C[X];T\mapsto \Psi^*\cdot T\cdot
  \Psi$ extends to an isomorphism $Ad_\Psi:C^*(\Si)\to C^*(X)$.

As a consequence of lemma \ref{lem-coarse-contain}, lemma \ref{lem-norm} admits the following generalisation to
spaces  that  contain
a net with bounded geometry.

\begin{lemma}\label{lem-norm-net} Let $X$ be a locally compact metric
  that contains a net with bounded geometry  and let $H_X$ be a standard non-degenerate
  $X$-module. Then for any positive number $r$,  there
  exists a  real $c_r$
 such that 
for any $*$-representation $\phi$  of $C[X]$ on a Hilbert
space $H_\phi$ and any $T$ in $C[X]$ with propagation less than $r$, then  $\|\phi(T)\|_{\L(H_\phi)}\leq c_r\|T\|_{H_X}$.
\end{lemma}

This allowed to define for $X$ the maximal Roe algebra as in the
discrete case.

\begin{definition} \cite{gwy} Let $X$ be a locally compact metric 
  space that   contains a  net with bounded geometry . The maximal Roe algebra of $X$,
  denoted by $C^*_\text{max}(X)$, is the
  completion of $C[X]$ with respect to the $*$-norm 
$$\|\phi(T)\|=\sup_{(\phi,H_\phi)}\|\phi(T)\|_{\L(H_\phi)},$$ when
$(\phi,H_\phi)$ runs through representation $\phi$ on of $C[X]$  a Hilbert space
$H_\phi$.
\end{definition}

\subsection{Maximal Roe algebra associated to a residually finite
  group}\label{sub-roe-resfin}
Let $\Ga$ be a residually finite group, finitely generated. Let
$\Ga_0\supset\Ga_1\supset\ldots\Ga_n\supset\ldots$ be a decreasing
sequence of finite index subgroups of $\Gamma$ such that
$\bigcap_{i\in\N}\Ga_i=\{e\}$. Let $d$ be a left invariant metric
associated to a given  finite set of generators. Let us endow $\Ga/\Ga_i$
with the metric $d(a\Ga_i,b\Ga_i)=\min\{d(a\ga_1,b\ga_2), \ga_1\text{
  and }\ga_2 \text{
  in } \Ga_i\}$.
We set $\x=\du_{i\in\N}\Ga/\Ga_i$ and we equip $\x$ with a metric $d$
such that,
\begin{itemize}
\item  on  $\Ga/\Ga_i$, then $d$ is the metric defined above;
\item 
$d(\Ga/\Ga_i,\Ga/\Ga_i)\geq i+j$ if $i\neq j$.
\item the group $\Ga$ acts on $\x$ by isometries. 
\end{itemize}

Let us define by $\K(H)$ the algebra of compact operators of $H$. Then
the $C^*$-algebra
$\l^\infty(\x,\K(H))$ acts on $\l^2(\x)\ts H$ by pointwise
action of $\K(H)$. This action is clearly by propagation zero locally
compact operators.
The group $\Ga$ acts diagonally on the Hilbert space $\l^2(\x)\ts H$
by finite propagation operators,
the action being on $\l^2(\x)$  induced by the action  on $\x$
and  trivial on $H$. From this, we get a covariant representation of
$(\li,\Ga)$ on $\l^2(\x)$, where the  action of $\Ga$ on $\l^\infty(\x,\K(H))$ is
induced by the action of $\Ga$ on  $\x$ by translations. This yields to a
$*$-homomorphism
$C_c(\Ga,\l^\infty(\x,\K(H)))\lto C[\x]$ (where
$C_c(\Ga,\l^\infty(\x,\K(H)))$ is equipped with the convolution
product) and thus, setting  $B_\Ga=\l^\infty(\x,\K(H))$,  to a $*$-homomorphism
$$B_\Ga\rtm\Ga\lto \cmx.$$ Under this map, the image of
$B_{\Ga,0}\stackrel{\text{def}}{=\!=}C_0(\x,\K(H))\rtm\Ga$ is contained in $\K(\l^2(\x)\ts H)$. Thus if we
set  $\a=B_\Ga/B_{\Ga,0}$, then we finally get a
$*$-homomorphism $$\Phi_{\Ga}:\a\rtm\Ga\lto \cmx/\K(\l^2(\x)\ts H).$$
\begin{proposition}\label{prop-roemax}
  $\Phi_\Ga$ is a $*$-isomorphism, i.e we have an short exact sequence
$$0\lto \K(\l^2(\x)\ts H)\lto\cmx \lto\a\rtm\Ga\lto 0.$$
\end{proposition}
\begin{proof}
Let us construct an inverse for $\Phi_\Ga$. Let $T$ be an element in
$C[\x]$ with propagation less than $r$. Let $n$ be any integer such
that $n\geq r$ and $B_\Ga(e,r)\cap \Ga_n=\{e\}$. Then there is a
decomposition $T=T'+T''$ with $T'$ in
$\K(\l^2(\coprod_{i=0}^{n-1}\Ga/\Ga_i)\ts H)$ and
$T''=(T''_i)_{i\geq n}\in\prod_{i\geq n} \K(\l^2(\Ga/\Ga_i)\ts
H)$. Let us denote for $\ga$ in $\Ga$ by $L_{\ga\Ga_i}$ the diagonal
  operator on $\l^2(\Ga/\Ga_i)\ts H$ given  by  left translation by
  $\ga\Ga_i$ on  $\l^2(\Ga/\Ga_i)$ and the identity  on $H$. For any
  integer $i$,  we have
  a unique decomposition $T''_i=\sum_{k=1}^p h_kL_{\ga_k\Ga_i}$, where
    $\ga_k\Ga_i$ belongs to $B_{\Ga/\Ga_i}(\Ga_i,r)$ and $h_k$ belongs
    to $C(\Ga/\Ga_i,\K(H))$ and is  viewed as an operator
    acting on $\l^2(\Ga/\Ga_i)\ts H$ by pointwise action of  operators of
    $\K(H)$. Since $B_\Ga(e,r)\cap \Ga_i=\{e\}$ for $i\geq n$, the
    element $\ga_k\Ga_i$ lifts to a unique element of
    $B_\Gamma(e,r)$. If $B_\Gamma(e,r)=\{g_1,\ldots,g_m\}$, then there
    is a unique decomposition $T''_i=\sum_{k=1}^m f_k^iL_{g_k\Ga_i}$,
    with $f_k^i$ in $C(\Ga/\Ga_i,\K(H)$. Let us denote for
    $k=1,\ldots,m$ by $\phi_k(T)$ the image of $( f_k^i)_{i\geq n}$
    under the projection $$\prod_{i\geq n} C(\Ga/\Ga_i,\K(H))\lto
    \prod_{i\geq n} C(\Ga/\Ga_i,\K(H)/\bigoplus_{i\geq n}
    C(\Ga/\Ga_i,\K(H)\cong \a.$$ It is then straightforward to check
    that we obtain in this way  a well defined map
$$\Lambda_r:C_r[\x]\lto C_c(\Ga,\a);T\mapsto \sum_{k=1}^m
\phi_k(T)\de_{g_k},$$  where $\de_g$   is the Dirac
function at an element $g$  of  $\Ga$. Moreover, if $r'\geq r$,then $\Lambda_{r'}$ restricts
to $\Lambda_r$ on $C_r[\x]$ and the maps $\Lambda_r$ extends to a
$*$-homomorphism
$C[\x]\lto C_c(\Ga,\a)$ and thus to a $*$-homomorphism $
\cmx\lto\a\rtm\Ga$. This homomorphism clearly factorizes through a
$*$-homomorphism 
$$
\cmx/\K(\l^2(\x)\ts H)\lto\a\rtm\Ga$$ which provides an inverse for
$\Phi_\Ga.$

\end{proof}
We shall denote by $\Psi_{\a,\Ga,\text{max}}:\cmx \lto\a\rtm\Ga$ the
projection map corresponding to the exact sequence of the previous
proposition. Let
$\lambda_{\Ga,A_\Ga}:\a\rtm\Ga\to\a\rtimes_{\text{red}}\Ga$ be the
homomorphism given
by the regular representation of the covariant system
$(A_\Ga,\Ga)$. The next lemma shows that
$\lambda_{\Ga,A_\Ga}\circ\Psi_{\a,\Ga,\text{max}}$ factorizes through
$C^*(\x)$ (see \cite{hls}).

\begin{lemma}
There exists a unique homomorphism $$\Psi_{\a,\Ga,\text{red}}:C^*(\x)
\to\a\rtimes_{\text{red}}\Ga$$ such that  $\lambda_{\x}\circ\Psi_{\a,\Ga,\text{red}}=\Psi_{\a,\Ga,\text{max}}\circ\lambda_{\Ga,A_\Ga}$. 
\end{lemma}
\begin{proof}
If such an homomorphism exists, it is clearly unique. Let us prove the
existence. Let $T$ be an element of $C_r(\x)$ such that
$\|T\|_{\L(\ell^2(\x)\ts H))}=1$ and let us set $x_T=\lambda_{\Ga,A_\Ga}\circ\Psi_{\a,\Ga,\text{max}}(T)$. Let us view $\a\rtimes_\text{red}\Ga$ as
an algebra of adjointable operator on the right $\a$-Hilbert module
$\a\ts\ell^2(\Ga)$. For any positive real $\eps$, let $\xi$ be an
element of $\a\ts\ell^2(\Ga)$ such that $\|\xi\|_{\a\ts\ell^2(\Ga)}=1$
and
$$\|x_T.\xi\|_{\a\ts\ell^2(\Ga)}\geq\|x_T\|_{\a\rtimes_\text{red}\Ga}-\eps.$$We
can  assume without loss of generality that $\xi$ lies indeed in
$C_c(\Ga,\a)$ and is supported in some $B_\Ga(e,s)$ for $s$ positive
real. Let
us fix an integer $k$ such that
$B_{\Ga}(e,2r+2s)\cap\Ga_k=\{e\}$. There is  a decomposition
$T=T'+T''$ with $T'$ in $\K(\ell^2(\sqcup_{i=1}^{k-1}\Ga/\Ga_i)\ts H)$
and $T''$ in ${\displaystyle \prod_{i\geq k} \K(\ell^2(\Ga/\Ga_i)\ts H)\cong
\prod_{i\geq k} C(\Ga/\Ga_i,\K(H))\rtimes \Ga/\Ga_i}$. Since this
decomposition is diagonal, we get $\|T''\|\leq 1$. Actually $T''$ can
be viewed as an adjointable operator on the right $\ell^\infty(\cup_{i\geq
  k}\Ga/\Ga_i,\K(H))$-Hilbert module $\prod_{i\geq k}
C(\Ga/\Ga_i,\K(H))\otimes \ell^2(\Ga/\Ga_i)$. Let us chose a lift
$\xi'$ in $C_c(\Ga,B_\Ga)$  of $\xi$ under the map
 induced by the canonical
projection $B_\Ga\to\a$ such that
$\|\xi'\|_{B_\Ga\otimes\ell^2(\Ga)}\leq 1+\eps$ and $\xi'$ is supported
in $B_\Ga(e,s)$. Under the identification $$C_c(\Ga,\prod_{i\geq k}
C(\Ga/\Ga_i,\K(H)))\cong \prod_{i\geq k} C_c(\Ga,C(\Ga/\Ga_i,\K(H))),$$
we can write $\xi'=(\xi'_i)_{i\geq k}$. Since $B_\Ga(e,2r)\cap
\Ga_i=\{e\}$ for $i\geq k$, the map $$B_\Ga(e,r)\to
B_{\Ga/\Ga_i}(\Ga_i,r);\ga\to \ga\Ga_i$$ is bijective. Hence, if we define 
  $\xi''_i$ in $C_c(\Ga/\Ga_i,\K(H))$ with support in
  $B_{\Ga/\Ga_i}(\Ga_i,r)$ by $\xi''_i(\ga\Ga_i)=\xi'_i(\ga)$ for any integer
    $i\geq k$ and $\ga$ in $B_{\Ga}(e,r)$, then  $\xi''=(\xi''_i)_{i\leq k}$ is an element of
    $\prod_{i\geq k} C(\Ga/\Ga_i,\ell^2(\Ga/\Ga_i))$ such that
    $\|\xi''\|=\|\xi'\|\leq 1+\eps$. Let us now define
    $\eta''=(\eta''_i)_{i\geq k}$ in $\prod_{i\geq
      k}C(\Ga/\Ga_i,C(\Ga/\Ga_i,\K(H)))$ by
    $\eta''=T''\cdot\xi''$. Then $\eta''_i\in
    C(\Ga/\Ga_i,C(\Ga/\Ga_i,\K(H)))$ has support in
    $B_{\Ga/\Ga_i}(\Ga_i,r+s)$. Let us  define $\eta'_i\in
    C(\Ga,C(\Ga/\Ga_i,\K(H))$ with support in $B_\Ga(e,r+s)$ by
    $\eta_i'(\ga)=\eta_i''(\ga\Ga_i)$ for $i\geq k$ and $\ga$ in
    $B_\Ga(e,r+s)$. Since  $B_\Ga(e,2r+2s)\cap
\Ga_i=\{e\}$ for $i\geq k$, the map $$B_\Ga(e,r+s)\to
B_{\Ga/\Ga_i}(\Ga_i,r+s);\ga\to \ga\Ga_i$$ is bijective and thus
$\|\eta''\|=\|\eta'\|$. It is then straightforward to check that the
image of $\eta'$ under the map $C_c(\Ga,B_\Ga)\to C_c(\Ga,\a)$ induced
by the canonical projection $B_\Ga\to\a$  is precisely $
x_T.\xi$. The result is
then a consequence of the following:
\begin{eqnarray*}
\|x_T.\xi\|&\leq&\|\eta''\|\\
&\leq&\|\xi''\|\\
&\leq&1+\eps\\
\end{eqnarray*}

\end{proof}


\begin{proposition}\label{prop-inj-cpt-max}
The inclusion $\K(\l^2(\x)\ts H)\hookrightarrow\cmx$ induced an
injection
$\Z\hookrightarrow K_0(\cmx)$
\end{proposition}
\begin{proof}
Let $p$ be a projector in $\K(\l^2(\x)\ts H)$ such that $[p]=0$ in
$K_0(\cmx)$. We can  assume without loss of generality that $p$ belongs to
$\K(\l^2(\coprod_{i=1}^n\Ga/\Ga_i)\ts H)$ for some $n$. This means that
$\begin{pmatrix} p&0&0\\
0&I_k&0\\
0&0&0
\end{pmatrix}$  is homotopic to $\begin{pmatrix} I_k&0&0\\
0&0&0\\
0&0&0
\end{pmatrix}$ in some $M_N(\widetilde{\cmx})$ by a homotopy of
projectors. Hence for every positive number $\eps<1/4$ there exists a
real $r$ such that $\begin{pmatrix} p&0&0\\
0&I_k&0\\
0&0&0
\end{pmatrix}$  and  $\begin{pmatrix} I_k&0&0\\
0&0&0\\
0&0&0
\end{pmatrix}$ are connected by a homotopy $(p_t)_{t\in[0,1]}$ of
$\eps$-almost projectors of  $M_N(\widetilde{\cmx})$ (i.e  selfadjoint
elements  satisfying   $\|p_t^2-p_t\|\leq \eps$)  such that 
$p_t$ has propagation less than $r$
for all $t$ in $[0,1]$.
Let us fix an integer $k\geq \max\{n,r\}$. Then for every  $t\in[0,1]$ we can
write $p_t=p'_t+ p''_t$, where
$(p'_t)_{t\in[0,1]}$ is a homotopy of selfadjoint elements in
$$M_N(\K(\l^2(\coprod_{i=1}^{k-1}\Ga/\Ga_i)\ts H)+\C
Id_{\l^2(\coprod_{i=1}^{k-1}\Ga/\Ga_i)\ts H}))$$  and   $(p''_t)_{t\in[0,1]}$ is a homotopy of selfadjoint elements in
$$M_N\left(\prod_{i\geq k} \K(\l^2(\Ga/\Ga_i)\ts H)+\C
Id_{\l^2(\coprod_{i\geq k}\Ga/\Ga_i)\ts H}\right).$$ 
Moreover, since
$p_t$ can be written diagonally as $p'_t\oplus p''_t$ in the decomposition
$\l^2(\x)\ts H=\l^2(\coprod_{i=1}^{k-1}\Ga/\Ga_i)\ts H\bigoplus
\l^2(\coprod_{i\geq k}\Ga/\Ga_i)\ts H$, then $p'_t$ is also a
$\eps$-projectors with propagation less than $r$. Let $\varphi:\R\to\R$ be a 
continuous function such that $\varphi(s)=0$ if $s\leq 1/2$ and
$\varphi(s)=1$ if $\frac{\sqrt{1-4\eps}+1}{2}\leq s$. Then
$(\varphi(p'_t))_{t\in[0,1]}$ is a homotopy of projectors in $M_N(\K(\l^2(\coprod_{i=1}^{k-1}\Ga/\Ga_i)\ts H)+\C
Id_{\l^2(\coprod_{i=1}^{k-1}\Ga/\Ga_i)\ts H})$ between
$\begin{pmatrix} p&0&0\\
0&I_k&0\\
0&0&0
\end{pmatrix}$  and $\begin{pmatrix} I_k&0&0\\
0&0&0\\
0&0&0
\end{pmatrix}$ and thus $p=0$.
\end{proof}
In conclusion, we get
\begin{corollary}\label{cor-suitex}
With previous notations, we have 
\begin{enumerate}
\item  a short exact sequence
$$0\lto \Z\lto K_0(\cmx)\stackrel{\Psi_{\a,\Ga,\text{max},*}}{\lto}  K_0(\a\rtm\Ga)\lto 0,$$
where the copy of $\Z$ in $K_0(\cmx)$ is generated by the class of a
rank one projector on $\ell^2(X)\otimes H$.
\item an isomorphism $$K_1(\cmx) \stackrel{\Psi_{\a,\Ga,\text{max},*}}{\lto}
  K_1(\a\rtm\Ga)$$
\end{enumerate}
\end{corollary}

\begin{remark}\label{prop-inj-cpt} The same proof  as for proposition
  \ref{prop-inj-cpt-max} applies to show that the injection
$\K(\ell^2(\x)\ts H)\hookrightarrow C^*(\x)$ induces an  injection
$\Z\hookrightarrow K_0(C^*(\x))$. But when the group $\Gamma$ has the
property $\tau$ with respects to the family $(\Ga_i)_{i\in\N}$, then
$\x$ is   a family of
expanders and it was proved in \cite{hls} that  the composition  $$\Z\lto
K_0(C^*(\x))\stackrel{\Psi_{\Ga,\text{red}},*}{\lto}K_0(\a\rtimes_{\text{red}}\Ga)$$
(and thus the composition $\K(\ell^2(\x)\ts H)\hookrightarrow C^*(\x)\stackrel{\Psi_{\Ga,\text{red}}}{\lto}
\a\rtimes_{\text{red}}\Ga$)
is not exact in the middle 
\end{remark}

\subsection{Assembly map for the Maximal Roe algebra}\label{subsection-roe}
Let $X$ be a locally compact metric space, then according to a result
of 
\cite{hr} that we shall recall  below,  the K-theory group
$K_*(C^*(X))$ up to canonical isomorphism does not
depend on the  chosen non-degenerated standard
$X$-module defining $C^*(X)$.
In \cite{hr} was defined an assembly map for the  Roe
algebra
$\widehat{\mu}_X:K_*(X)\to K_*(C^*(X))$ in the following way. Let
$H_X$ be
 a  standard $X$-Hilbert module with respect to a non-degenerated
 representation $\rho_X:C_0(X)\to \L(H_X)$.
Let us define the following subalgebras of $\L(H_X)$
\begin{eqnarray*}
\D(X)&=&\{T\in \L(H_X)\text{ such that } [f,T]\in\K(H_X)\text{ for all
} f\in C_0(X)\},\\
\CC(X)&=&\{T\in \L(H_X)\text{ such that } f\cdot T \in\K(H_X)\text{ and
  }T\cdot f \text{ for all
} f\in C_0(X)\},\\
D^*(X)&=&\{T\in \D(X) \text{ and } T \text{ is in the closure of
  finite propagation operator} \}.
\end{eqnarray*}
Then every element in $K_*(X)$ can be represented by a  $K$-cycle
$(\rho_X,H_X,T)$, with $T\in\D(X)$. This operator then defines a class
$[T]$ in  $K_{*+1}(\D(X)/\CC(X))$  and  we get in this way an isomorphism
called the Paschke duality \cite{hr}
$$K_*(X)\lto K_{*+1}(\D(X)/\CC(X)); [(\rho_X,H_X,T)]\mapsto [T].$$
According to  \cite[Lemma 12.3.2]{hr}, the
$C^*$-algebras inclusions $C^*(X)\hookrightarrow \CC(X)$ and
$D^*(X)\hookrightarrow \D(X)$  induce  an isomorphism 
\begin{equation}\label{iso-corona}
D^*(X)/C^*(X)\stackrel{\cong}{\lto}\D(X)/\CC(X).
\end{equation} Using the inverse
of this isomorphism, we get finally an isomorphism
$$K_*(X)\stackrel{\cong}{\lto} K_{*+1}(D^*(X)/C^*(X))$$
  which when composed with
the boundary map  in $K$-theory 
associated to the short exact sequence
$$0\to C^*(X)\to D^*(X)\to D^*(X)/C^*(X)\to 0$$ gives rise to the 
assembly map $$\widehat{\mu}_{X,*}:K_*(X)\lto K_*(C^*(X))$$
\begin{remark}\label{rem-ind}
Every element $x$ in  $K_*(X)$ can be indeed, represented by a  K-cycle
$(\rho_X,H_X,T)$, with $T\in D(X)$. In this case, $(T,C^*(X))$ is
$K$-cycle for $K_*(C^*(X))=KK_(\C,C^*(X))$ and thus defines an element
of $K_*(C^*(X))$ we shall denote by $\Ind_X T$. It is then
straightforward to check that $\widehat{\mu}_{X,*}(x)=\Ind_X  T$.
\end{remark}
In order to define the coarse Baum-Connes assembly maps, we shall
recall some functoriality results of the Roe algebras under coarse
maps.

\medskip

Let $\phi:X\mapsto Y$ be a coarse map between locally compact metric
spaces. Let $H_X$ (resp. $H_Y$) be a non-degenerated standard
$X$-Hilbert module (resp. $Y$-Hilbert module) with respect to a
representation $\rho_X$ (resp. $\rho_Y$). Recall from \cite{hr} that
there is  an isometry $V:H_X\to H_Y$ that covers $\phi$, i.e there exists
a real $s$ such that for any $x$ and $y$ in $X$ with 
$d(\phi(x),y)>s$, we can find $f$ in $C_0(Y)$ and 
$g$ in $C_0(X)$ that satisfy $f(y)\neq 0$, $g(x)\neq 0$ and
$\rho_Y(f)V\rho_X(g)=0$.  The map 
$\L(H_X)\to\L(H_Y);\,T\mapsto VTV^*$ then restricts  to a $*$-homomorphism
$C[X]\to C[Y]$ and thus to a homomorphism
$\ad V:C^*(X)\to C^*(Y)$.
The crucial point, due to \cite{hr} is that the homomorphism $\ad_*
V:K_*(C^*(X))\to K_*(C^*(Y))$ induced in K-theory by $\ad V$ does not
depend on the choice of the isometry $V$ covering $\phi$. Hence, we define $\phi_*=\ad_* V: K_*(C^*(X))\to K_*(C^*(Y))$, where $V:H_X\to
H_Y$ is any isometry covering $\phi$. 
\begin{remark}\label{rem-funct}\
\begin{enumerate}
\item If $\phi:X\to X$ is a coarse  map such that for some real $C$,
  we have $d(x,\phi(x)< C$ for all $x$ in $X$, then $\Id_{H_X}$ covers
  $\phi$ and hence $\phi_*=\Id_{K_*(C^*(X))}$.
\item If $\phi:X\to Y$ and $\psi:Y\to Z$ are two  coarse maps,
  then $(\phi\circ\psi)_*=\phi_*\circ\psi_*$.
\item In consequence  if $\phi:X\to Y$ is a coarse
  equivalence,
 then $\phi_*$ is an isomorphism. Moreover, if $\phi':X\to Y$ is
 another  coarse
  equivalence, then $\phi_*=\phi'_*$.
\item In particular, by chosing for   two non-degenerated
  standard $X$-modules $H_X$ and $H'_X$ an isometry  $V:H_X\to H'_X$  that
  covers $\Id_X$, we see that up to canonical isomorphisms, the
  K-group $K_*(C^*(X))$ does not depend on the choice of a non-degenerated
  standard $X$-module.
\end{enumerate}
\end{remark}
The previous construction can be extended to maximal Roe algebras by
using the following lemma.
\begin{lemma}
With above notations, assume that  $X$ and $Y$ a locally compact
metric spaces both containing nets with bounded geometry. For any isometry   $V:H_X\to
H_Y$ covering a coarse map $\phi:X\to Y$, we have
\begin{enumerate}
\item $C[X]\to C[Y];\,T\mapsto VTV^*$ extends to  a homomorphism
$C[X]\to C[Y]$ and thus to a homomorphism
$\ad_{\text{max}} V:C^*_{\text{max}}(X)\to C^*_{\text{max}}(Y)$.
\item The homomorphism $\ad_{\text{max},*} V:K_*(C^*_{\text{max}}(X))\to
  K_*(C^*_{\text{max}}(Y))$ induced by $\ad_{\text{max}} V$ in K-theory does not
  depend on the choice of the isometry   $V:H_X\to
H_Y$ covering $\phi$.
\end{enumerate}
\end{lemma}
\begin{proof}
The first item is just a consequence of the universal properties for
$C^*_{\text{max}}(X)$ and $C^*_{\text{max}}(Y)$. For the second point, assume
that $V':H_X\to
H_Y$ is another isometry covering $\phi$  and let us set $W=V\cdot V'^*$. Then $W$ is
 a partial isometry of $H_Y$ with finite propagation.
Then $Wx$ and $(\Id_{H_X}-W^*W)x$ are in $C[Y]$  for any $x$ in $C[Y]$,  and since
\begin{eqnarray*}
 x^*x
&=&x^*W^*Wx+x^*(\Id_{H_X}-W^*W)^2x\\
&=& (Wx)^*Wx+(x(\Id_{H_X}-W^*W))^*(\Id_{H_X}-W^*W)x,
\end{eqnarray*}
we get that  $(Wx)^*Wx\leq x^*x$ in $C^*_{\text{max}}(Y)$ and thus
$\|Wx\|\leq \|x\|$. Hence 
 $C[Y]\to
C[Y];x\mapsto W\cdot x$ extends to  a bounded linear map
$C^*_\text{max}(X)\to C^*_\text{max}(X)$ we shall denote again by
$W$. But $W^*$ is also a partial isometry of $H_Y$ with finite support
and it is straightforward to
check that for any $x$ and $y$ in $C^*_\text{max}(Y)$ then
${W^*}(x)y=x W(y)$ and thus $W$ is a multiplier of
$C^*_\text{max}(Y)$. Since $\ad_{\text{max},*} V=W\cdot
\ad_{\text{max},*} V\cdot W^*$, we get the result by using \cite[Lemma 4.6.2]{hr}.
\end{proof}
 With above notation, this allowed to define for   a coarse map
 $\phi:X\to Y$ 
 with $X$ and $Y$ both containing nets with bounded geometry,    
 $$\phi_{\text{max},*}=\ad_{\text{max},*} V:K_*(C^*_\text{max}(X))\lto K_*(C^*_\text{max}(Y)),$$ where $V:H_X\to
H_Y$ is any isometry covering $\phi$.
Notice that the remark \ref{rem-funct}  obvioulsy admits a analogous
formulation for  maximal Roe algebra.
In view of remark  \ref{rem-ind} and  in order to define the
maximal assembly map, we will need the following result.

\begin{lemma}\label{lem-pseudo-mult}
Let $H_X$ a standard non-degenerated $X$-module. Then, for any
pseudo-local operator  $T$ on $H_X$ with finite propagation, the map $$C[X]\to C[X];\,x\mapsto Tx$$ extends in a unique way to a multiplier of
  $C^*_\text{max}(X)$ we shall again denote by $T$. Moreover, for any 
positive real $r$, there exists a real constant $c_r$ such that if $T$
has propagation less than $r$, then 
 $$\|T\|_{M(C^*_\text{max}(X)) }\leq c_r \|T\|_{\L(H_X) },$$ where 
  $M(C^*_\text{max}(X))$ stands for the algebra of multiplier of $C^*_\text{max}(X)$. 
\end{lemma}
\begin{proof}We can assume without loss of generallity that $T$ as an
  operator on $H_X$ is norm $1$.
\begin{itemize}
\item
Let $(f_i)_{i\in I}$ be a partition of unit, whose supports have
uniformally bounded diameters. Since $T-\sum_{i\in I}f_i^{1/2}T
f_i^{1/2}$ is in $C[X]$ and has norm operator for $H_X$ less than $2$,
then according to lemma \ref{lem-norm-net} it is enought  to prove the result for 
$\sum_{i\in I}f_i^{1/2}T
f_i^{1/2}$ instead $T$.
For an element $x$ of $C[X]$, let us set  $$x'=
\sum_{i\in I}f_i^{1/2}T
f_i^{1/2}x,$$ $$A=(\Id_{H_X}-T^*T)^{1/2}$$ and $$y=(\sum_{i\in I}f_i^{1/2}A
f_i^{1/2})^*(\sum_{i\in I}f_i^{1/2}A
f_i^{1/2})-\sum_{i\in I}f_i^{1/2}A^*A
f_i^{1/2}.$$ 
Then $A$ is a positive pseudo-local operator of $H_X$ and  $y$ is a
self-adjoint element in $C[X]$ and
\begin{eqnarray*}
x'^*x'-x^*x&=& \sum_{i\in I}x^*f_i^{1/2}(T^*T-\Id_{H_X})
f_i^{1/2}x+x^*yx\\
&=& -(\sum_{i\in I}f_i^{1/2}A
f_i^{1/2}x)^*(\sum_{i\in I}f_i^{1/2}A
f_i^{1/2}x)+x^*(y+y')x,
\end{eqnarray*}
where $y'=(\sum_{i\in I}f_i^{1/2}A
f_i^{1/2})^2-\sum_{i\in I}f_i^{1/2}A^2
f_i^{1/2}$ lies in  
$C[X]$. According to lemma \ref{lem-norm-net},  since $y$ and $y'$ has
operator norm on $H_X$ bounded by $2$ and have propagation less than
$r$,  there exists a real
constant $c'_r$, depending only on $H_X$ and $r$ and such that
$\|y+y'\|\leq c'_r$ in  $C^*_\text{max}(X)$. Hence 
$x'^*x'-x^*x\leq c_r x^*x$ and hence $\|x'\|\leq (1+c_r)^{1/2}\|x\|$ in
$C^*_\text{max}(X)$.
In consequence, the map $C[X]\to C[X];\,x\mapsto \sum_{i\in I}f_i^{1/2}T
f_i^{1/2}x,$ is bounded for the norm of  $C^*_\text{max}(X)$ and thus
extends to a bounded linear map $C^*_\text{max}(X)\to
C^*_\text{max}(X)$.
\item Applying the preceding point also to $T^*$, we get that
$(T^*x)^*y=x(Ty)$ for all $x$ and $y$ in $C^*_\text{max}(X)$ (check it on $C[X]$) and thus $T$ is a multiplier for
 $C^*_\text{max}(X)$.
\end{itemize}
\end{proof}
\begin{remark}
The set of pseudo-local operator of $H_X$ of finite propagation is a
$*$-subalgebra of $\L(H_X)$ which contains $C[X]$ as an ideal.
From preceding lemma, we get then a $*$-homomorphism
from the algebra of  pseudo-local operator of $H_X$ of finite
propagation to the multiplier of $C^*_\text{max}(X)$ whose restriction
to $C[X]$ is just the inclusion $*$-homomorphism $C[X]\hookrightarrow
C^*_\text{max}(X)$.
\end{remark}
\begin{corollary}\
\begin{itemize}
\item If  $(\rho_X,H_X,T)$ is a  $K$-cycle   for $K_*(X)$  with $T$
of finite propagation. Then  $(T,C^*_\text{max}(X))$ is a 
$K$-cycle for $K_*(C^*_\text{max}(X))=KK(\C,C^*_\text{max}(X))$  we shall denote 
 by $\Ind_{X,\text{max}} T$.
\item  
$(\rho_X,H_X,T)\mapsto \Ind_{X,\text{max}} T$ gives rise to a
homomorphism
$$\widehat{\mu}_{X,\text{max},*}:K_*(X)\to K_*(C^*_\text{max}(X)).$$
\end{itemize}
\end{corollary}
\begin{proof}
We only have to check that the definition of $\Ind_{X,\text{max}} T$
only  depends on the class of $(\rho_X,H_X,T)$ in $K_*(X)$. But if
$(\rho_X,H_X,T)$ and $(\rho_X,H_X,T')$ are two K-cycle for $K_*(X)$
with $T$ and $T'$ of finite propagation then
\begin{itemize}
\item if $f(T-T')$ is compact for all $f$ in $C_0(X)$, then
  $\Ind_{\text{max}} T=\Ind_{\text{max}} T'$;
\item if $T$ and $T'$ are connected by a homotopy of operators $(T_s)_{s\in
  [0,1]}$ such that $(\rho_X,H_X,T_s)$ is a K-cycle for all $s$ in
  $[0,1]$, then according to the preceding point, we can replace $(T_s)_{s\in
  [0,1]}$ by $(\sum_{i\in\N}f_i^{1/2}T_sf_i^{1/2})_{s\in
  [0,1]}$, where $(f_i)_{i\in\N}$ is a partition of unit with support
of uniformally bounded diameter. 
\end{itemize}
The result is then a consequence of the
second item of lemma \ref{lem-pseudo-mult}
 \end{proof}

\begin{remark}
 Let $x$ be an element of $K_0(X))$ represented by an even $K$-cycle
 $(\rho_X,H_X,T)$  as in the previous corollary. Let us set
$$W=\begin{pmatrix}\Id_{H_X}&T\\0&\Id_{H_X}\end{pmatrix}\begin{pmatrix}\Id_{H_X}&0\\-T&\Id_{H_X}\end{pmatrix}\begin{pmatrix}\Id_{H_X}&T\\0&\Id_{H_X}\end{pmatrix}\begin{pmatrix}0&-\Id_{H_X}\\
  \Id_{H_X}&0\end{pmatrix}.$$ 
Then $$\left[W
\begin{pmatrix}\Id_{H_X}&0\\0&0\end{pmatrix}W^{-1}\right]-\left[
\begin{pmatrix}\Id_{H_X}&0\\0&0\end{pmatrix}\right]$$ defines an element in
$K_0(C^*_\text{max}(X))$ which is precisely 
 $\widehat{\mu}_{X,*,\text{max}}(x)$.
\end{remark}


We are now in position to define the coarse Baum-Connes assembly maps.
Recall that for  a  proper metric set  $\Si$ and a real $r$, the Rips
complex of order $r$ is the set $P_r(\Si)$ of probability measures  on $\Si$ with
support of diameter less than $r$.  Recall that $P_r(\Si)$ is a locally
finite simplicial
complex that can be provided with a proper metric extending the
euclidian metric on each  simplex. Moreover, by viewing an element of
$\Sigma$ as a Dirac measure, we get an inclusion
$\Si\hookrightarrow P_r(\Si)$, which turns out to be a coarse
equivalence.
If we fix for each real $r$ a coarse equivalence
$\phi_r:P_r(\Si)\to\Si$, then the collections  of homomorphisms given by the
  compositions
$$K_*(P_r(\Si))\stackrel{\widehat{\mu}_{P_r(\Si)},*}{\lto}K_*(C^*(P_r(\Si)))\stackrel{\phi_{r,*}}{\lto}
K_*(C^*(\Si))$$ and
$$K_*(P_r(\Si))\stackrel{\widehat{\mu}_{P_r(\Si)},\text{max},*}{\lto}K_*(C^*_\text{max}(P_r(\Si)))\stackrel{\phi_{r,\text{max},*}}{\lto}
K_*(C^*_\text{max}(\Si))$$ give  rise respectivelly  to the
the Baum-connes coarse assembly map  
$$\mu_{\Si,*}:\lim_{r>0}K_*(P_r(\Si))\lto K_*(C^*(\Si))$$ and  to the maximal
Baum-Connes assembly map
$$\mu_{\Si,*,\text{max}}:\lim_{r>0} K_*(P_r(\Si))\lto K_*(C^*_\text{max}(\Si)).$$
Moreover, if 
 $z$ in $\lim_r
K_*(P_r(\Si))$ comes  from a K-cycle $(\rho_{P_r(\Si)}, H_{P_r(\Si)},
T)$ for some  $K_*(P_r(\Si))$, where  $T$ is a finite propagation
operator on the
non-degenerated standard $P_r(\Si)$-module  $ H_{P_r(\Si)}$  then
$$\mu_{\Si,*}(z)=\phi_{r,*} \Ind_{P_r(\Si)} T$$ and $$\mu_{\Si,\text{max},*}(z)=\phi_{r,\text{max},*} \Ind_{P_r(\Si),\text{max}} T.$$

\begin{remark}\label{rem-regular} Let $\lambda_\Si:C^*_\text{max}(\Si)\to C^*(\Si)$ be the
  homomorphism induced from the representation $C[X]\hookrightarrow
  B(H_\Si)$. Then $\mu_{\Si,*}=\mu_{\Si,\text{max},*}\circ
  \lambda_{\Si,*}$

\end{remark}

\section{The Baum-Connes assembly map}\label{sec-assembly-map}
We gather this section with result we will need later on concerning
the Baum-Connes assembly map and its left hand side. For a proper
$\Ga$-space $X$ and a $\Ga$-algebra $A$, we shall denote for short
$KK^\Ga_*(X,A)$ instead of $KK^\Ga_*(C_0(X),A)$. 

\subsection{Definition of the maximal assembly map}\label{sub-ass}
Let $X$ be a locally compact proper and cocompact $\Ga$-space and let
$(\rho,\E,T)$ be a K-cycle for $KK^\Ga_*(X,A)$. Up to averaging
with a cut-off function, we can assume that the operator $T$ is
$\Ga$-equivariant. Let $\E_\Ga$ be the separated completion of $C_c(X)\cdot \E$
with respect to the $A\rtm\Ga$-valued scalar product defined by
$\langle e/e'\rangle_{\E_\Ga}(\ga)=\langle \xi/\ga(\xi')\rangle_{\E}$ for
$\xi$ and $\xi'$ in $C_c(X)\cdot \E$ and  $\ga$ in $\Ga$ (recall that
the separated completion is obtained by first  divide out by the submodule of vanishing elements for
the pseudo-norm associated to the inner product and then by completion
of the
quotient with respect to the induced norm).
Up to replace $T$ by $\sum_{\ga\in\Ga}\ga(f^{1/2})T\ga(f^{1/2})$, for
$f\in C_c(X,[0,1])$ a cut-off function with respect to the action of 
$\Ga$ of $X$, the map $C_c(X)\cdot \E\to  C_c(X)\cdot\E;\, \xi\mapsto
T\xi$ extends to an adjointable operator $T_\Ga:\E_\Ga\to\E_\Ga$. Then
we can check that $(\E_\Ga,T_\Ga)$ is a K-cycle for
$KK_*(\C,A\rtm\Ga)=K_*(A\rtm\Ga)$ whose class
$\Ind_{\Ga, A,\text{max}}T$ only depends on the class of $(\rho,\E,T)$
in $KK^\Ga_*(X,A)$.
The left hand side of the maximal assembly map is then the topological
K-theory for $\Ga$  with coefficients in $A$
$$K^\text{top}_*(\Ga,A)=\lim_{r>0}KK^\Ga_*(P_r(\Ga),A)$$ and the
assembly map 
$$\mu_{\Ga, A,\text{max}}:K^\text{top}_*(\Ga,A)\longrightarrow K_*(A\rtm\Ga)$$ is
defined for an element $x$ in $K^\text{top}_*(\Ga,A)$ coming  from
the class of a K-cycle $(\rho,\E,T)$ for $KK^\Ga_*(P_r(\Ga),A)$ by
$\mu_{\Ga, A,\text{max}}(x)=\Ind_{\Ga, A,\text{max}}T$.

\subsection{Induction}\label{subsec-lhs-ind}
We recall now from \cite{oyono} the description of  induction to a group, from the action of one
of its subgroup on a $C^*$-algebra, and    the behaviour of the left-hand side of
the Baum-Connes assembly map under this transformation.

\medskip

Let $\Ga'$ be a subgroup of a discrete group $\Ga$, and let $A$ be a
$\Ga'$-$C^*$-algebra.
Define
\begin{equation*}\begin{split}  \I_{\Ga'}^\Ga A=\{f:\Ga\to A \text{
      such that } \ga'\cdot
f(\ga\ga')=f(\ga)&\text{ for all }\ga\in\Ga\,,\ga'\in\Ga'\\ &\text{ and }\ga\Ga'\mapsto \|f(\ga)\| \text{ is in } C_0(\Ga/\Ga')\}.\end{split}\end{equation*}
Then $\Ga$ acts on $\I_{\Ga'}^\Ga A$ by left translation and it is a
standard fact that the dynamical systems $(\I_{\Ga'}^\Ga  A,\Ga)$ and
$(A,\Ga')$ have equivalent covariant representations. In particular,
the $C^*$-algebras $A\rtm\Ga'$ and $\I_{\Ga'}^\Ga A\rtm \Ga$ are
Morita equivalent (the same holds for reduced crossed products).
Notice that if the action of $\Ga'$ on $A$ is indeed the restriction
of an action of $\Ga$, then 
$$\I_{\Ga'}^\Ga A\to C_0(\Ga/\Ga',A);\, f\mapsto [\ga\Ga'\mapsto
\ga\cdot f(\ga)]$$ is a $\Ga$-equivariant isomorphism, where 
$C_0(\Ga/\Ga',A)\cong C_0(\Ga/\Ga')\ts A$ is provided with the
diagonal action of $\Ga$.
In \cite{oyono} was defined an induction homomorphism
$$\I_{\Ga', A,*}^{\Ga,\text{top}}:K^\text{top}(\Ga',A)\to K^\text{top}(\Ga
,\I_{\Ga'}^\Ga A),$$ which turned out to be an isomorphism. If $\Ga'$
has finite index in $\Ga$, then $\I_{\Ga', A,*}^{\Ga,\text{top}}$ can be
described quite easily as follows.
Recall first that if $\Ga'$ is a subgroup of $\Ga$ with finite index,
then the family of inclusions $C_0(P_r(\Ga'))\hookrightarrow
C_0(P_r(\Ga))$ gives rise to an isomorphism
\begin{equation}\label{equ-subgroup}
\lim_r
KK^{\Ga'}_*(P_r(\Ga'),A)=K^\text{top}(\Ga',A)\stackrel{\cong}{\longrightarrow}\lim_r
KK^{\Ga'}_*(P_r(\Ga),A),\end{equation} and under this
identification, the assembly map is defined as before : if 
  $x$ in an element in $K^\text{top}_*(\Ga,A)$ coming  from
the class of a K-cycle $(\rho,\E,T)$ for
$KK^{\Ga'}_*(C_0(P_r(\Ga),A)$, then 
$\mu_{\Ga, A,\text{max},*}(x)=\Ind_{\Ga', A,\text{max}}T$.

Now let $(\rho,\E,T)$ be a K-cycle for some $KK^{\Ga'}_*(X,A)$,
where $X$ is a  proper and cocompact $\Ga'$-space. Let us
define 
 $$\I_{\Ga'}^\Ga \E=\{\xi:\Ga\to A,\, \ga'\cdot
\xi(\ga\ga')=\xi(\ga)\text{ for all }\ga\in\Ga\}.$$ The
pointwise  right $A$-Hilbert module structure provides a
right $\I_{\Ga'}^\Ga A$-Hilbert module structure for $\I_{\Ga'}^\Ga
\E$ which is covariant for the action of $\Ga$ by  left
translations. If $T$ is chosen $\Ga'$-equivariant, then $\Ga\to\E;\,\ga\mapsto
T.\xi$ lies in $\I_{\Ga'}^\Ga \E$ for $\xi$ in $\I_{\Ga'}^\Ga \E$ and
 we get in this way a $\Ga$-equivariant and adjointable operator 
$\I_{\Ga'}^\Ga T:\I_{\Ga'}^\Ga \E\to\I_{\Ga'}^\Ga \E$. Finally, for $f$ in
$C_0(X)$, pointwise left multiplication by $\ga\mapsto \ga(\rho(f))$
defines a covariant  representation $\I_{\Ga'}^\Ga\rho$ of $C_0(X)$ on the
Hilbert right $\I_{\Ga'}^\Ga A$-Hilbert module  $\I_{\Ga'}^\Ga
\E$. It is straightforward to check that $(\I_{\Ga'}^\Ga\rho,\I_{\Ga'}^\Ga
\E,\I_{\Ga'}^\Ga T)$ 
is a K-cycle for $KK^\Ga(X,\I_{\Ga'}^\Ga A)$ whose class only
depends on the class of $(\rho,\E,T)$ in  $KK^{\Ga'}_*(X,A)$.
On the other hand, if $(\rho',\E',T')$ is a K-cycle for
$KK^\Ga(X,\I_{\Ga'}^\Ga A)$, and  if we consider the
$\Ga'$-equivariant homomorphism 
$\psi:\I_{\Ga'}^\Ga A\to A;\, f\mapsto f(e)$, then 
$\E=\E'\ts_\psi A$ is a $\Ga'$-covariant  right $A$-Hilbert module. Let us
set $T=T'\ts_\psi\Id_A$ and $\rho(f)=\rho'(f)\ts_\psi\Id_A$ for all
$f$ in $C_0(X)$. Then  $(\rho,\E,T)$ is a $K$-cycle for
$KK^{\Ga'}_*(X,A)$ and we can check that  $(\I_{\Ga'}^\Ga\rho,\I_{\Ga'}^\Ga
\E,\I_{\Ga'}^\Ga T)$ is a K-cycle unitary equivalent to
$(\rho',\xi',T')$. Finally, we get an isomorphism
$$\I_{\Ga',A,*}^{\Ga,X,*}:KK^{\Ga'}_*(X,A)\stackrel{\cong}{\longrightarrow}KK^\Ga_*(X,\I_{\Ga'}^\Ga
A)$$ which maps the class of a K-cycle $(\rho,\E,T)$ for
$KK^{\Ga'}_*(X,A)$  to the class of  the K-cycle $(\I_{\Ga'}^\Ga\rho,\I_{\Ga'}^\Ga
\E,\I_{\Ga'}^\Ga T)$ in $KK^\Ga_*(X,\I_{\Ga'}^\Ga
A)$.
Under the identification of equation \ref{equ-subgroup}, the family of 
isomorphisms  $(\I^{\Ga,\pr}_{\Ga',A,*})_{r>0}$ gives rise to an 
isomorphism $$\I_{\Ga',A,*}^{\Ga,\text{top}}:K^\text{top}_*(\Ga',
A)\stackrel{\cong}{\longrightarrow}K^\text{top}_*(\Ga,\I_{\Ga'}^\Ga
A).$$ Moreover, up to the identification $K_*(A\rtm\Ga')\cong
K_*(\I_{\Ga'}^\Ga A\rtm\Ga)$ induced by the Morita equivalence, we have $$\mu_{\Ga',A,\text{max}}=\mu_{\Ga,\I_{\Ga'}^\Ga
A,\Ga,\text{max}}\circ \I^{\Ga,\text{top}}_{\Ga',A,*}.$$
\begin{remark}
Let $A$ be a $\Ga$-algebra, let $\Ga'$ be a subgroup of $\Ga$ with
finite index and  let $(\rho,\E,T)$ be a  K-cycle for
$KK^{\Ga'}_*(X,A)$, such that the action of $\Ga'$ on $\E$ is
indeed the restriction of a covariant action of $\Ga$. Under the
identification $\I_{\Ga'}^\Ga A\cong  C(\Ga/\Ga',A)$ we have seen
before, then 
 
$$\I_{\Ga'}^\Ga \E \to C(\Ga/\Ga',\E);\, \xi\mapsto [\ga\Ga'\mapsto
\ga\cdot \xi(\ga)]$$ is a $\Ga$-equivariant isomorphism of right
$C(\Ga/\Ga',A)$-Hilbert module, where we have equipped
$C(\Ga/\Ga',\E)\cong C(\Ga/\Ga')\ts \E$   with the
diagonal action of $\Ga$. Moreover, under this identification,
$\I_{\Ga'}^\Ga\rho$ is given pointwise by the representation $\rho$
and $\I_{\Ga'}^\Ga T$ is the pointwise multiplication by $\ga\mapsto
\ga(T)$.
\end{remark}

\subsection{The left hand side for product of stable algebras}

As it was proved in \cite{ce}, the topological
$K$-theory for a group is a functor with respect to the coefficients which
 commutes with direct sums, i.e $K^\text{top}(G,\oplus_{i\in I}
A_i)=\oplus_{i\in I}K^\text{top}(G,A_i)$ for every locally compact group
$G$ and  every family $(A_i)_{i\in I}$ of $C^*$-algebras $A_i$ equipped
with an action of $G$ by automorphisms.
The aim of this section is to prove a similar result for product of a
family of stable  $C^*$-algebras.

\medskip
Let us first prove the result for usual K-theory.
\begin{lemma}\label{lem-prod}
Let $\A=(A_i)_{i\in I}$ be a family of  unital $C^*$-algebras. Let
$$\Theta_*^\A:K_*(\Pi_{i\in I}(A_i\ts \K(H))\lto \prod_{i\in I}K_*(A_i\ts
\K(H))\cong \prod_{i\in I}K_*(A_i)$$ be the homomorphism  induced  on the
$j$-th factor  by the   projection 
$$\prod_{i\in I}(A_i\ts \K(H))\lto A_j\ts \K(H).$$ Then $\Theta_*^\A$ is an
isomorphism.
\end{lemma}
\begin{proof}
Is clear that $\Theta^\A_*$ is onto. The injectivity of $\Theta^\A_*$ is
then a consequence of the next lemma.
\end{proof}
\begin{lemma}
There exists a map $\phi:(0,+\infty[\to(0,+\infty[$ such that
for any unital $C^*$-algebra  $A$, the following properties hold:
\begin{enumerate}
\item If $p$ and $q$ are projectors in some  $M_n(A)$ connected by a
  homotopy of projectors. Then there exists integers $k$ and $N$ with
   $n+k\leq N$ and a homotopy of projectors $(p_t)_{t\in[0,1]}$  in $M_N(A)$ connecting
  $\diag(p,I_k,0)$ and $\diag(q,I_k,0)$  and such that for any
 positive real  $\eps$ and any $s$ and $t$ in $[0,1]$ with
 $|s-t|\leq\phi(\eps)$, then $\|p_s-p_t\|\leq \eps$.
\item If $u$ and $v$ are homotopic unitaries in $U_n(A)$, then there
  exists an integer $k$ and a homotopy $(w_t)_{t\in[0,1]}$  in
  $U_{n+k}(A)$  connecting  $\diag(u,I_k)$ and $\diag(v,I_k)$
 such that for any
 positive real  $\eps$ and any $s$ and $t$ in $[0,1]$ with
 $|s-t|\leq\phi(\eps)$, then $\|w_s-w_t\|\leq \eps$.
\end{enumerate}

\end{lemma}

\begin{proof}
We can assume without loss of generality that $n=1$.
\begin{itemize}
\item Let us notice first that using \cite[Proposition 5.2.6, page 90
  ]{we}, then  there exists a positive real $\alpha$ such that for any
  unital
  $C^*$-algebra $A$ and any projectors $p$ and $q$ in $A$ such that
  $\|p-q\|\leq\alpha$, then $q=u\cdot p\cdot u^*$ for some unitary $u$
  of $A$ with $\|u-1\|\leq 1/2$. Hence there is a self-adjoint element
  $h$ of $A$ with $\|h\|\leq \ln 2$ such that $u=\exp ih$. Considering
  the homotopy of projectors $(\exp \imath th\cdot p\cdot\exp
  - \imath th)_{t\in[0,1]}$, we see that there exists a map
    $\phi_1:(0,+\infty[\to(0,+\infty[$ such that for any
  $C^*$-algebra $A$ and any projectors $p$ and $q$ in $A$ such that
  $\|p-q\|\leq\alpha$, then $p$ and $q$ are connected by a homotopy of
  projectors $(p_t)_{t\in[0,1]}$ and such that for any
 positive real  $\eps$ and any $s$ and $t$ in $[0,1]$ with
 $|s-t|\leq\phi_1(\eps)$, then $\|p_s-p_t\|\leq \eps$.
\item By considering for a projector $p$ in $A$  the homotopy of projectors
  $$\begin{pmatrix}\cos ^2 \pi t/2\cdot p&\sin  \pi t/2\cos  \pi t/2\cdot p\\
\sin  \pi t/2\cos  \pi t/2\cdot p&\sin ^2 \pi t/2\cdot p+1-p\end{pmatrix}_{t\in[0,1]}
$$ in $M_2(A)$, we also get  that there exists a map
    $\phi_2:(0,+\infty[\to(0,+\infty[$ such that for any
  $C^*$-algebra $A$ and any projector $p$  in $A$, then $\diag(1,0)$
  and $\diag(p,1-p)$ are connected by a homotopy of
  projectors $(q_t)_{t\in[0,1]}$ such that for any
 positive real  $\eps$ and any $s$ and $t$ in $[0,1]$ with
 $|s-t|\leq\phi_2(\eps)$, then $\|q_s-q_t\|\leq \eps$.
\item To prove the general case, let $p$ and $q$ be two homotopic
  projectors in a $C^*$-algebra $A$, and let $p=p_0\,,p_1,\cdots,p_m=q$
  be $m+1$ projectors in $A$ such that $\|p_{i+1}-p_i\|\leq\alpha$ for
  $i=0,\cdots,m-1$. Let us consider the following projectors in
  $M_{2m+1}(A)$:
\begin{eqnarray*}
q_0&=&\diag(p_0,I_{2k-1},0)\\
q_1&=&\diag(p_0,1,0,\cdots,1,0)\\
q_2&=&\diag(p_0,1-p_1,p_1,\cdots,1-p_m,p_m)\\
q_3&=&\diag(p_0,1-p_0,p_1,1-p_1,\cdots,p_{m-1},1-p_{m-1},p_m)\\
q_4&=&\diag(1,0,\cdots,1,0,p_m)\\
q_5&=&\diag(0,I_{2k-1},p_m)\\
q_6&=&\diag(p_m,I_{2k-1},0)
\end{eqnarray*}
Since $\|q_{3}-q_2\|\leq\alpha$, if we set $\phi=\min\{\phi_1,\phi_2\}$
and if we use the previous cases, we get for every $l$ in $\{0,5\}$
homotopies $(q^l_t)_{t\in[l,l+1]}$ between $q_l$ and $q_{l+1}$
such that  for any
 positive real  $\eps$ and any $s$ and $t$ in $[0,1]$ with
 $|s-t|\leq\phi(\eps)$, then $\|q^l_s-q^l_t\|\leq \eps$. Hence, by
 considering the total homotopy, we get the result.
\end{itemize}
The proof for unitaries is similar.
\end{proof}
\begin{proposition}\label{prop-prod}
Let $\Ga$ be a discrete group. Let $\A=(A_i)_{i\in\N}$ be a family of
$C^*$-algebras equipped with an action of $\Ga$ by automorphisms. Let
us equip $A_i\ts\K(H)$ with the diagonal action, the action of $\Ga$
on $\K(H)$ being trivial and let us then consider the induced action on
$\prod_{i\in I}(A_i\ts\K(H))$. Let
\begin{equation*}\begin{split}\Theta^{\Ga,\A}_*:KK^{\Ga}_*(P_{r}(\Ga),\Pi_{i\in I}&(A_i\ts \K(H)))\lto\\
 &\prod_{i\in I}KK^{\Ga}_*(P_{r}(\Ga),A_i\ts \K(H))\cong\prod_{i\in
   I}KK^{\Ga}_*(P_{r}(\Ga),A_i)\end{split}\end{equation*}  be the homomorphism  induced  on
 the $k$-th factor by  the projection
$$\prod_{i\in I}(A_i\ts \K(H))\to A_k\ts \K(H).$$
Then $\Theta^{\Ga,\A}_*$ is an isomorphism.
\end{proposition}
\begin{proof} Let us set $B_i=A_i\ts \K(H)$ for $i$ in $I$.
We can define an analogous morphism
 $$\Theta^X_*:KK^{\Ga}_*(X,\Pi_{i\in I}B_i)\to
 \prod_{i\geq n}KK^{\Ga}_*(X,B_i)$$ for any locally compact space $X$ equipped
 with an action of $\Ga$. Let us denote by
 $\Theta^X_{*,k}:KK^{\Ga}_*(X,\Pi_{i\in I} B_i)\to
 KK^{\Ga}_*(X,B_k)$ the homomorphism  induced  by  the projection on  the $k$-th factor. Up to
 take a barycentric subdivision of $P_{r}(\Ga)$, we can assume that
 $P_{r}(\Ga)$ is a locally finite and finite dimension typed
 simplicial complex, equipped with a simplicial and type preserving
 action of $\Ga$. Let $Z_0,\cdots,Z_n$ be the skeleton decomposition
 of $P_r(\Ga)$. Then $Z_j$ is a simplicial complex of dimension $j$,
 locally finite and equipped with a proper, cocompact  and type
 preserving simplicial action of $\Ga$. Let us prove by induction on
 $i$ that
$\Theta^{Z_j}_*$ is an isomorphism. The $0$-skeletton $Z_0$ is a
finite union of orbits and thus, for $j=0$, it is enought to prove that
$\Theta^{\Ga/F}_*$ is an isomorphism when $F$ is a finite subgroup of
$\Ga$. Let us recall from \cite{oyono}  that for every
$C^*$-algebra $B$ equipped with an action of $\Ga$, there is a natural
restriction isomorphism $\res^B_{F,\Ga}:KK^\Ga_*(\Ga/F,B)\lto
KK^F_*(\C,B)\cong K_*(B\rtimes F)$. We get by naturality  the following
commutative diagram
$$
\begin{CD}
 KK^\Ga_*(\Ga/F,\Pi_{i\in I}B_i)  @>\Theta^{\Ga/F}_{*,k}>>KK^{\Ga}_*(\Ga/F,B_k)\\
        @V\res^{\Pi_{i\in I}B_i}_{F,\Ga}VV            @VV\res^{B_k}_{F,\Ga} V\\
 K_*((\Pi_{i\in I}B_i)\rtimes F)   @>>> K_*(B_k\rtimes F)
\end{CD},
$$
where the bottom row is induced by the homomorphism  $\Pi_{i\in I}(B_i\rtimes
F)\to B_k\rtimes F$ arising from the projection on the
$k$-th factor $\Pi_{i\in I}B_i\to B_k$. Since $F$ is finite, $\Pi_{i\in I}(B_i\rtimes
F)$ is naturally isomorphic to  $(\Pi_{i\in I}B_i)\rtimes F$,
and under this identification, the bottom row homomorphism induces by
lemma \ref{lem-prod} an isomorphism  $$K_*((\Pi_{i\in I}B_i)\rtimes F)\lto \prod_{i\in I}K_*(B_i\rtimes
F).$$ Hence $\Theta^{\Ga/F}_{*}$ is an isomorphism.

Let us assume that we have proved that  $\Theta^{Z_{j-1}}_{*}$ is an
isomorphism. Then the short exact sequence 
$$0\lto C_0(Z_j\setminus Z_{j-1}) \lto C_0(Z_j)\lto C_0(Z_{j-1})\lto 0  $$
gives rise to an natural long exact sequence
$$\lto KK^\Ga_*(Z_{j-1},\bullet)\lto KK^\Ga_*(Z_{j},\bullet)\lto
KK^\Ga_*(Z_{j}\setminus Z_{j-1},\bullet)\lto
KK^\Ga_{*+1}(Z_{j-1},\bullet)$$ and thus by naturallity, we get a
diagram
{\tiny{$$
\begin{CD}
 KK^\Ga_*(Z_{j-1},\Pi_{i\in I}B_i) @>>> KK^\Ga_*(Z_{j},\Pi_{i\in
   I}B_i) @>>> KK^\Ga_*(Z_{j}\setminus Z_{j-1},\Pi_{i\in I}B_i)\\
        @V\Theta^{Z_{j-1}}_{*}VV    
        @V\Theta^{Z_{j}}_{*}VV   @V\Theta^{Z_{j-1}}_{*+1}VV\\
{\displaystyle \Pi_{i\in I}} KK^\Ga_*(Z_{j-1},B_i) @>>>\Pi_{i\in I} KK^\Ga_*(Z_{j},B_i) @>>>\Pi_{i\in I} KK^\Ga_*(Z_{j},B_i)\\
&&& 
@>>>
KK^\Ga_{*+1}(Z_{j-1},\Pi_{i\in I}B_i)\\
&&&&&  @V\Theta^{Z_{j}\setminus Z_{j-1}}_{*}VV \\
&&&  @>>>
\Pi_{i\in I}KK^\Ga_{*+1}(Z_{j-1},B_i)
\end{CD},
$$}}
Let $\intsi_j$  be the interior of the standard $j$-simplex. Since the
action of $\Ga$ is type preserving, then $Z_{j}\setminus Z_{j-1}$ is
equivariantly homeomorphic  to $\intsi_j\times\Sigma_j$, where $\Sigma_j$ is the
set of center of $j$-simplices  of $Z_{j}$, and where $\Ga$ acts
trivially on $\intsi_j$. This identification, together with Bott
periodicity provides a commutative diagram

$$
\begin{CD}
KK^\Ga_*(Z_{j}\setminus Z_{j-1},\Pi_{i\in I}B_i) @>>>KK^\Ga_{*+1}(\Sigma_j,\Pi_{i\in I}B_i)\\
  @V\Theta^{Z_{j}\setminus Z_{j-1}}_{*}VV   
        @V\Theta^{\Sigma_j}_{*+1}VV\\
 \prod_{i\in I} KK^\Ga_*(Z_{j}\setminus Z_{j-1},B_i) @>>>
\prod_{i\in I}KK^\Ga_{*+1}(\Sigma_j,B_i) \end{CD},
$$

By the first step of induction, $\Theta^{\Sigma_j}_{*}$ is an isomorphism,
and hence  $\Theta^{Z_{j}\setminus Z_{j-1}}_{*}$ is an
isomorphism. Using the induction hypothesis and the five lemma, we get
then that $\Theta^{Z_{j}}_{*}$ is an isomorphism.\end{proof}
Let $(A_i)_{i\in \N}$ be a family of $\Gamma$-algebras, let $H$ be an
Hilbert space and let $x$ be an
element of $KK^*(\pr,\prod_{i\in\N} A_i\ts\K(H))$ (the action of
$\Ga$ on $H$ being trivial) represented by a
K-cycle $(\phi,\E,T)$. Let $p_k:\prod_{i\in\N} A_i\ts\K(H)\to
A_k\ts\K(H)$ be the canonical projection on the $k$-th factor, and let
us set $\E_k=\E\ts_{p_k} A_k\ts\K(H)$, $T_k=T\ts_{p_k}\Id_{A_k\ts\K(H)}$ and let us
define the $\Ga$-equivariant representation of $C_0(\pr)$ on $\E_k$ by
$\phi_k(f)=\phi(f)\ts_{p_k}\Id_{A_k\ts\K(H)}$ for all $f$ in
$C_0(\pr)$.
Then $\prod_{i\in\N}\E_i$ provided with the diagonal action 
is a $\Ga$-equivariant right $\prod_{i\in\N} A_i\ts\K(H)$-Hilbert
module. Moreover, if $S$ is a compact operator on $\E$, then for every
$\eps>0$, there exists a finite rank operator $S'$ on $\E$ such that
$\|S-S'\|\leq\eps$. Then $(S_i')_{i\in\N}$ provides a finite rank
operator on $\prod_{i\in\N}\E_i$ such that $\|S_i-S'_i\|\leq\eps$ for
all integer $i$. Hence  $(S_i)_{i\in\N}$ gives rise to a compact 
operator on $\prod_{i\in\N}\E_i$. Consequently,
$((\phi_i)_{i\in\N},\prod_{i\in\N}\E_i,(T_i)_{i\in\N})$ is a
K-cycle for $KK_*^\Ga(\pr,\prod_{i\in\N} A_i\ts\K(H))$ which  in
view of the isomorphism of proposition \ref{prop-prod} represents also
$x$. Using the imprimitivity bimodule implementing the Morita
equivalence between $A_i$ and $A_i\ts\K(H)$, we can actually replace
 $\E_i$ by $\K(A_i\otimes H,H_i)$ where
 $H_i=\E_i\otimes_{A_i\ts\K(H)}\otimes A_i$. Hence we 
 obtain that 
every element $x$ in  $KK^*(\pr,\prod_{i\in\N} A_i\ts\K(H))$ can be
represented by a K-cycle
$((\phi_i)_{i\in\N},\prod_{i\in\N}\K_{A_i}(A_i\ts
H,\mathcal{H}_i),(T_i)_{i\in\N})$ such that for every integer $i$,
\begin{itemize}
\item $\mathcal{H}_i$ is a $\Ga$-equivariant right $A_i$-Hilbert
  module;
\item  $\phi_i$ is a $\Ga$-equivariant representation of $C_0(\pr)$ on
  $\mathcal{H}_i$;
\item  $T_i$ is a $\Ga$-equivariant operator  on 
  $\mathcal{H}_i$;
\item the action of $T_i$ and of $\phi_i(f)$ for $f$ in $C_0(\pr)$ on 
$\K_{A_i}(A_i\ts
H,\mathcal{H}_i)$ being by left composition.
\end{itemize}
 Moreover, we can assume that $\|T_i\|\leq 1$ for all positive integer
 $i$.

\smallskip

As a consequence of proposition \ref{prop-prod} we get
\begin{corollary}\label{cor-ass-prod}
If $\Ga$ admits a universal example which is a finite dimension and
cocompact simplicial complex (equipped with a simplicial action of
$\Ga$), then we have an isomorphism
 $$K^{\text{top}}_*(\Ga,\Pi_{i\in I}(A_i\ts \K(H)))\to
 \prod_{i\in I}K^{\text{top}}_*(\Ga,A_i\ts \K(H))\cong\prod_{i\in I}K^{\text{top}}_*(\Ga,A_i)$$  induced  on
 the $k$-th factor by  the projection
$$\Pi_{i\in I}(A_i\ts \K(H))\to A_k\ts \K(H).$$
\end{corollary}

\subsection{The case of coverings}\label{subsec-lhs-cov}
Recall from \cite{val} that for a cocompact covering $\widetilde{X}\to
X$ of group $\Ga$, we have a natural isomorphism
$${\Upsilon}_{\widetilde{X},*}^\Ga:KK^{\Ga}_*(\widetilde{X},\C){\lto}K_*(X)$$ which can be described as follows.
 Let $(\rho,H,T)$ be a $K$-cycle for $KK^{\Ga}_*(\widetilde{X},\C)$. We
can assume without loss of generality that the representation
$\rho:C_0(\widetilde{X})\to \L(H)$ is non-degenerated. We can also assume that the operator $T$ of
$\L(H)$ is $\Ga$-equivariant and that 
$$T\cdot C_c(\widetilde{X})\cdot H\subset C_c(\widetilde{X})\cdot H.$$ If $\langle\bullet,\bullet\rangle$ is the
scalar product on the Hilbert space $H$, then we can define on
$C_c(\widetilde{X})\cdot H$ the inner product
$$\langle\langle\xi,\eta\rangle\rangle=\sum_{\ga\in\Ga}\langle\xi,\ga(\eta)\rangle.$$
Then $\langle\langle\bullet,\bullet\rangle\rangle$ is positive and
thus by taking the separated completion of $C_c(\widetilde{X})\cdot
H$, we get a Hilbert space $\widehat{H}$.
The operator $T$ being equivariant, its restriction to
$C_c(\widetilde{X})\cdot H$ extends to a continuous operator
$\widehat{T}$ on
$\widehat{H}$. Since $\rho$ is non-degenerated, it extends to a
representation of $C({X})$  (viewed  as an algebra of
multiplier for $C_0(\widetilde{X})$)  on $H$ by equivariant
operator. Moreover, since this representation  preserves
$C_c(\widetilde{X})\cdot H$, it induces a representation
$\widehat{\rho}$ of $C({X})$  on $\widehat{H}$.
It is straightforward to check that 
$(\widehat{\rho},\widehat{H},\widehat{T})$ is a $K$-cycle for
$K_*({X})$ and  we get
is this way a homomorphism
\begin{equation}\label{equ-iso-valette}
{\Upsilon}^\Ga_{\widetilde{X},*}:KK^{\Ga}_*(C_0(\widetilde{X}),\C){\lto}K_*({X})
\end{equation} which maps the class of $({\rho},{H},{T})$ to the  class
of
$(\widehat{\rho},\widehat{H},\widehat{T})$.

\begin{theorem}\cite{val}
${\Upsilon}^\Ga_{\widetilde{X},*}$ is a isomorphism.
\end{theorem}

We want  now to study how the  propagation behave under the above transformation.
So assume that $\widetilde{X}$ is a locally compact metric space equipped with 
 a free, proper, isometric and cocompact action of $\Ga$. Let
$\eta$ be a $\Ga$-invariant measure, and let $\widehat{\eta}$ be the
measure induced on $X=\widetilde{X}/\Ga$. Let us set $H_{\widetilde{X}}=L^2(\eta)\ts
H$ and $H_{{X}}=L^2(\widehat{\eta})\ts
H$. We can view $C({X})$ as the algebra of
$\Ga$-invariant continuous and bounded functions on $\widetilde{X}$ and according
to this, for any continuous and compactly supported function $f:\widetilde{X}\to
\C$, then $\widehat{f}=\sum_{\ga\in\Ga}\ga(f)$ belongs to
$H_{\widehat{X}}$.
It is straightforward to check that $f\mapsto\widehat{f}$ extends to
a unitary map $\widehat{H_{\widetilde{X}}}\to H_X $.
\begin{lemma}\label{lem-cov-lc}
If $T$ is a locally compact equivariant operator on 
$H_{\widetilde{X}}$ with propagation less than $r$. Then, under the
above identification
between $\widehat{H_{\widetilde{X}}}$ and $H_X$, the operator
$\widehat{T}$
is a compact operator with propagation less than $r$.
\end{lemma}
\begin{proof}
Since $T$ is equivariant and since $\widetilde{X}$ is cocompact, the operator $T$
is given by a kernel $K:\widetilde{X}\times \widetilde{X}\to \K(H)$ such that $K(\ga x,\ga
y)=K(x,y)$ for almost all $(x,y)$ in $\widetilde{X}\times\widetilde{X}$ and with cocompact
support (for the diagonal action of $\Ga$ on $\widetilde{X}\times\widetilde{X}$) of diameter
less than $r$. 
Under the above identification between $\widehat{H_{\widetilde{X}}}$
and $H_X$, then for any 
  continuous and compactly supported function $f:\widetilde{X}\to
\C$, we have $$\widehat{T}\cdot \widehat{f}=\sum_{\ga\in\Ga}\ga(T\cdot
f)=\sum_{\ga\in\Ga}T\cdot\ga( f).$$
By viewing $X$ as a  borelian fundamental domain for the action
of  $\Ga$ on $\widetilde{X}$, we get
\begin{eqnarray*}
\widehat{T}\cdot \widehat{f}(x)&=&\sum_{\ga\in\Ga}\int_{\widetilde{X}}
K(x,y)f(\ga y)d\eta(y)\\
&=&\sum_{(\ga,\ga')\in\Ga^2}\int_{X}
K( x,\ga'y)f(\ga\ga' y)d\widehat{\eta}(y)\\
&=&\sum_{(\ga,\ga')\in\Ga^2}\int_{X}
K(\ga'^{-1} x,y)f(\ga y)d\widehat{\eta}(y)\\
&=&\int_{X}\sum_{\ga'\in\Ga}K(\ga'^{-1}
x,y)\widehat{f}(y)d\widehat{\eta}(y)\\
&=&\int_{X}F(x,y)
\widehat{f}(y)d\eta(y),\\
\end{eqnarray*}
with $$F:\widetilde{X}\times \widetilde{X};(x,y)\mapsto \sum_{\ga\in\Ga}K(\ga
x,y).$$
The kernel $F$ is $\Ga\times\Ga$-invariant and thus can be viewed as a
kernel on ${X}\times{X}$ and thus we get
$\widehat{T}\cdot \widehat{f}(x)=\int_{{X}}F(x,y)\widehat{f}(y)d\widehat{\eta}(y)$. Hence $\widehat{T}$ is a compact
operator and  since $F(x,y)=0$ for
almost every $(x,y)$ in ${X}\times {X}$ such that $d(x, y)\geq r$,
we see that $\widehat{T}$ as propagation less than $r$.
\end{proof}
The previous lemma can be extended to pseudo-local operators on
$H_{\widetilde{X}}$ with finite propagation. Recall from \cite{val}
  that if $T$ is a pseudo-local equivariant operator on $H_{\widetilde{X}}$ with
  finite propagation, then $\widehat{T}$ is a pseudo-local  operator on
  $\widehat{H_{\widetilde{X}}}\cong H_X$
\begin{lemma}\label{lem-cov-psl}
With notation of lemma \ref{lem-cov-lc}, if 
 $T$ is pseudo-local  $\Ga$-equivariant  operator  on 
$H_{\widetilde{X}}$ with propagation less than $r$. Then, the operator
$\widehat{T}$
is a pseudo-local  operator with propagation less than $r$
\end{lemma}
\begin{proof}
Let $f_1,\ldots,f_n$ be a partition of unit for $X$ with support of
diameter less than $r$. Let us set $T'=\sum_{i=1}^n f_i^{1/2}\circ
q \cdot T \cdot f_i^{1/2}\circ
q$, where $q:\tilde{X}\to X$ is the projection map of the
covering. Then $T'-T$ is an equivariant  locally compact on
$H_{\widetilde{X}}$ with propagation less than $r$ and  thus,
according to lemma \ref{lem-cov-lc}, we get that
$\widehat{T}-\widehat{T'}$ is compact and has propagation less than
$r$.
Since $\widehat{T'}=\sum_{i=1}^n f_i^{1/2} \cdot \widehat{T}\cdot  f_i^{1/2}$,
then $T'$ is a pseudo-local operator of propagation less than $r$ and
hence we get the result.
\end{proof}


\section{The left hand side for the coarse space associated  to a residually
  finite group}
 \label{sec-khom}
The aim of
this section is to state for the sources of the assembly maps the
analogous  of proposition \ref{prop-roemax}, i.e the existence of  a
group homomorphism $$\Psi_{\x,*}:\lim_r
KK^\Ga(P_r(\x),\C)\lto K^\text{top}(\Ga,A_\Ga),$$ such that 
\begin{equation}\label{equ-comp}\Psi_{\Ga,\a,*}\circ\mu_{\x,\text{max},*}=\mu_{\Ga,A_\Ga,\text{max},*}\circ
\Psi_{\x,*}.\end{equation}

\subsection{Rips complexes  associated to a residually
  finite group} \label{subsec-rips}
Let $\Ga$ be a residually finite group, finitely generated. Let
$\Ga_0\supset\Ga_1\supset\ldots\Ga_n\supset\ldots$ be a decreasing
sequence of normal finite index subgroups of $\Gamma$ such that
$\bigcap_{i\in\N}\Ga_i=\{e\}$. Recall  notations of  section
\ref{sub-roe-resfin}, if  $d$ be a left invariant metric
associated to any finite set of generators for $\Ga$, then we  endow $\Ga/\Ga_i$
with the metric $d(a\Ga_i,b\Ga_i)=\min\{d(a\ga_1,b\ga_2), \ga_1\text{
  and }\ga_2 \text{
  in } \Ga_i\}$.
We set $\x=\du_{i\in\N}\Ga/\Ga_i$ and we equip $\x$ with a metric $d$
such that on  $\Ga/\Ga_i$, then $d$ is the metric defined above and 
$d(\Ga/\Ga_i,\Ga/\Ga_i)\geq i+j$ if $i\neq j$.

\medskip

For every integer $n$ such that $r>n$, then
$$P_r(\x)=P_r(\coprod_{i=1}^{n-1} \Ga/\Ga_i)\coprod\left(
\,\coprod_{i\geq n}P_r(\Ga/\Ga_i)\right),$$where $P_r(\coprod_{i=1}^{n-1}
\Ga/\Ga_i)$ and $\coprod_{i\geq n}P_r(\Ga/\Ga_i)$ can be viewed as
distinct open subsets of $P_r(\x)$. Hence we have a splitting
\begin{eqnarray}\label{equ-phi1}
\nonumber K_*(P_r(\x))&\cong& K_*(P_r(\coprod_{i=1}^{n-1}
\Ga/\Ga_i))\bigoplus K_*(\coprod_{i\geq n}P_r(\Ga/\Ga_i))\\
&\cong& K_*(P_r(\coprod_{i=1}^{n-1}
\Ga/\Ga_i))\bigoplus \prod_{i\geq n}K_*(P_r(\Ga/\Ga_i))
\end{eqnarray}
corresponding to the inclusion of the disjoint open subsets 
$P_r(\coprod_{i=1}^{n-1}
\Ga/\Ga_i)$ and $\coprod_{i\geq n}P_r(\Ga/\Ga_i)$ into $P_r(\x)$.

Let us show that in the inductive limit  when $r$ runs through positive real, $K_*(P_r(\Ga/\Ga_i))$, behave
like $K_*(P_r(\Ga)/\Ga_i)$ (recall that $\Ga_i$ acts properly
on $P_r(\Ga)$).
For $f$ in $P_r(\Ga)$, let us define
$$\widetilde{f}:\Ga/\Ga_i\to[0,1];\,\ga\Ga_i\mapsto \sum_{g\in\Ga_i}f(\ga
g).$$ Then $\widetilde{f}$ is a probability on $\Ga/\Ga_i$. Let $\ga$
and $\ga'$ be elements of $\Ga$ such that $\widetilde{f}(\ga\Ga_i)\neq 0$ and
$\widetilde{f}(\ga'\Ga_i)\neq 0$. Then there exists $g$ and $g'$ in $\Ga_i$ such
that $f(\ga g)\neq 0$ and
$f(\ga'g)\neq 0$. Since $f$ is in $P_r(\Ga)$, we get that $d(\ga
g,\ga'g)\leq r$ and hence $d(\ga\Ga_i
,\ga'\Ga_i)\leq r$. Thus  $\widetilde{f}$ belongs to $P_r(\Ga/\Ga_i)$,
and since $\widetilde{\ga\cdot f}=\widetilde{f}$ for any $\ga$ in $\Ga_i$, we finally obtain a
continuous map $\upsilon_{r,i}:P_r(\Ga)/\Ga_i\to
P_r(\Ga/\Ga_i);\dot{f}\mapsto \widetilde{f}$, where $\dot{f}$ is the
class in $P_r(\Ga)/\Ga_i$ of  $f$ in $P_r(\Ga)$.
For a positive real $r$ and an integer $n$, let
$$ \Lambda_{*,r,n}:\prod_{k\geq n}K_*(P_{r}(\Ga)/\Ga_k)\lto
\prod_{k\geq n}K_*(P_r(\Ga/\Ga_k))$$
be the homomorphism  induced  on the $k$-th factor by the map 
$$P_r(\Ga)/\Ga_k\to P_{r}(\Ga/\Ga_k);\, \dot{f}\mapsto \widetilde{f}.$$
\begin{lemma}\label{lem-diam}
Let $i$ be an integer such that $B_\Ga(e,2r)\cap\Ga_i=\{e\}$. Let
$\{\ga_1,\cdots,\ga_n\}$ and $\{\ga'_1,\cdots,\ga'_n\}$ be subsets of 
$\Ga$ of diameter less than $r$ and such that $\ga_j{\ga'}_j^{-1}$ is in
$\Ga_i$ for all $j$ in $\{1,\cdots,n\}$, then
$\ga_j{\ga'}_j^{-1}=\ga_k{\ga'}_k^{-1}$ for all  $j$ and $k$ in $\{1,\cdots,n\}$.
\end{lemma}
\begin{proof}
We have $d(\ga_1,\ga_j)\leq r$ and $d(\ga'_1,\ga'_j)\leq r$ for all
$j$ in $\{1,\cdots,n\}$. Let us set $g=\ga_1{\ga'}_1^{-1}$.
Then
\begin{eqnarray*}
d(\ga_j,g\ga'_j)&\leq&d(\ga_j,\ga_1)+d(\ga_1,g\ga'_j)\\
&\leq&d(\ga_j,\ga_1)+d(g\ga'_1,g\ga'_j)\\
&\leq&d(\ga_j,\ga_1)+d(\ga'_1,\ga'_j)\\
&\leq&2r.
\end{eqnarray*}
Hence, since $\Ga_i$ is normal,
$\ga_j^{-1}g{\ga}'_j=(\ga_j^{-1}g\ga_j)\ga_j^{-1}(\ga_j'\ga_j^{-1})\ga_j$
belongs to $B_\Ga(e,2r)\cap\Ga_i$ and thus $\ga_j=g\ga'_j$.
\end{proof}
Let $i$ be an integer such that $B_\Ga(e,4r)\cap\Ga_i=\{e\}$. Let $h$
be an element of  $P_r(\Ga/\Ga_i)$. We can choose a finite subset
$\{\ga_1,\cdots,\ga_n\}$ of diameter less than $2r$  such that the
support of $h$ lies in  $\{\ga_1\Ga_i,\cdots,\ga_n\Ga_i\}$. According to
lemma \ref{lem-diam}, a such subset is unique up to left translations by an
element of $\Ga_i$. Let us define $\hat{h}$ in $P_{2r}(\Ga)/\Ga_i$ as the
class of the probability  of $P_r(\Ga)$ with support in
$\{\ga_1,\cdots,\ga_n\}$ with value  on an element $\ga$ in that set 
$\hat{h}(\ga)=h(\ga\Ga_i)$. It is straightforward to check that if $h$
is in $P_r(\Ga/\Ga_i)$, then $\widetilde{\hat{h}}$ is the image of $h$
under the inclusion map $P_r(\Ga/\Ga_i)\hookrightarrow
P_{2r}(\Ga/\Ga_i)$. If $f$ is an element of $P_r(\Ga)$, then since
$B_\Ga(e,r)\cap\Ga_i=\{e\}$, the intersection of the support of $f$
with any $\ga\Ga_i$ for $\ga$ in $\Ga$ has at most one element. Hence,
according to
lemma \ref{lem-diam}, $\widehat{\widetilde{f}}$ is the image of the class
of  $f$ in $P_{r}(\Ga)/\Ga_i$ under the inclusion
$P_{r}(\Ga)/\Ga_i\hookrightarrow P_{2r}(\Ga)/\Ga_i$.

\begin{lemma}
 Let $r$ be a positive real and let $i$ be an integer such that
 $B_\Ga(e,r)\cap\Ga_i=\{e\}$. Then the action of $\Ga_i$ on $P(\Ga)$
 is free.
\end{lemma}
\begin{proof}
Let $f$ be an element of $P_r(\Ga)$. If $\ga\cdot f=f$ with $\ga$ in
$\Ga_i$, then the support of $f$ is invariant under the action of
$\ga$. In particular, $\ga=g_1\cdot g_2^{-1}$ with $g_1$ and $g_2$ in the
support of $f$. Hence $\ga$ is in  $B_\Ga(e,r)\cap\Ga_i$ and  thus
$\ga=e$.
\end{proof}
In consequence, with condition of the lemma above, $P_r(\Ga)\to
P_r(\Ga)/\Ga_i$ is a covering map and since $\Ga_i$ has finite index
in $\Ga$, this covering is cocompact.
\subsection{Construction of $\Psi_\x$}
For a positive real $r$ and an integer $n$, such that $r\geq n$ and
$B_\Ga(e,4r)\cap\Ga_n=\{e\}$ let us define
\begin{itemize}
\item 
$\Psi_{*,r,n}^1:K_*(P_r(\x))\to\prod_{i\geq
  n}K_*(P_r(\Ga/\Ga_i))$  the projection homomorphism corresponding to the 
decomposition in equation \ref{equ-phi1} of  section \ref{subsec-rips}.
\item 
$ \Psi_{*,r,n}^2:\prod_{k\geq n}K_*(P_{r}(\Ga/\Ga_k))\lto
\prod_{k\geq n}K_*(P_{2r}(\Ga)/\Ga_k)$
 the homomorphism induced  on the $k$-th factor by the maps 
$P_r(\Ga/\Ga_k)\to P_{2r}(\Ga)/\Ga_k;\, h\mapsto \hat{h}$.
\item 
$$\Psi_{*,r,n}^3:\prod_{k\geq n}K_*(P_{r}(\Ga)/\Ga_k)\lto
\prod_{k\geq n}KK^{\Ga_k}_*(P_{r}(\Ga),\C)$$  the homomorphism given on
the $k$-th factor by the inverse of the isomorphism
$\Upsilon^{\Ga_k}_{P_{r}(\Ga),*}:KK^{\Ga_k}_*(P_{r}(\Ga),\C)\stackrel{\cong}{\lto} K_*(P_{r}(\Ga)/\Ga_k)$ (see section \ref{subsec-lhs-cov}).
\item  $\Psi_{*,r,n}^4:\prod_{k\geq n}KK^{\Ga_k}_*(P_{r}(\Ga),\C)\lto
\prod_{k\geq n}KK^{\Ga}_*(P_{r}(\Ga),C(\Ga/\Ga_k))$  the homomorphism
given on the $k$-th factor by the induction homomorphism
$\I_{\Ga_k,*}^{\Ga, P_{r}(\Ga)}$ (see section \ref{subsec-lhs-ind}).
\item $\Psi_{*,r,n}^5:\prod_{i\geq n}KK^{\Ga}_*(P_{r}(\Ga),C(\Ga/\Ga_i))\lto
KK^{\Ga}_*(P_{r}(\Ga),\ell^\infty(\sqcup_{i\geq n} \Ga/\Ga_i,\K(H)))$
the inverse of the isomorphism $\Theta_{*}^{\Ga,\A}$ of proposition
\ref{prop-prod} applied to the family $\A=(C(\Ga/\Ga_i))_{i\in\N}$;
\item  $\Psi_{*,r,n}^6:
KK^{\Ga}_*(P_{r}(\Ga),\ell^\infty(\sqcup_{i\geq n}
\Ga/\Ga_i,\K(H)))\lto KK^{\Ga}_*(P_{r}(\Ga),A_\Ga)$ the homomorphism 
induced by the $\Ga$-equivariant epimorphism 
$$\ell^\infty(\sqcup_{i\geq n} \Ga/\Ga_i,\K(H)){\lto} A_\Ga.$$ 
\end{itemize}
\begin{remark}\label{rem-restriction}
\begin{enumerate} 
\item  Let us
also define for any real $r$ and $r'$ such that $0\leq r\leq r'$ the
homomorphisms $$\iota_{*,r,r'}^{k\geq n}:\prod_{k\geq
  n}K_*(P_{r}(\Ga/\Ga_k))\lto \prod_{k\geq
  n} K_*(P_{r'}(\Ga/\Ga_k))$$
and $${\iota'}_{*,r,r'}^{k\geq n}:\prod_{k\geq
  n}K_*(P_{r}(\Ga/)\Ga_k)\lto \prod_{k\geq
  n} K_*(P_{r'}(\Ga/\Ga_k)$$
 respectively induced on the $k$-th factor by the
inclusions $P_{r}(\Ga/\Ga_k)\hookrightarrow P_{r'}(\Ga/\Ga_k)$ and 
$P_{r}(\Ga)/\Ga_k\hookrightarrow P_{r'}(\Ga)/\Ga_k$. According to
the discussion that follows lemma \ref{lem-diam} and with notations of
section \ref{subsec-rips}, if $n$ is chosen  such that 
$n\leq r$ and $B_\Ga(e,4r)\cap\Ga_n=\{e\}$, then 
\begin{eqnarray}\label{equ-psi1} \Psi_{*,r,n}^2\circ\Lambda_{*,r,n}&=&{\iota'}_{*,r,2r}^{k\geq n}
\end{eqnarray} and 
\begin{eqnarray}\label{equ-psi2}\Lambda_{*,2r,n}\circ\Psi_{*,r,n}^2&=&\iota_{*,r,2r}^{k\geq n}.
\end{eqnarray}
\item
Using the same
argument as in the proof of proposition \ref{prop-prod}, we get that 
$\Psi_{*,r,n}^5$ restricts to an isomorphism $$\bigoplus_{i\geq n}KK^{\Ga}_*(P_{r}(\Ga),C(\Ga/\Ga_i))
\stackrel{\cong}{\lto}
KK^{\Ga}_*(P_{r}(\Ga),C_0(\sqcup_{i\geq n}
\Ga/\Ga_i,\K(H))).$$
\end{enumerate}
\end{remark}
According to  the next lemma, $\Psi_{*,r,n}^6$ is an epimorphism.
\begin{proposition}\label{prop-surj}
The equivariant short exact sequence 
$$0\lto C_0(\sqcup_{i\geq n}
\Ga/\Ga_i,\K(H))\lto\ell^\infty(\sqcup_{i\geq n} \Ga/\Ga_i,\K(H))\lto
A_\Ga\lto 0$$ gives rise to a short exact sequence 
\begin{equation*}\begin{split} 0\lto KK^{\Ga}_*(P_{r}(\Ga),&C_0(\sqcup_{i\geq n}
\Ga/\Ga_i,\K(H)))\lto \\&KK^{\Ga}_*(P_{r}(\Ga),\ell^\infty(\sqcup_{i\geq
  n} \Ga/\Ga_i,\K(H)))
\lto
KK^{\Ga}_*(P_{r}(\Ga),A_\Ga)\lto 0.\end{split}\end{equation*}
\end{proposition}
\begin{proof}
Using the six-term exact sequence associated to an equivariant short
exact sequence of $C^*$-algebras, this amounts to show that the
inclusion
$$\iota:C_0(\sqcup_{i\geq n}
\Ga/\Ga_i,\K(H))\hookrightarrow \ell^\infty(\sqcup_{i\geq n}
\Ga/\Ga_i,\K(H))$$ induces a monomorphism 
$$\iota_*: KK^{\Ga}_*(P_{r}(\Ga),C_0(\sqcup_{i\geq n}
\Ga/\Ga_i,\K(H)))\hookrightarrow
KK^{\Ga}_*(P_{r}(\Ga),\ell^\infty(\sqcup_{i\geq n} \Ga/\Ga_i,\K(H))).$$
According to remark \ref{rem-restriction} and to proposition
\ref{prop-prod}, we have a splitting
$$KK^{\Ga}_*(P_{r}(\Ga),C_0(\sqcup_{i\geq n}
\Ga/\Ga_i,\K(H)))\cong\bigoplus_{i\geq n} KK^{\Ga}_*(P_{r}(\Ga),C(
\Ga/\Ga_i))$$ and an isomorphism
$$KK^{\Ga}_*(P_{r}(\Ga),\ell^\infty(\sqcup_{i\geq n}
\Ga/\Ga_i,\K(H)))\cong\prod_{i\geq n} KK^{\Ga}_*(P_{r}(\Ga),C(
\Ga/\Ga_i)).$$ Up to these identifications, 
$\iota_*$ is the inclusion
$$\bigoplus_{i\geq n} KK^{\Ga}_*(P_{r}(\Ga),C(
\Ga/\Ga_i)) \hookrightarrow\prod_{i\geq n} KK^{\Ga}_*(P_{r}(\Ga),C(
\Ga/\Ga_i)).$$
\end{proof}
\begin{remark}\label{rem-intertwin}
Taking the inductive limit over all the $\pr$, we get
a short exact sequence 
\begin{equation*}\begin{split}0\lto K^{\text{top}}_*(\Ga,C_0(\sqcup_{i\geq n}
\Ga/\Ga_i,&\K(H)))\lto \\&K^{\text{top}}_*(\Ga,\ell^\infty(\sqcup_{i\geq
  n} \Ga/\Ga_i,\K(H))) 
\lto
 K^{\text{top}}_*(\Ga,A_\Ga)\lto 0\end{split}\end{equation*}
In the same way, since the composition \begin{equation*}\begin{split} \bigoplus_{i\geq n} K_*(C(
\Ga/\Ga_i,\K(H))\rtm\Ga)\to K_*(\prod_{i\geq n}
&C(\Ga/\Ga_i,\K(H))\rtm\Ga)\\ &\to \prod_{i\geq n}
K_*(C(\Ga/\Ga_i,\K(H))\rtm\Ga)\end{split}\end{equation*}  is injective, where the second map is
induced on the $k$-th factor by the projection $\prod_{i\geq n}
C(\Ga/\Ga_i,\K(H))\to C(\Ga/\Ga_k,\K(H))$, the 
exact sequence  for maximal cross product
$$0\lto C_0(\sqcup_{i\geq n}
\Ga/\Ga_i,\K(H))\rtm\Ga\lto\ell^\infty(\sqcup_{i\geq n} \Ga/\Ga_i,\K(H))\rtm\Ga\lto
A_\Ga\rtm\Ga\lto 0$$ gives rise to a short exact sequence 
\begin{equation*}\begin{split} 0\lto \oplus_{i\geq n} K_*(C_0(
\Ga/\Ga_i,&\K(H))\rtm\Ga)\lto\\
& K_*(\ell^\infty(\sqcup_{i\geq n}
\Ga/\Ga_i,\K(H))\rtm\Ga)\lto
K_*(A_\Ga\rtm\Ga)\lto 0\end{split}\end{equation*} and moreover, the assembly maps intertwin
the corresponding above exact sequences
\end{remark}
Let $r$ be a positive real and let $n$ be an integer such that $n\leq
r$ and 
$B_\Ga(e,4r)\cap \Ga_n=\{e\}$.
Let us define
 $$\Psi_{*,r,n}:
KK_*(P_r(\x),\C)\lto KK^{\Ga}_*(P_{2r}(\Ga),A_\Ga)$$ by
$$\Psi_{*,r}=\Psi^6_{*,2r,n}\circ\Psi_{*,2r,n}^5\circ\Psi^4_{*,2r,n}\circ\Psi^3_{*,2r,n}\circ\Psi^2_{*,r,n}\circ\Psi^1_{*,r,n}.$$
Notice that $\Psi_{*,r,n}$ does not depend on the choice of the integer $n$ such
that $n\leq
r$ and 
$B_\Ga(e,4r)\cap \Ga_n=\{e\}$.  

\medskip

For any positive real $r$ and $r'$ such that $r\leq r'$, let $$\iota^{\Ga,\a}_{*,r,r'}: KK^{\Ga}_*(P_{r}(\Ga),A_\Ga)\to
KK^{\Ga}_*(P_{r'}(\Ga),A_\Ga)$$ be  the homomorphism induced by the
inclusion $P_{r}(\Ga)\subset P_{r'}(\Ga)$.
\begin{lemma}\label{lem-surj}
For every element $y$ in $KK^{\Ga}_*(P_{r}(\Ga),A_\Ga)$, there exists
an element $x$ in
$KK_*(P_r(\x),\C)$ such that $\Psi_{*,r}(x)=\iota^{\Ga,\a}_{*,r,2r}(y)$.
\end{lemma}

\begin{proof}According to proposition
  \ref{prop-surj}, the homomorphism $\Psi^6_{*,2r,n}$  is onto. Since $\Psi^5_{*,2r,n}$,$\Psi^4_{*,2r,n}$ and
  $\Psi^3_{*,2r,n}$ are isomorphisms, there exists a $z$ in
  $\prod_{i\geq n} K_*(P_{r}(\Ga)/\Ga_i)$ such that
  $y=\Psi^6_{*,r,n}\circ\Psi_{*,r,n}^5\circ\Psi^4_{*,r,n}\circ\Psi^3_{*,r,n}(z)$, 
Using equation \ref{equ-psi1}, 
we have $$\Psi_{*,r,n}^2\circ\Lambda_{*,r,n}(z)={\iota'}_{*,r,2r}^{k\geq
  n}(z)$$ and since $\Psi_{*,r,n}^1$ is onto, there exists an element $x$
in $K_*(P_r(\x))$ such that
$\Lambda_{*,r,n}(z)=\Psi_{*,r,n}^1(x)$. The lemma is then a consequence of
the equality
$$\Psi^6_{*,2r,n}\circ\Psi_{*,2r,n}^5\circ\Psi^4_{*,2r,n}\circ\Psi^3_{*,2r,n}\circ{\iota'}_{*,r,2r}^{k\geq
  n}=\iota^{\Ga,\a}_{*,r,r'}\circ\Psi^6_{*,2r,n}\circ\Psi_{*,2r,n}^5\circ\Psi^4_{*,2r,n}\circ\Psi^3_{*,2r,n}.$$
\end{proof}
Let us denote for a pair of real $r$ and $r'$ such that $0\leq r\leq
r'$ by $\iota_{*,r,r}^\x:K_*(P_r(\x))\lto K_*(P_{r'}(\x))$ the
morphism induced by the inclusion $P_r(\x)\subset P_{r'}(\x)$. The
class $\Ga_0$ of $\Ga/\Ga_0$  can be
viewed as an element of $P_r(\x)$ and this inclusion induces a homomorphism
 $$\kappa_{*,r}:\Z\cong K_*(\{[\Ga_0]\})\lto K_*(P_r(\x)).$$
\begin{lemma}\label{lem-inj}
Let $x$ be an element of  $K_*(P_r(\x))$ such that
$\Psi_{*,r}(x)=0$, then there exists a real $r'$ such that $r\leq r'$
such that $\iota_{*,r,r'}^\x(x)$ is in the range of $\kappa_{*,r'}$.
\end{lemma}
\begin{proof}
Let us fix a integer $n$ such that $n\geq r$ and
$B_\Ga(e,4r)\cap\Ga_n=\{e\}$.
According to proposition \ref{prop-surj},
\begin{equation*}\begin{split}
\Psi_{*,2r,n}^5\circ&\Psi^4_{*,2r,n}\circ\Psi^3_{*,2r,n}\circ\Psi^2_{*,r,n}\circ\Psi^1_{*,r,n}(x)\in\\
&
KK^{\Ga}_*(P_{2r}(\Ga),C_0(\sqcup_{i\geq n}
\Ga/\Ga_i,\K(H)))\subset KK^{\Ga}_*(P_{2r}(\Ga),\li).\end{split}\end{equation*} In view of remark \ref{rem-restriction}, we
get that $$\Psi^4_{*,2r,n}\circ\Psi^3_{*,2r,n}\circ\Psi^2_{*,r,n}\circ\Psi^1_{*,r,n}(x)\in\bigoplus_{k\geq n} KK^{\Ga}_*(P_{2r}(\Ga),C(
\Ga/\Ga_k)).$$ Since $\Psi^4_{*,2r,n}$ and $\Psi^3_{*,2r,n}$ restricts on 
direct summands  to isomorphisms, then we get that $\Psi^2_{*,r,n}\circ\Psi^1_{*,r,n}(x)$
lies in  $\bigoplus_{k\geq n} KK_*(P_{2r}(\Ga)/\Ga_k,\C)$. According
to equation \ref{equ-psi2}, 
$${\iota}_{*,r,2r}^{k\geq
  n}\circ\Psi_{*,r,n}^1(x)=\Psi_{*,2r,n}^1\circ{\iota}_{*,r,2r}^{\x}(x)
 \in \bigoplus_{k\geq n} KK_*(P_{2r}(\Ga/\Ga_k),\C).$$ But then 
 ${\iota}_{*,r,2r}^{\x}(x)$ lies in  a  finite sum  of  summands of
 $\bigoplus_{k\geq n}KK_*(P_{2r}(\Ga)/\Ga_k,\C)$ and thus  we get that 
 for some integer $m$ and some real $s$ with   $s\geq 2r$ and $m\geq
 \sup\{n,s\}$, then
${\iota}_{*,r,s}^{\x}(x)$ belongs to $KK_*(P_{s}(\sqcup_{0\leq k\leq
  m}\Ga/\Ga_k),\C)$. But since $\coprod_{0\leq k\leq
  m}\Ga/\Ga_k$ is finite, $P_{s}(\sqcup_{0\leq k\leq
  m}\Ga/\Ga_k)$ is  compact and  up to choose a bigger $s$ is also
convex. Hence,
${\iota}_{*,r,r'}^{\x}(x)$ lies in the range of $\kappa_{*,r'}$ for $r'$
big enough.
\end{proof} 
It is straightforward to check that 
$\Psi_{*,r'}\circ{\iota}_{*,r,r'}^{\x}={\iota}_{*,r,r'}^{\Ga,\a}\circ
\Psi_{*,r}$ and thus the family of homomorphism $(\Psi_{*,r})_{r\leq 0}$
gives rise to a homomorphism
$$\Psi_{\x,*}:\lim_r KK_*(P_r(\x),\C)\lto K_*^\text{top}(\Ga,A_\Ga).$$ Let
$x_0$ be the image of (any) $\kappa_{*,r}(1)$ in $\lim_r
KK_*(P_r(\x),\C)$. As a consequence of lemmas  \ref{lem-surj} and
\ref{lem-inj}, we get
\begin{theorem}\label{thm-iso-psi}\
\begin{itemize}
\item $\Psi_{\x,*}$ is onto;
\item In odd degree, $\Psi_{\x,*}$ is an isomorphism;
\item In even degree, $\ker\Psi_{\x,*}$ is the infinite cyclic group
  generated by $x_0$;
\end{itemize}
\end{theorem}
\subsection{Compatibility of $\Psi_{\x,*}$ with the assembly maps}\label{sub-sec-compatibility}The
proof of equation \ref{equ-comp} require some preliminary work.
For   a locally compact and proper $\Ga$-space $X$, the notion of
standard-$X$-module, was extended to the equivariant case in
\cite{tzanev} as follows.
\begin{definition}
Let  $X$ be a locally compact and proper $\Ga$-space, let $H$ be a
$\Ga$-Hilbert space. A non-degenarated $\Ga$-equivariant
representation  $\rho:C_0(X)\to \L(H)$ is called $X$-$\Ga$-ample if when
extended to $C_0(X)\rt\Ga$, then $\rho(C_0(X)\rt\Ga)\cap \K(H)=\{0\}$.
\end{definition}
\begin{example}\label{ex-ample}
If $\eta$ is a $\Ga$-invariant measure on $\pr$ fully supported i.e
with support $\pr$ and if
$H$ is  a separable Hilbert space, then $L^2(\eta)\ts H$ equipped with the
diagonal action of $\Ga$, trivial on $H$ together with the
representation $$\rho_r:C_0(\pr)\to \L(L^2(\eta)\ts H);\,f\mapsto f\ts
\Id_H$$ is an $X$-$\Ga$-ample representation. The reason is that $\pr$
contains as a $\Ga$-space a copy of $\Ga\times Y$, where $Y$ is a open
subset of $\pr$, and where $\Ga$ acts diagonally, by left translations
on $\Ga$ and trivially on $Y$.
\end{example}
\begin{lemma}\label{lem-ample}
Let $X$ be locally compact and proper $\Gamma$-space, let $H_0$ and $H_1$
be two $\Gamma$-Hilbert spaces and let $\rho_i:C_0(X)\to \L(H_i)$ for $i=0,1$
be two non-degenerated and $\Ga$-equivariant representation. Assume that
$\rho_0$ is 
$X$-$\Ga$-ample. Then there exists
\begin{itemize}
\item $H_2$ a $\Gamma$-Hilbert space;
\item $\rho_2:C_0(X)\to \L(H_2)$
a non-degenerated
$\Ga$-equivariant representation;
\item $U:H_1\oplus H_2\to H_0$ a unitary 
\end{itemize}
 such that for every $f$ in $C_0(X)$,
$$U\cdot(\rho_1\oplus\rho_2)(f)- \rho_0(f)\cdot U\in\K(H_1\oplus
H_2,H_0).$$
\end{lemma}
\begin{proof}
Up to replace $\rho_1$ by$\rho_0\oplus\rho_1$, we can assume without
loss of generality that $\rho_1$ is also $X$-$\Ga$-ample.
Then, according to \cite{tzanev}, there exists an $\Ga$-equivariant
isometry $W:H_1\to H_0$ such that $W\cdot\rho_1(f)- \rho_0(f)\cdot W$
is in $\K(H_1,H_0)$ for every $f$ in $C_0(X)$. Let us set $P=\Id_{H_0}-W\cdot W^*$.
Now, by using the completly positive map 
$$C_0(X)\to \L (P\cdot H_0);\,f\mapsto P\cdot \rho_0(f)\cdot P,$$ we can
use the proof of \cite[Theorem 3.4.6]{hr} to conclude.
\end{proof}
From this and by using next lemma, we can prove that if $X$ is a
locally compact and proper $\Ga$-space, then every element in
$KK^\Ga_*(X,\C)$ can be represented by a $K$-cycle supported on a
prescribed  $X$-$\Ga$-ample and non-degenerated representation.
\begin{lemma}\label{lem-equ}
Let $G$ be a locally compact group and let $A$ and $B$ be two
$G$-algebras. Let $(\rho,\mathcal{E},T)$ be a $K$-cycle for
$KK^G_*(A,B)$ and let $\rho'$ be an equivariant representation of $A$
on the right $B$-Hilbert module $\E$ such that $\rho(a)-\rho'(a)$ is
compact for all $a$ in $A$. Then $(\rho',\mathcal{E},T)$ is a $K$-cycle for
$KK^G_*(A,B)$ equivalent to $(\rho,\mathcal{E},T)$.
\end{lemma}
\begin{proof}
It is clear that  $(\rho',\mathcal{E},T)$ is a $K$-cycle for
$KK^G_*(A,B)$. Then 
$$\left( \begin{smallmatrix}\cos t\pi/2&-\sin t\pi/2\\
\sin t\pi/2& \cos t\pi/2 \end{smallmatrix}\right)
\cdot \left(\begin{smallmatrix}T&0\\
0&\ \Id_E \end{smallmatrix}\right)
\cdot
\left(\begin{smallmatrix}\cos t\pi/2&\sin t\pi/2\\
-\sin t\pi/2&\ \cos t\pi/2 \end{smallmatrix}\right)_{t\in [0,1]}$$
provides a homotopy between the  $K$-cycles
$(\rho\oplus\rho',\mathcal{E}\oplus\mathcal{E},T\oplus \Id_E)$ and
 $(\rho\oplus\rho',\mathcal{E}\oplus\mathcal{E}, \Id_E\oplus T)$.
\end{proof}
\begin{corollary}\label{cor-kcycle-ample}
Let  $X$ be a locally compact and proper $\Gamma$-space and let $\rho_X$
be a $X$-$\Ga$-ample representation of $C_0(X)$ on a $\Ga$-Hilbert
space $H_X$. Then every element of $KK^\Ga_*(X,\C)$ can be
represented by a K-cycle $(\rho_X,H_X, T)$ where $T$ is a
$\Ga$-equivariant operator on $H_X$.
\end{corollary}

We fix once for all a separable Hilbert space $H$ and for each real $r$ a $\Ga$-invariant measure
$\eta_r$ on $\pr$ fully supported. Let us consider 
$H_\pr=L^2(\eta_r)\otimes H$  with the 
$X$-$\Ga$-ample representation $\rho_r$ defined in example
\ref{ex-ample}.
Define  $\Psi^{\Ga_i}(H_{\pr})$ as  the $*$-algebra  of pseudo-local,
$\Ga_i$-equivariant and finite propagation operators  on
 $H_{\pr}$. An element of  $\Psi^{\Ga_i}(H_{\pr})$ is called a K-cycle
 if it satisfies the K-cycle condition with respect to $\rho_r$.

\begin{lemma}
Let $x$ be an element of  $KK^\Ga(\pr,\li)$. Then there exists a
real $s$ and a family $(T_i)_{i\in\N}$ of bounded operators on $H_\pr$
such that
\begin{enumerate}
\item $T_i$ is a K-cycle of $\Psi^{\Ga_i}(\pr)$  of propagation
  less than $s$ and $\|T_i\|\leq 1$  for every integer $i$;
\item Under the identification
  $$\ell^\infty(X,\K(H,H_\pr))\cong\prod_{i\in\N}
  C(\Ga/\Ga_i,\K(H,H_\pr))\cong\prod_{i\in\N}\I_{\Ga_i}^\Ga\,\K(H,H_\pr),$$ then
$$((\I_{\Ga_i}^\Ga\,\rho_r)_{i\in\N},\,\ell^\infty(X,\K(H,H_\pr)),\,(\I_{\Ga_i}^\Ga\,T_i)_{i\in\N})$$
is a K-cycle that represents $x$.
\item If $x_i$ is the class of $(\rho_r,H_\pr,T_i)$ in
  $KK^{\Ga_i}_*(\pr,\C)$,
  then $\Theta^{\Ga,\A}_*(x)=(\I_{\Ga_i}^\Ga\,x_i))_{i\in\N}$, with $\A=(C(\Ga/\Ga_i))_{i\in\N}$,
\end{enumerate}
\end{lemma}
\begin{proof}
Item (iii) is a consequence of item (ii) together with proposition
\ref{prop-prod}.
According to the discussion following propositon \ref{prop-prod}, we
can assume that $x$ is represented by a K-cycle
$((\phi_i)_{i\in\N},\prod_{i\in\N}\K_{C(\Ga/\Ga_i)}(C(\Ga/\Ga_i, H),\mathcal{E}_i),(T_i)_{i\in\N})$
such that for all integer $i$
\begin{itemize}
\item $\mathcal{E}_i$ is a $\Ga$-equivariant $C(\Ga/\Ga_i)$-Hilbert
  module;
\item $T_i$ is a $\Ga$-equivariant adjointable operator on
  $\mathcal{E}_i$ with $\|T_i\|\leq 1$;
\item $\phi_i$ is a $\Ga$-equivariant representation of $C_0(\pr)$ on  $\mathcal{E}_i$;
\item   $C_0(\pr)$ and  $T_i$  act then  on
  $\K_{C(\Ga/\Ga_i)}(H,\mathcal{E}_i)$ by left
  composition.
\end{itemize}
But for every integer $i$, there exist a $\Ga_i$-Hilbert space $H_i$,
a  $\Ga_i$-equivariant representation $\psi_i$ of $C_0(\pr)$ on  $H_i$
and a $\Ga_i$-equivariant  bounded operator   $F_i$ on $H_i$ such that
$\mathcal{E}_i=\I^\Ga_{\Ga_i} H_i$, $\phi_i=\I^\Ga_{\Ga_i} \psi_i$ and $T_i=\I^\Ga_{\Ga_i} F'_i$.
Up to replace $H_i$ by the $\Ga$-Hilbert space induced by the
inclusion $\Ga_i\hookrightarrow \Ga$, we can assume that $H_i$ is a
$\Ga$-Hilbert which is up to add the degenerated K-cycle
$(\rho_r,H_\pr,\Id_{H_\pr})$ can be chosen $X$-$\Ga$-ample. By adding
the degenerated K-cycle
$$(\bigoplus_{k\in\N,\,k\neq i}\psi_k,\,\bigoplus_{k\in\N,\,k\neq i}H_k,
\bigoplus_{k\in\N,\,k\neq i}\Id_{H_k}),$$ we can assume that $H_i=H_0$ and
$\psi_i=\psi_0$. According to lemma \ref{lem-ample}, by taking an
unitary equivalence of K-cycle, we can assume that $H_0=H_\pr$ and
that
$\psi_0(f)-\rho_r(f)\in\K(H_\pr)$ for every integer $f$ in  $C_0(\pr)$.
It is straightforward to check that 
$$((\I_{\Ga_i}^\Ga\rho_r)_{i\in\N},\,\ell^\infty(X,\K(H,H_\pr)),\,(\I_{\Ga_i}^\Ga
F_i)_{i\in\N})$$
is a K-cycle for $KK^G_*(\pr,\li)$, which is by proposition
\ref{prop-prod} and lemma \ref{lem-equ} equivalent to
$$((\I_{\Ga_i}^\Ga\,\psi_0)_{i\in\N},\,\ell^\infty(X,\K(H,H_\pr)),\,(\I_{\Ga_i}^\Ga\,F_i)_{i\in\N}).$$
If $f\in C_c(\pr,[0,1])$ is a cut-off function for the action of $\Ga$
on $\pr$, then up to replace $F=(\I_{\Ga_i}^\Ga\,F_i)_{i\in\N}$ by
$\sum_{\ga\in\Ga}\gamma(f)F  \gamma(f)$,  we can assume that there
exists a real $s$ such that for all integer $i$ the operator $F_i$ has
propagation less than $s$.
\end{proof}
 Set
$\zeta=\ell^\infty(\x,\K(H,\ell^2(\Gamma\otimes H))$ and  let $\zeta_\Ga$ be  the right $\li\rtm\Ga$-Hilbert module
constructed from $\zeta$ in section \ref{sub-ass}. 
Viewing  $\ell^2(\Gamma)\otimes \ell^\infty(\x,\K(H))$ as a
right 
$\ell^\infty(\x,\K(H))$-Hilbert submodule of $
\ell^\infty(\x,\K(H,\ell^2(\Gamma)\otimes H))$, we see that 
\begin{eqnarray*}C_c(\Ga)\cdot
\ell^\infty(\x,\K(H,\ell^2(\Gamma)\otimes H))&=&C_c(\Ga)\cdot\left(
\ell^2(\Gamma)\otimes \ell^\infty(\x,\K(H))\right)\\
&\cong& C_c(\Ga)\otimes \ell^\infty(\x,\K(H))\end{eqnarray*}
and  under this identification, we get that
 $$C_c(\Ga)\otimes \ell^\infty(\x,\K(H))\to\li\rtm\Ga;\, \delta\ga\otimes
 a\mapsto \ga^{-1}(a)\delta\ga^{-1}$$ extends to isomorphism of 
right $\li\rtm\Ga$-Hilbert module
\begin{equation}\label{equ-ass-map}\zeta_\Ga\stackrel{\cong}{\lto}\li\rtm\Ga.\end{equation}
Define  for $i$ integer  $\Psi^{\Ga_i}(\Ga)$  as
   the $*$-algebras of pseudo-local,
$\Ga_i$-equivariant and finite propagation operators  on
$\ell^2(\Ga)\otimes H$.
\begin{lemma}\label{lem-si}
Let $(S_i)_{i\in\N}$ be a family in $\prod_{i\in\N}\Psi^{\Ga_i}(\Ga)$
uniformally bounded and with propagation uniformally bounded by a real
$s$.
Then under the identification $$\zeta\cong\prod_{i\in\N}
  C(\Ga/\Ga_i,\K(H,\ell^2(\Gamma)\otimes H))\cong\prod_{i\in\N}\I_{\Ga_i}^\Ga
\K(H,\ell^2(\Gamma)\otimes H)),$$
\begin{enumerate}
\item there exists a unique multiplier 
  $\lambda_\Gamma(S_i)_{i\in\N}$ of $\li\rtm\Ga$ which  under the
  identification $\zeta_\Ga{\cong}\li\rtm\Ga$    restricts to
  $(\I_{\Ga_i}^\Ga S_i)_{i\in\N}$  on  $C_c(\Ga)
  \cdot\zeta$.
\item The multiplier image of  $\lambda_\Ga(S_i)_{i\in\N}$ under the
  canonical projection
  $$\li\rtm\Ga{\longrightarrow}\a\rtm\Ga$$
and  the multiplier image of   $\bigoplus_{i\in\N}
  \widehat{S_i}$ under the map
  $$\cmx\stackrel{\Psi_\Ga}{\longrightarrow}\a\rtm\Ga$$ coincide.
\end{enumerate}
\end{lemma}
\begin{proof}
Let us prove first the lemma  for  a family $(S_i)_{i\in\N}$ of locally
compact  operators. Since such families are algebraically generated by
families
\begin{itemize}
\item $(f_i)_{i\in\N}$ uniformally bounded with  $f_i$ in $C(\Ga/\Ga_i,\K(H))$ acting by
  pointwise multiplication;
\item $(R_\ga)_{i\in\N}$, for $\ga$ in $\Ga$, where $R_\ga$ is induced by the right
  regular representation on $\ell^2(\Ga)\ts H$,
\end{itemize}
this amounts to prove the lemma  for these families.

Then the elements of   $\li\rtm\Ga$ given by
\begin{itemize}
\item  the image of $(f_i)_{i\in\N}$ viewed as element of
  $\li$ under the inclusion 
  $\li\hookrightarrow\li\rtm\Ga$ for the first case.
\item the element $\delta_\gamma\in \li\rtm\Ga$ for the family
  $(R_\ga)_{i\in\N}$, where  $\delta_\gamma$ be the Dirac function at
$\ga$.
\end{itemize}
satisfy the required property. For a family $(S_i)_{i\in\N}$ of
pseudo-local  operators,  let us set for $i$ integer
$S'_i=\sum_{\ga\in\Ga}\de_\ga S_i \de_\ga$. Then $S'_i-S_i$ is
$\Ga_i$-equivariant and locally compact for all integer $i$. Moreover,
as already mention in subsection \ref{sub-ass},  
$(\I_{\Ga_i}^\Ga S'_i)_{i\in\N}=\sum_{\ga\in\Ga}\de_\ga(\I_{\Ga_i}^\Ga S_i)_{i\in\N} \de_\ga$ 
extends to an adjointable operator of $\zeta_\Ga$ and thereby to a
multiplier $\lambda_\Ga(S'_i)_{i\in\N}$ of $\li\rtm\Ga$. We then set  
$\lambda_\Ga(S_i)_{i\in\N}=\lambda_\Ga(S'_i)_{i\in\N}+\lambda_\Ga(S_i-S'_i)_{i\in\N}$.
Unicity is quite obvious.
Since pseudo-local operator on $\ell^2(\Ga)\ts H$ are multiplier for
locally compact operator, we get item (ii) by multiplicativity of 
$(S_i)_{i\in\N}\mapsto\lambda_\Ga(S_i)_{i\in\N}$ and of 
$(S_i)_{i\in\N}\mapsto\bigoplus_{i\in\N}
  \widehat{S_i}$.
\end{proof}

we are now in position to prove the main theorem of the section.

\begin{theorem}\label{thm-iso-khom}
Let $\Ga$ be a residually finite and  finitely
generated discrete group with respect to a sequence  $\Ga_0\subset\cdots\Ga_n\subset\cdots$   of
finite index normal subgroups. Then if we set
$\x=\coprod_{n\in\N}\Ga/\Ga_i$, we have
$$\Psi_{\Ga,\a,*}\circ\mu_{\x,\text{max},*}=\mu_{\Ga,\a,\text{max},*}\circ\Psi_{\x,*}.$$
\end{theorem}
\begin{proof}
Let $r$ be a real and let $n$ be any integer such that $2r\leq n$ and
$B(e,4r)\cap\Ga_n=\{e\}$.
Then\begin{equation}\label{equ-khom}
K_*(P_r(\x))\cong K_*(P_r(\coprod_{i=1}^{n-1}\Ga/\Ga_i))\oplus
\prod_{i\geq n}K_*(P_r(\Ga/\Ga_i)).
\end{equation}
If under this identification, $x$ comes from
$K_*(P_r(\coprod_{i=1}^{n-1}\Ga/\Ga_i))$, then $\Psi_{\x,*}(x)=0$ and using the
naturality of the assembly map, we can see that $\mu_{\x,\text{max},*}$ lies in
$K_*(\cm(\coprod_{i=1}^{n-1}\Ga/\Ga_i))\subset K_*(\cmx)$ and hence
$\Psi_{\Ga,\a,*}\circ\mu_{\x,\text{max},*}(x)=0$. Thereby, we have to prove that 
$\Psi_{\Ga,*}\circ\mu_{\x,\text{max},*}(x)=\mu_{\Ga,\a,\text{max},*}\circ\Psi_*(x)$, for element
$x$ coming under the identification
of equation \ref{equ-khom} from $x'$ in $\prod_{i\geq n}K_*(P_r(\Ga/\Ga_i))$. According to equation \ref{equ-psi2}, and up to replace $r$
by $2r$, we can assume that $x'=\Lambda_{*,r,n}(y)$, with $y$ in 
 $\prod_{i\geq n}K_*(P_r(\Ga)/\Ga_i)$. We can assume indeed without loss of generality
that $n=0$ and that the action of $\Ga_i$ on $\pr$ is free for all
integer $i$.  From now on, we  will write  $\Lambda_{*}$
(resp. $\Psi^i_{*}$, $i=1,\cdots,6$) instead of
$\Lambda_{*,r,0}$ (resp. $\Psi^i_{*,r,0}$, $i=1,\cdots,6$). Let us set
$z={\Psi^3_{*}}(y)$ in $\prod_{i\in\N}KK^{\Ga_i}_*(C_0(P_r(\Ga)),\C)$. The proof of the
 theorem is divided in the following steps.
\begin{description}
\item[First Step]
Assume that $z$ is given by a family of K-cycles
$(\ror,H_{P_r(\Ga)},T_i)_{i\in\N}$ such that  for a real $s$, then
for all integer  $i$ the operator
$T_i$ is $\Ga_i$-equivariant with $\|T_i\|\leq 1$ and has propagation
  less than $s$.
Let us  set
$H_{P_r(\Ga)/\Ga_i}=L^2(\eta_{r,i})\ts H$ and $\rori:C_0(\pr/\Ga_i)\to\L(L^2(\eta_{r,i})\ts H);\,f\mapsto f\ts\Id_H$
where for all integer $i$, the measure $\eta_{r,i}$ is  induced by $\eta_r$
on $P_r(\Ga)/\Ga_i$.
Let us choose a $\Ga$-equivariant coarse map
$\widetilde{\phi_r}:\pr\to\Ga$. Then $\widetilde{\phi_r}$ is a coarse
equivalence and induces  a coarse equivalence
$$\phi_r:\coprod_{i\in\N}\pr/\Ga_i\lto
\coprod_{i\in\N}\Ga/\Ga_i=\x.$$ Let us show that  with notations of
section \ref{subsection-roe},
\begin{equation}\label{equ-form-ass}\mu_{\x,\text{max},*}(x)=\phi_{r,\text{max},*}\Ind_{\text{max},\x}\oplus_{k\in\N}\widehat{T^k},\end{equation}
where $\oplus_{k\in\N}\widehat{T^k}$ is viewed as an operator on  the
non-degenerated standard $\coprod_{k\in\N}\pr/\Ga_k$-module
$\bigoplus_{k\in\N}H_{P_r(\Ga)/\Ga_k}$ (for the representation
$\oplus_{i\in\N}\rho_{H_{P_r(\Ga)/\Ga_i}}$).

\medskip

Let $\upsilon_{r,k}:P_r(\Ga)/\Ga_k\to
P_r(\Ga/\Ga_k);\dot{h}\mapsto \widetilde{h}$ be the map defined in section
\ref{sec-khom}. Notice that the family $(\upsilon_{r,k})_{\in \N}$
induces a coarse equivalence $\upsilon_r:\coprod_{k\in\N}\pr/\Ga_k\lto\coprod_{k\in\N}
P_r(\Ga/\Ga_k)$. Moreover, if we set
$$\phi_k:C_0(P_r(\Ga/\Ga_k))\to\L(H_{P_r(\Ga)/\Ga_k});\, f\mapsto
\rork(f\circ \upsilon_{r,k}),$$ then  $x'=\Lambda_*(y)$ is the class of the K-cycle
$$\left(\bigoplus_{k\in\N}\phi_k,\,
\bigoplus_{k\in\N}H_{P_r(\Ga)/\Ga_k},\,\bigoplus_{k\in\N}\widehat{T_k}\right)$$ in
$K_*(\coprod_{k\in\N}(P_r(\Ga/\Ga_k))$. For any
non-degenerated standard $P_r(\Ga/\Ga_k)$-module $H_k$ given
by a representation $\rho_k$, then $\phi_k\oplus\rho_k$ also provides
a non-degenerated standard $P_r(\Ga/\Ga_k))$-Hilbert module structure
for $H_{P_r(\Ga)/\Ga_k}\oplus H_k$. Since the K-cycles 
$$\left(\bigoplus_{k\in\N}\phi_k,\,
\bigoplus_{k\in\N}H_{P_r(\Ga)/\Ga_k},\,\bigoplus_{k\in\N}\widehat{T_k}\right)$$
and 
$$\left(\bigoplus_{k\in\N}\phi_k\oplus \rho_k,\,
\bigoplus_{k\in\N}H_{P_r(\Ga)/\Ga_k}\oplus
H_k,\,\bigoplus_{k\in\N}\widehat{T_k}\oplus \Id_{H_k}\right)$$ are equivalent,
we get that
$$\mu_{\x,\text{max},*}(x)=\psi_{r,\text{max},*}\Ind_{\text{max},\coprod_{k\in\N}
  P_r(\Ga/\Ga_k)}\bigoplus_{k\in\N}\widehat{T_k}\oplus \Id_{H_k},$$
where $\psi_r:\coprod_{i\in\N}
  P_r(\Ga/\Ga_k)\to\x$ is any coarse equivalence. Since the inclusion
$\bigoplus_{k\in\N}H_{P_r(\Ga)/\Ga_k}\hookrightarrow\bigoplus_{k\in\N}H_{P_r(\Ga)/\Ga_k}\oplus
H_k$ covers the coarse map $\upsilon_r:\coprod_{k\in\N}
  P_r(\Ga)/\Ga_k\to \coprod_{k\in\N}
  P_r(\Ga/\Ga_k)$, we get that 

$$\Ind_{\text{max},\coprod_{k\in\N}
  P_r(\Ga/\Ga_k)}\,\left(\oplus_{k\in\N}\widehat{T_k}\oplus \Id_{H_k}\right)=\upsilon_{r,\text{max},*}\Ind_{\text{max},\coprod_{k\in\N}
  P_r(\Ga)/\Ga_k}\,\left(\oplus_{k\in\N}\widehat{T_k}\right).$$
Notice that since $\phi_r$ and $\psi_r\circ \upsilon_r$ are both
coarse equivalence between $\coprod_{k\in\N}
  P_r(\Ga)/\Ga_k$ and $\x$, then
  $\phi_{r,\text{max},*}=\psi_{r,\text{max},*}\circ
  \upsilon_{r,\text{max},*}$ and  hence we get the equality  of
  equation \ref{equ-form-ass}.
\item[Second step]
According to \cite{tzanev}, there exists an $\Ga$-equivariant
isometrie $W_r:H_\pr\to \ell^2(\Ga)\otimes H$ that covers
$\widehat{\phi_r}:\pr\to\Ga$.
Then,  if 
$W_{r,k}:H_{P_r(\Ga/\Ga_k)}\to \ell^2(\Ga/\Ga_k)\otimes H$ stands for  the
isometry induced by $W_r$ for all integer $k$,  then
$$\bigoplus_{k\in\N}W_{r,k} :\bigoplus_{k\in\N}H_{P_r(\Ga/\Ga_k)}\lto\bigoplus_{k\in\N} \ell^2(\Ga/\Ga_k)\otimes H$$
is an isometrie that covers $\phi_r$ and thus  
\begin{equation*}\begin{split}
\phi_{r,\text{max},*}\Ind_{\text{max},\coprod_{k\in\N}
  P_r(\Ga)/\Ga_k}\,\bigoplus_{k\in\N}\widehat{T_k}&=\\
\Ind_{\text{max},\coprod_{k\in\N}
 \Ga/\Ga_k}&
\,\bigoplus_{k\in\N}W_{r,k}\widehat{T_k}W_{r,k}^*+\Id_{\ell^2(\Ga/\Ga_k)\otimes
    H}-W_{r,k}W_{r,k}^*
\end{split}
.\end{equation*}

 Finally we get that
 \begin{eqnarray*}
\mu_{\x,\max,*}(x)&=&\Ind_{\text{max},\coprod_{k\in\N}
 \Ga/\Ga_k}\,\bigoplus_{k\in\N}W_{r,k}\widehat{T_k}W_{r,k}^*+\Id_{\ell^2(\Ga/\Ga_k)\otimes H}-W_{r,k}W_{r,k}^*.\\
&=&\Ind_{\text{max},\coprod_{k\in\N}
 \Ga/\Ga_k}\,\bigoplus_{k\in\N}\widehat{W_{r}T_kW_{r}^*}+\Id_{\ell^2(\Ga/\Ga_k)\otimes H}-\widehat{W_{r}W_{r}^*}
\end{eqnarray*}

\item[Third step]

By  naturallity of the assembly map,
we get that 
$$\mu_{\Ga,A_\Ga,\text{max},*}\circ\Psi_{*} (x)=
\Psi^6_{\Gamma,*}\circ\mu_{\Ga,\li,\text{max},*}(z'),$$ where
\begin{itemize}
\item  $z'$ is the element in
$K^{top}_*(\Gamma,\li)$ coming from
$\Psi^5_{*}\circ\Psi^4_{*}(z)\in KK^\Ga_*(\pr,\li)$;
\item $\Psi^6_{\Gamma}:\li\rtm\Ga\to\a\rtm\Ga$ is induced by the
  projection 
  $\Psi^6:\li\to\a$.
\end{itemize}
Let us compute $\mu_{\Ga,\li,\mx,*}(z')$.

\medskip

According to lemma \ref{lem-ample}, under the identification 
$$\ell^\infty(\x,\K(H,H_\pr))\cong\prod_{i\in\N}
  C(\Ga/\Ga_i,\K(H,H_\pr))\cong\prod_{i\in\N}\I_{\Ga_i}^\Ga\,\K(H,H_\pr),$$ 
the element $\Psi^5_{*}\circ\Psi^4_{*}(z)$ of
$KK^\Ga_*(\pr, \li)$ can  be represented by a  K-cycle
$$\left((\I_{\Ga_i}^\Ga\,\rho_r)_{i\in\N},\,\ell^\infty(X,\K(H,H_\pr)),\,(\I_{\Ga_i}^\Ga\,F_i)_{i\in\N}\right)$$
where
\begin{itemize}
\item $F_i$ is a K-cycle of $\Psi^{\Ga_i}(\pr)$ with $\|F_i\|\leq 1$ 
  for all integer $i$;
\item there exists a real $s$ such that $F_i$ has propagation less
  than $s$ for all integer $i$;
\item if  $x_i$ is the  class of $(\rho_r,H_\pr,F_i)$ in
  $KK^{\Ga_i}_*(\pr,\C)$ 
  then $z=(x_i)_{i\in\N}$.
\end{itemize}
Moreover, if we set $\E=\ell^\infty(\x,\K(H,H_\pr))$ and  $F=(\I_{\Ga_i}^\Ga\,F_i)_{i\in\N}$, we can
assume by averaging by a cut-of function for the action of $\Ga$ on
$\pr$ that  $F\cdot  C_c(\pr)\cdot \E\subset
C_c(\pr)\cdot \E$. Let us also set
$$\zeta=\ell^\infty(\x,\K(H,\ell^2(\Gamma)\otimes H))$$ and let $\E_\Ga$
and $\zeta_\Ga$ be  the right $\li\rtm\Ga$-module
constructed in section \ref{sub-ass} respectively from $\E$ and
$\zeta$. Since the isometrie $$W_\Ga:H_\pr\to \ell^2(\Gamma)\otimes H$$
has finite propagation, it induces a map
$$C_c(\pr)\cdot \E\to C_c(\Gamma)\cdot\zeta; f\mapsto
W_r\circ f,$$ which extends to an isometrie
$W_\Ga:\E_\Ga\to\zeta_\Ga$.  
As we have seen before, $\zeta_\Ga$ is a right-$\li\rtm\Ga$ module
isomorphic to $\li\rtm\Ga$ and  in view of this, $\E_\Ga$ is a direct factor of
$\li\rtm\Ga$. Moreover, if 
 $F_\Ga$ is  the
operator of $\E_\Ga$ extending $C_c(\Ga) \cdot\E\to C_c(\Ga)
\cdot\E;f\mapsto T\circ f$, then  we get with notations  of lemma \ref{lem-si},
that  $$\lambda_\Ga(W_rF_iW_r^
*+\Id_{\ell^2(\Ga)\otimes H}-W_rW_r^
*)_{i\in\N}=W_\Ga F_\Ga\cdot W_\Ga^*+\Id_{\zeta_\Ga}-W_\Ga
W_\Ga^*.$$ Hence  $\mu_{\Ga,\li,*}(z')$ is the class in $K_*(\li\rtm\Ga)$
of the K-cycle $(\li\rtm\Ga, \lambda_\Ga(W_rF_iW_r^
*+\Id_{\ell^2(\Ga)\otimes H}-W_rW_r^
*)_{i\in\N})$.
The theorem is then a consequence of lemma \ref{lem-si}.

\end{description}
\end{proof}

\subsection{Applications}
We end this section with application concerning injectivity and
bijectivity of the maximal coarse Baum-Connes assembly map.
Let $\Ga$ be a residually finite and  finitely
generated discrete group with respect to a fixed sequence
$\Ga_0\subset\cdots\Ga_n\subset\cdots$  of
finite index normal subgroups. Recall that we have defined
$\x=\du_{i\in\N}\Ga/\Ga_i$ and $\a=\l^\infty(\x,\K(H))/C_0(\x,\K(H))$.
We can formulate  corollary  \ref{cor-suitex}, theorem \ref{thm-iso-psi} and
theorem \ref{thm-iso-khom} together as follows:
We have a commutative diagram 
$$\begin{CD}
  0 @>>>\Z @>>>\lim_r
K_0(P_r(\x)) @>\Psi_{\x,*}>> K^\text{top}_0(\Ga,A_\Ga) 
  @>>>  0\\
    & &   @V=VV   
         @V\mu_{\x,\text{max},*}VV  @V\mu_{\Ga,\a,\text{max},*}VV\\
0 @>>>\Z @>>> K_0(\cmx) @>\Psi_{\Ga,\a,*}>> K_0(\a\rtm\Ga)
  @>>>  0
\end{CD}$$ with exact rows and a commutative diagram
$$\begin{CD}
\lim_r KK_1(P_r(\x)) @>\Psi_{\x,*}>\cong> K^\text{top}_1(\Ga,A_\Ga)  \\  
         @V\mu_{\x,\text{max},*}VV @V\mu_{\Ga,\a,\text{max},*}VV\\K_1(\cmx) @>\Psi_{\Ga,\a,*}>\cong> K_1(\a\rtm\Ga)
\end{CD}$$
From this commutative diagram, we can deduce the following series of
results concerning injectivity and bijectivity of assembly maps.
\begin{theorem}\label{thm-equ} The following assertions are equivalent:
\begin{itemize}
\item[(i)]   
The maximal coarse assembly map $$\mu_{\x,\text{max},*}:\lim_r
K_*(P_r(\x),\C)\to  K_*(\cmx)$$ is an
isomorphism.
\item[(ii)] the maximal assembly map   
$$\mu_{\Ga,\a,\text{max}}:K^\text{top}_*(\Ga,A_\Ga)\to K_*(\a\rtm\Ga)$$
is an isomorphism.
\end{itemize}\end{theorem}
Example of groups that satisfies item (ii) of the theorem are provided
by groups that satisfy the so called {\em strong Baum-Connes
  conjecture}.  Recall first  that
a $\Ga$-algebra $D$ is said to be a {\em proper} $\Ga$-algebra, if $D$ is a $C_0(Z)$-algebra for some proper $\Ga$-space $Z$ in such a way that the structure map $\Phi:C_0(Z)\to ZM(D)$
is $\Ga$-equivariant. A group $\Ga$ satisfies the strong Baum-Connes conjecture
if there exist, a proper $\Ga$-algebra $D$, an element $\alpha$ in
$KK_*^\Ga(D,\C)$  and a element $\beta$ in 
$KK_*^\Ga(\C, D)$  such that $\beta\otimes_D\alpha$ is the unit of  $KK_*^\Ga(\C,\C)$.
It is well know  (see \cite{tumoy} for instance) that if $\Ga$
satisfies the strong Baum-Connes conjecture, then
$\mu_{\Ga,B,\text{max}}:K^\text{top}_*(\Ga,B)\to K_*(B\rtm\Ga)$ is an
isomorphism for every $\Ga$-algebra $B$.
As a consequence, we get
\begin{corollary}\label{cor-strongBC}
If $\Ga$ satisfies the strong Baum-Connes conjecture, then 
 $$\mu_{\x,\text{max}}:\lim_r
K_*(P_r(\x),\C)\to  K_*(\cmx)$$ is an
isomorphism.
\end{corollary}
As a particular case, we obtain

\begin{corollary}
If $\Ga$ is a group with the Haagerup property, then
 $$\mu_{\x,\text{max},*}:\lim_r
K_*(P_r(\x),\C)\to  K_*(\cmx)$$ is an
isomorphism.
\end{corollary}
\begin{example}
If $\Ga=SL_2(\Z)$ and $\Ga_k=\ker:SL_2(\Z)\lto SL_2(\Z/k\Z)$, then 
$\mu_{\x,\text{max},*}$ is an isomorphism while
$\mu_{\x,*}$ is not surjective.
\end{example}
\begin{remark}
If the group $\Ga$ has the Kazdhan property (T), then the family of
projectors corresponding to the  $0$-eigenvalue of the Laplacians of
the family $(\Ga/\Ga_i)_{i\in\N}$ provides a projector $p$ in $\cmx$ \cite{lubotzky}.
We know from \cite{hls} that the image in $K_0(C^*(\x))$ of the class
of $p$ under the homomorphism $$\lambda_{\x,*}:K_0(\cmx)\to
K_0(C^*(\x))$$ is not in the range of the coarse assembly map $$\mu_{\x,*}:\lim_r
K_*(P_r(\x),\C)\to  K_*(C^*(\x)).$$ Hence, according to remark
\ref{rem-regular}, the assembly map $\mu_{\x,\text{max},*}$ is not surjective.
\end{remark}
Recall from \cite{tumoy} that if the group $\Ga$ satisfies the strong Baum-Connes
conjecture, then $\Ga$ is K-amenable and in particular, the K-theory
of reduced and maximal crossed product coincide. This allows to get
explicit computation for $K_*(\cmx)$ in the following situation.
\begin{corollary}\label{cor-exact-sequence}
If $\Ga$ satisfies the strong Baum-Connes conjecture and admits a
universal example for proper action which is simplicial, with
simplicial and cocompact action of $\Ga$, then we have a short exact
sequence
$$0\to \Z\to K_0(\cmx)\to \prod_{i\in\N}
K_0(C_{\text{red}}^*(\Ga_i))/\bigoplus_{i\in\N}
  K_0(C_{\text{red}}^*(\Ga_i))\to 0$$ and an isomorphism

$$ K_1(\cmx)\stackrel{\cong}{\longrightarrow}\prod_{i\in\N}
K_1(C_{\text{red}}^*(\Ga_i))/\bigoplus_{i\in\N}
  K_1(C_{\text{red}}^*(\Ga_i)).$$
\end{corollary}
\begin{proof}
First notice that since $\Ga$ is K-amenable, then
$$\lambda_{\Ga,\a,*}:K_*(\a\rtm\Ga)\to K_*(\a\rtimes_{\red}\Ga).$$
Let us show that  we have an isormorphim 
$$K_*(\a\rtimes_{\red}\Ga)\stackrel{\cong}{\longrightarrow}\prod_{i\in\N}
K_*(C_{\text{red}}^*(\Ga_i))/\bigoplus_{i\in\N}
  K_*(C_{\text{red}}^*(\Ga_i)).$$
Let us consider the following commutative diagram
$$\begin{CD}
 K_*^\text{top}(\Ga,\prod_{i\in\N}C(\Ga/\Ga_i,\K(H))) @>\mu_{\Ga,\prod_{i\in\N}C(\Ga/\Ga_i)\Ga),\text{red},*}>>K_*((\prod_{i\in\N}C(\Ga/\Ga_i,\K(H)))\rtimes_\text{red}\Ga)\\
   @VVV    @VVV\\
 \prod_{i\in\N}K_*^\text{top}(\Ga,C(\Ga/\Ga_i))
 @>\prod_{i\in\N}\mu_{\Ga,C(\Ga/\Ga_i,\K(H))\Ga),\text{red},*}>>\prod_{i\in\N}K_*(C(\Ga/\Ga_i,\K(H))\rtimes_\text{red}\Ga)
\end{CD},$$ where the vertical arrow are induced on the $k$-th factor
by the projection $$\prod_{i\in\N}C(\Ga/\Ga_i,\K(H))\to
C(\Ga/\Ga_k,\K(H)).$$ But since  the group $\Ga$ admits a
universal example for proper action which is simplicial, with
simplicial and cocompact action of $\Ga$, then the left vertical arrow
is an isomorphism. Since $\Ga$ satisfies the Baum-Connes conjecture,
then the horizontal map are also isomorphism. Hence the right vertical
map is also an isomorphism and hence the result is a consequence of
remark \ref{rem-intertwin} and of the Morita equivalence between
$C(\Ga/\Ga_i)\rtimes_\text{red}\Ga$ and $C^*_\red(\Ga_i)$.
\end{proof}

Regarding injectivity, we have similar results.
\begin{theorem}The following assertions are equivalent:
\begin{itemize}
\item[(i)]   
The maximal coarse assembly map $$\mu_{\x,\text{max},*}:\lim_r
K_*(P_r(\x),\C)\to  K_*(\cmx)$$ is injective.
\item[(ii)] the maximal assembly map   
$$\mu_{\Ga,\a,\text{max},*}:K^\text{top}_*(\Ga,A_\Ga)\to K_*(\a\rtm\Ga)$$
is injective.
\end{itemize}\end{theorem}

We can also deduce the following result concerning the (usual) coarse
Baum-Connes conjecture.

\begin{theorem}\label{thm-inj} Assume that the assembly map 
$$\mu_{\Ga,\a,\text{red},*}:K^\text{top}_*(\Ga,A_\Ga)\to K_*(\a\rtimes_{\text{red}}\Ga)$$
is injective . Then the   
 coarse assembly map $$\mu_{\x,*}:\lim_r
K_*(P_r(\x),\C)\to  K_*(C^*(\x))$$ is also injective.
\end{theorem}
\begin{proof}
 In the even case, let us consider the
following diagram
$$\begin{CD}
 0 @>>>\Z @>>>\lim_r
KK^\Ga_0(P_r(\x),\C) @>\Psi_{\x,*}>> K^\text{top}_0(\Ga,A_\Ga) 
  @>>>  0\\
&&   @V=VV  
         @V\mu_{\x,\text{max},*}VV  @V\mu_{\Ga,\a,\text{max},*}VV\\
 0 @>>>\Z @>>> K_0(\cmx) @>\Psi_{\Ga,\a,\mx,*}>> K_0(\a\rtm\Ga)
  @>>>  0\\
 & &   @V=VV   
         @V\lambda_{\x,*}VV  @V\lambda_{\Ga,A_\Ga,*}VV\\
   &&\Z @>>> K_0(C^*(\x)) @>\Psi_{\Ga,\a,\text{red},*}>> K_0(\a\rtimes_{\text{red}}\Ga)
\end{CD}$$
where the bottom left corner horizontal arrow is induced by the inclusion
$$\K(\ell^2(\x)\ts H)\hookrightarrow C^*(\x)$$ and is according to
remark \ref{prop-inj-cpt} injective. Thereby, since the top row is
exact, we get that injectivity of
$\mu_{\Ga,\a,\max,*}=\lambda_{\Ga,A_\Ga}\circ \mu_{\Ga,\a,\text{max},*}$
implies injectivity of $\mu_{\x,*}=\lambda_{\x}\circ
\mu_{\x,\text{max},*}$
\end{proof}
It was proved in \cite{sty} that for a group $\Ga$ which embeds
uniformally in a Hilbert space, then $\mu_{\Ga,B,*}$ is injective for
any $\Ga$-algebra $B$.
As a consequence we obtain
\begin{corollary}
Let $\Ga$ be  a group uniformally embeddable in a Hilbert space, then
the   
 coarse assembly map $$\mu_{\x,*}:\lim_r
K_*(P_r(\x),\C)\to  K_*(C^*(\x))$$ is  injective.
\end{corollary}
The last application is to rational injectivity of
$\mu_{\x,\text{max},*}$. Theorem \ref{thm-inj} admits an obvious rational
version. This allowed to recover the following result of \cite{gwy}
\begin{theorem}
Assume that $\Ga$ admits a universal example for proper action which
is simplicial and with simplicial and cocompact action of $\Ga$. If
$\mu_{\Ga,\C,\text{max},*}$ is rationnaly injective, then
$\mu_{\x,\text{max},*}$ is also rationnaly injective.
\end{theorem}
\begin{proof}
Since rationnal injectivity of $\mu_{\bullet,\C,\text{max},*}$ is
inherited by finite index subgroups, we get under the hypothesis of
the theorem that $\mu_{\Ga_i,\C,\text{max},*}$ is rationnaly injective
for all integer $i$. Since  assembly maps are compatible with
induction, we get that
$\mu_{\Ga,C(\Ga/\Ga_i),\text{max},*}$ is rationnaly injective.
According to corollary \ref{cor-ass-prod} and since we have the
commutative diagram
\begin{equation}
\label{diag-ratio-inj}
\begin{CD}
 K^{\text{top}}_*(\Ga,\li)@>>\mu_{\Ga,\li,\text{max},*}> K_*(\li\rtm\Ga)\\
   @VVV   @VVV\\
\prod_{i\in \N}K^{\text{top}}_*(\Ga,C(\Ga/\Ga_i))
@>(\mu_{\Ga_i,C(\Ga/\Ga_i),\text{max},*})_{i\in \N}>> \prod_{i\in \N} K_*(C(\Ga/\Ga_i)\rtm\Ga) 
\end{CD}\end{equation} where the vertical arrows are induced on 
 the $k$-th factor up to Morita equivalence by  the projection
$\Pi_{i\in \N}\li\to C(\Ga/\Ga_k,\K(H))$, we see that

$$\mu_{\Ga,\li,\text{max},*}:K^{\text{top}}_*(\Ga,\li)\to K_*(\li\rtm\Ga)$$
is rationnaly injective.
As we have already seen before, the assembly is also compatible with
direct sum of coefficients. Hence we get that
$$\mu_{\Ga,C_0(\x,\K(H)),\text{max},*}:K^{\text{top}}_*(\Ga,C_0(\x,\K(H)))\to
K_*(C_0(\x,\K(H)\rtm\Ga)$$ is also rationnaly injective. By using the
maps  induced for each integer $k$ by the $k$-th factor $\Pi_{i\in \N}\li\to C(\Ga/\Ga_k,\K(H))$, with see
that
the inclusion $C_0(\x,\K(H))\hookrightarrow \li$  induces inclusions
$$K^{\text{top}}_*(\Ga,C_0(\x,\K(H)))\hookrightarrow
K^{\text{top}}_*(\Ga,\li)$$ and $$K_*(C_0(\x,\K(H))\rtm\Ga)\hookrightarrow
K_*(\li\rtm\Ga)$$ and thus we get a commutative diagram 
{\tiny{$$\begin{CD}
 0 @>>>K^{\text{top}}_*(\Ga,C_0(\x,\K(H))) @>>>K^{\text{top}}_*(\Ga,\li)\\
&&   @VVV  
         @VVV \\
 0 @>>> K_*(C_0(\x,\K(H))\rtm\Ga) @>>>K_*(\li\rtm\Ga)\\  
&&&@>>> K^\text{top}_*(\Ga,A_\Ga) 
 @>>>  0\\
&&&&& @VVV\\
&&&@>>> K_*(\a\rtm\Ga)
  @>>>  0
\end{CD}$$}}
 with exact rows and where the vertical arrows are given by
the assembly maps.
Using once again the commutativity of diagram \ref{diag-ratio-inj}, we
get that if $\mu_{\Ga,\a,\text{max},*}(x)$ comes rationnally from an
element in $K_*(C_0(\x,\K(H))\rtm\Ga)$, then $x$ comes
rationnaly from an element in $K^{\text{top}}_*(\Ga,C_0(\x,\K(H)))$. Hence $$\mu_{\Ga,\a,\text{max},*}: K^\text{top}_*(\Ga,A_\Ga) \to
K_*(\a\rtm\Ga)$$ is rationally injective.
\end{proof}
\section{Asymptotic  quantitative Novikov/Baum-Connes conjecture}
Corollary \ref{cor-exact-sequence} suggest that the property of the
coarse assembly map  $$\mu_{\x,*}:\lim_r
K_*(P_r(\x),\C)\to  K_*(C^*(\x))$$ should be closely related to the
family fo assembly maps
$$\left(\mu_{\Ga_i,\C,\text{max},*}:K^{\text{top}}_*(\Ga_i,\C)\to
  K_*(C_\mx^*(\Ga_i)\right)_{i\in\N}.$$ We this, we introduce some
quantitative assembly maps which take into account the propagation. The
relevant propagation here is indeed the one induced by $\Ga$ under the
Morita equivalence between $C_\mx^*(\Ga_i)$ and
$C(\Ga/\Ga_i)\rtm\Ga$. We then give asymptotic statements for these
quantitative assembly maps and give examples of group for which they
are satisfied.
\subsection{Almost projectors, almost unitaries  and propagation}
\begin{definition}
Let $A$ be a unital $C^*$-algebra and let $\eps$ in $(0,1/4)$.
\begin{itemize}
\item An element $p$ in $A$ is called {\bf an $\eps$-projector}  if  $p=p^*$ and $\|p^2-p\|<\eps$.
\item An element $u$ in $A$ is called {\bf an $\eps$-unitary} if
  $\|u^*u-1\|<\eps$ and  $\|uu^*-1\|<\eps$
\end{itemize}
\end{definition}
Notice that if $p$ is an $\eps$-projector of a $C^*$-algebra $A$ and
$\phi_{\eps}:\R\to\R$ is  any continuous function such that
$\phi_{\eps}(t)=0$ for $t<\frac{1-\sqrt{1-4\eps}}{2}$ and $\phi_{\eps}(t)=1$ for
$t>\frac{1+\sqrt{1-4\eps}}{2}$, then $\phi_{\eps}(p)$ is a projector.  Moreover, we have $\|\phi_\eps(p)-p\|\leq 2\eps$.
\begin{remark}\label{rem-closed}If $A$ is a $C^*$-algebra, then for every
$\eps$ in $(0,1/4)$  
\begin{itemize}
\item if $p$ is an $\eps$-projector of $A$,   then any  element $q$ in
  $A$ such that $\|p-q\|<\frac{\eps-\|p^2-p\|}{4}$ is an
$\eps$-projector. In this case $(tp+(1-t)q)_{t\in[0,1]}$ is a homotopy
  of $\eps$-projectors between $p$ and $q$  and in consequence
  $\phi_\eps(p)$ and $\phi_\eps(q)$ are homotopic  projectors.
\item if $A$ is unital and if $u$ is an  $\eps$-unitary of $A$,   then any
  element $v$ such that 
$\|u-v\|<\frac{\eps-\max\{\|u^*u-1\|,\|uu^*-1\|\}}{3}$ is
$\eps$-unitary and $(tu+(1-t)v)_{t\in[0,1]}$ is a homotopy of
  $\eps$-unitary connecting $u$ and $v$.
\end{itemize}
\end{remark}

\begin{definition}
Let $A$ be a $\Ga$-algebra. An element $x$ of $A\rtm\Ga$ is said to be of
finite propagation if $x$  lies in $C_c(\Ga,A)$. We say that  $x$ has
propagation less than $r$ if the support of $x$ as an element of
$C_c(\Ga,A)$ is  in $B_\Ga(e,s)$.
These definitions have an obvious extension to
$\widetilde{A\rtm\Ga}$ by requiring the unit to be  of  propagation zero.
\end{definition}
For $\eps$ in $(0,1/4)$, $A$ a unital $\Ga$-algebra,  and $p_0$ and $p_1$ two
$\eps$-projectors of 
$A\rtm\Ga\ts\K(H)$ with propagation less than $s$  and  $n_0$ and $n_1$   positive integers,
we write $(p_0,n_0)\sim_{s,\eps}(p_1,n_1)$ if there is an integer $k$
and a $\eps$-projector homotopy 
$(q_t)_{t\in[0,1]}$   in 
$C([0,1],(A\rtm\Ga)\ts\K(H))$ between $\left(\begin{smallmatrix} p_0&0\\
  0&I_{k+n_1}\end{smallmatrix}\right)$
 and  $\left(\begin{smallmatrix} p_1&0\\
  0&I_{k+n_0}\end{smallmatrix}\right)$  such that $q_t$
has propagation less than $s$ for every $t$ in $[0,1]$.
Similarly if   $u_0$ and $u_1$ are $\eps$-unitaries in
$\widetilde{A\rtm\Ga\ts\K(H)}$ with propagation less than $s$, 
we write $u_0\sim_{s,\eps}u_1$ if there is a $\eps$-unitary homotopy 
$(v_t)_{t\in[0,1]}$   in 
$C([0,1],\widetilde{(A\rtm\Ga)\ts\K(H))}$ between $u_0$ and $u_1$ and such that $v_t$
has propagation less than $s$ for every $t$ in $[0,1]$.

Notice that if $p$ and $q$ are two $\eps$-projectors in
$A\rtm\Ga\ts\K(H)$ such that $p_0\sim_{s,\eps}p_1$ then
$\phi_\eps(p_0)$ and $\phi_\eps(p_1)$ are homotopic projectors.

\subsection{Propagation and assembly map}\label{subsec-prop}
As before 
$\Ga$ is  a  finitely generated group which is residually finite  with
respect to a family 
$\Ga_0\supset\Ga_1\supset\ldots\Ga_n\supset\ldots$ of  normal finite
index subgroups of $\Gamma$.

\smallskip

Recall that   $\Psi^{\Ga_i}(\Ga)$ and  $\Psi^{\Ga_i}({\pr})$ are
respectively   the $*$-algebras of pseudo-local,
$\Ga_i$-equivariant and finite propagation operators  on
$\ell^2(\Ga)\otimes H$ and $H_{\pr}$.
 Recall from section \ref{subsec-lhs-ind}
and corollary \ref{cor-kcycle-ample} that any element of
$KK^\Ga_*(\pr,C(\Ga/\Ga_i))$ can be represented by K-cycle
$(\I_{\Ga_i}^\Ga\,\rho_r,C(\Ga/\Ga_i, H_\pr),\I_{\Ga_i}^\Ga\,F)$, where 
\begin{itemize}
\item $\rho_r$ is the standard representation of $C_0(\pr)$ on
  $H_\pr=L^2(\pr,\eta_r)\ts H$;
\item $F$ is a K-cycle of  $\Psi^{\Ga_i}(H_{\pr})$.
\item We have identified $\I_{\Ga_i}^\Ga\,H_\pr$ with
  $C(\Ga/\Ga_i,H_\pr)\cong  C(\Ga/\Ga_i)\otimes H_\pr $ provided with the
  diagonal action of $\Ga$;
\item Under this identification, $\I_{\Ga_i}^\Ga\,\rho_r$ is the
  pointwise representation $\rho_r$ on $H_\pr$ and $\I_{\Ga_i}^\Ga\,F$
  is the pointwise action by $\Ga/\Ga_i\to B(H_\pr);\,\ga\Ga_i\mapsto
  \ga(F)$.
\end{itemize}
For a K-cycle $F$  of $\Psi^{\Ga_i}(\pr)$, let us denote by $x_F$ the corresponding element in
$K^\text{top}(\Ga,C(\Ga/\Ga_i))$ coming from the K-cycle
$(\I_{\Ga_i}^\Ga\,\rho_r,C(\Ga/\Ga_i, H_\pr),\I_{\Ga_i}^\Ga\,T_i)$.

Let us set  
$\zeta_i=\I_{\Ga_i}^\Ga\,\ell^2(\Ga)\otimes
H\cong C(\Ga/\Ga_i,\ell^2(\Ga)\otimes H)$ and let $\zeta_{i,\Ga}$ be
the right $C(\Ga/\Ga_i)\rtm\Ga$-Hilbert module constructed from
$\zeta_i$ in section
\ref{sub-ass}. Notice that $\zeta_{i,\Ga}$ as a  right
$C(\Ga/\Ga_i)\rtm\Ga$-Hilbert module is isomorphic to
$H\otimes C(\Ga/\Ga_i)\rtm\Ga$ (compare with isomorphism of equation
\ref{equ-ass-map}).

 Proceeding as we did for proving  lemma
\ref{lem-si} and denoting  the
multiplier algebra of  $C(\Ga/\Ga_i,\K(H))\rtm\Ga$  by $M(C(\Ga/\Ga_i,\K(H))\rtm\Ga)$,  we get  with notations of section
\ref{sub-sec-compatibility}
\begin{lemma}\label{lem-lambda-i}
For every integer $i$, there is a $*$-homomorphism
$$\lambda_i:\Psi^{\Ga_i}(\Ga)\to M(C(\Ga/\Ga_i,\K(H)\rtm\Ga))$$ such
that
\begin{itemize}
\item   Under the identification $\zeta_{i,\Ga}\cong H\otimes
  C(\Ga/\Ga_i)\rtm\Ga$, then $\lambda_i(S)$  restricts to
  $\I_{\Ga_i}^\Ga\,S$ on $C_c(\Ga)\cdot\zeta_{i,\Ga}$;
\item For any $f$ in $C(\Ga/\Ga_i,\K(H))$, viewed as a locally compact
  operator on $\ell^2(\Ga)\otimes H$, then $\lambda_i(f)$ is  the
  image of $f$ under the inclusion $C(\Ga/\Ga_i,\K(H))\hookrightarrow C(\Ga/\Ga_i,\K(H))\rtm\Ga$;
\item For every $\ga$ in $\Ga$, then $\lambda_i(R_\ga)$ is the left
  multiplication by $\de_\ga\in C(\Ga/\Ga_i)\rtm\Ga$ (viewed as a
  multiplier of  $C(\Ga/\Ga_i,\K(H)\rtm\Ga)$);
\end{itemize}
\end{lemma}
\begin{remark}\label{rem-norm-equ}
Let us denote by $C[\Ga]^{\Ga_i}$ the set of  $\Ga_i$-equivariant
operators 
of $C[\Ga]$. Then $C[\Ga]^{\Ga_i}$ is a $*$-algebra isomorphic to
$C_c(\Ga,C(\Ga/\Ga_i))$ (equiped with convolution product) and thus
$\lambda_i$ induces by restriction a homomorphism 
$$C((\Ga/\Ga_i)\rtm\Ga\to M(C(\Ga/\Ga_i,\K(H)\rtm\Ga)),$$ which is in
fact the natural inclusion. According to \cite[lemma 4.13]{gwy}, lemma
\ref{lem-norm} can be generalised to the equivariant case and 
hence 
for every real $t$ and any integer $i$, there exists a positive   real $C_{t,i}$
 such that for any element $S$ of  $C[\Ga]^{\Ga_i}$ with propagation less
 than $t$, then 
$\|\lambda_i(S)\|_{C(\Ga/\Ga_i,\K(H)\rtm\Ga)}\leq C_{t,i}
\|S\|_{\ell^2(\Ga)\otimes H}$.
\end{remark}
Let us fix until the end of this subsection
\begin{itemize}
\item a $\Ga$-equivariant coarse equivalence
  $\widehat{\phi_r}:\pr\to\Ga$;
\item a isometry $W_r:H_r\to\ell^2(\Ga)\otimes H$ that covers
  $\widehat{\phi_r}$.
\end{itemize}
By using the same argument as in the third step of the proof of theorem \ref{thm-iso-khom}  we get the
following result.
\begin{proposition}\label{prop-ass-max}
Let $x$ be  in $K^\text{top}_*(\Ga,C(\Ga/\Ga_i))$ coming from
an element $x_F$ in  $KK^\Ga_*(\pr,C(\Ga/\Ga_i)))$ for a  K-cycle $F$ in
some 
$\Psi^{\Ga_i}(\pr)$. Then $\mu_{\Ga,C(\Ga/\Ga_i),\text{max},*}(x)$ is the class in
$K_*(C(\Ga/\Ga_i)\rtm\Ga)$ of  the K-cycle $$(H\otimes
C(\Ga/\Ga_i)\rtm\Ga, \lambda_i(W_r F W_r^*+ \Id_{\ell^2(\Gamma)\otimes
  H}-W_r  W_r^*)).$$

\end{proposition}
Define  $F_{r,i}=\lambda_i(W_r F W_r^*+
\Id_{\ell^2(\Gamma)\otimes H}-W_r  W_r^*
)$ for $F$ a K-cycle  of $\Psi^{\Ga_i}(\pr)$  and let us set  in
the even case
\begin{equation*} V_{F}=\begin{pmatrix} \Id_{H}&F_{r,i}\\0& \Id_{H}\end{pmatrix}
\cdot \begin{pmatrix} \Id_{H}&0\\-F_{r,i}& \Id_{H}\end{pmatrix}
\cdot \begin{pmatrix} \Id_{H}&F_{r,i}\\0& \Id_{H}\end{pmatrix}\cdot\begin{pmatrix}0&- \Id_{H}\\
   \Id_{H}&0\end{pmatrix},\end{equation*}where $\Id_H$ is viewed as the unit of $C(\Ga/\Ga_i,\K(H))\widetilde{\rtm}\Ga$.
Since $\lambda_i$ is a $*$-homomorphism, we see that the matrix 
$$V_F\begin{pmatrix} \Id_{H}&0\\0&0\end{pmatrix}V_F^{-1}-
\begin{pmatrix} \Id_{H}&0\\0&0\end{pmatrix}$$ has coefficients in
$C(\Ga/\Ga_i,\K(H))\rtm\Ga$ and moreover  we get
\begin{proposition}\label{prop-ass-max-even}
With notations of proposition \ref{prop-ass-max}, if  $x$   in
$K^\text{top}_*(\Ga,C(\Ga/\Ga_i))$ comes from
an element $x_F$ in  $KK^\Ga_*(\pr,C(\Ga/\Ga_i))$ for a  K-cycle $F$ in
some 
$\Psi^{\Ga_i}(\pr)$ and
if we set $e_F=V_F\begin{pmatrix} \Id_{H}&0\\0&0\end{pmatrix}V_F^{-1}$,
then $$\mu_{\Ga,C(\Ga/\Ga_i),\text{max},0}(x)=\left[e_F\right]-\left[\left(
\begin{smallmatrix} \Id_{H}&0\\0&0\end{smallmatrix}\right)\right].$$
\end{proposition}
The crucial point is that with notations of proposition
\ref{prop-ass-max-even}, then the coefficients of the idempotent $e_F$ have indeed finite
propagation depending only on the propagation of $F$.
Since with notation of lemma
\ref{lem-lambda-i}, the algebra $C[\Ga]^{\Ga_i}$ is   generated by 
$R_\ga$ for $\ga$ in $\Ga$ an by functions $f$ in $C(\Ga/\Ga_i,\K(H))$,
it is straightforward to check that $\lambda_i$ is propagation
preserving. Using this, we obtain for every  positive real $r$ the existence of a non-decreasing function
$h_r:\R^+\to\R^+$,(which is in fact affine),  independant on  $i$,  such that for every $s$
and every K-cycle  $F$ in  $\Psi^{\Ga_i}(\Ga)$ with propagation less
than $s$, then with notation of lemma
\ref{prop-ass-max-even}, the idempotent $e_F$ has  propagation less than $h_r(s)$.
Notice that $e_F$ has operator norm less than $\alpha_{r,i}=(1+\|F_{r,i}\|)^6$. Recall
that if we set
$e'_F=(1+(2e_F-1)(2e^*_F-1))^{-1/2}e(1+(2e_F-1)(2e^*_F-1))^{1/2}$,
then $e'_F$ is a projector equivalent to $e_F$.
Fix once for all  two  sequences of real polynomial functions
$(P_j)_{j\in\N}$ and $(Q_j)_{j\in\N}$ such that $P_j$ and  $Q_j$ have
degre $j$ for all $j\in\N$ and on every compact
subset of $\R^+$, 
\begin{itemize}
\item $(P_j)_{j\in\N}$ converges uniformally to $t\mapsto \sqrt{1+t}$;
\item $(Q_j)_{j\in\N}$ converges uniformally to $t\mapsto
  \frac{1}{\sqrt{1+t}}$.
\end{itemize}
Let us define $\Psi_1^{\Ga_i}(\pr)=\{T\in\Psi^{\Ga_i}(\pr)\text{ such
  that }\|T\|\leq 1\}$.
 For   $F$ a K-cycle of $\Psi_1^{\Ga_i}(\pr)$,  a positive real $r$
 and $\eps$ in $(0,1/4)$ let $j_{\eps,F,i}$ be the smallest integer such that
$$|P_j(t)-\sqrt{1+t}|\leq ((8+4\alpha_{r,i}\|)\|+2)^{-2} 
\text{ and } \left|Q_j(t)-\frac{1}{\sqrt{1+t}}\right|\leq ((8+4\alpha_{r,i}))^{-2}$$ for all integer $j\geq j_{\eps,F,i}$ and all
$t\in[0,4\alpha_{r,i}]$. For   $F$ a K-cycle of $\Psi_1^{\Ga_i}(\pr)$ with
propagation less than $s$, let us set
$$\tilde{p}_{F,\eps}=1/2Q_{j_{\eps,F,i}}((2e_F-1)(2e_F^*-1))(e_F+e_F^*)P_{j_{\eps,F,i}}((2e_F-1)(2e_F^*-1)).$$
Then $\|e'_F-\tilde{p}_{F,\eps}\|<\eps/8$ and according to remark
\ref{rem-closed}, then  $\tilde{p}_{F,\eps}$ is a $\eps$-projection and has propagation less than
$(2j_{\eps,F,i}+1)h_r(s)$.
Moreover, for any continuous function $\phi_\eps:\R\to\R$ such that
$\phi_\eps(t)=0$ for $t<\frac{1-\sqrt{1-4\eps}}{2}$ and $\phi_\eps(t)=1$ for
$t>\frac{1+\sqrt{1-4\eps}}{2}$, then $\phi_\eps(\tilde{p}_{\eps})$ is a projector
equivalent to $e_F$.
Now fix an identification between $\K(H)$ and the closure of
$\cup_{n\in\N}M_n(\C)$ and consider $q_n$ the rank $2n$ projector  of $M_2(C(\Ga/\Ga_i,\K(H))\widetilde{\rtm}\Ga))$
corresponding to the identity of $M_2(M_n(\C))$.
For a  K-cycle $F$ of $\Psi_1^{\Ga_i}(\pr)$, let $n_{F,\eps}$
be the smaller integer $n$  such that
$\|q_n\tilde{p}_{F,\eps/2}q_n-\tilde{p}_{F,\eps/2}+\left(\begin{smallmatrix}
    \Id_{H}-I_n&0\\0&0\end{smallmatrix}\right)\|<\eps/8$ and set 
${p}_{F,\eps}=q_{n_{F,\eps}}\tilde{p}_{F,\eps/2}q_{n_{F,\eps}}$.
Then $\|{p}_{F,\eps}+\left(\begin{smallmatrix}
    \Id_{H}-I_n&0\\0&0\end{smallmatrix}\right)-e'_F\|<\eps/4$, and according to remark
\ref{rem-closed}, ${p}_{F,\eps}$ is as a summand of ${p}_{F,\eps}+\left(\begin{smallmatrix}
    \Id_{H}-I_n&0\\0&0\end{smallmatrix}\right)$ a $\eps$-projector in
$M_{2n_{F,\eps}}(C(\Ga/\Ga_i)\rtm\Ga)$. Moreover, we have
$$\mu_{\Ga,C(\Ga/\Ga_i),\text{max},0}(x_F)=[\phi_\eps({p}_{F,\eps})]-[I_{n_{F,\eps}}].$$

\smallskip

In the odd case, if   $F$ is a  K-cycle  of $\Psi_1^{\Ga_i}(\pr)$ with
propagation less than $s$, let us set using the
notations of the discussion following proposition
\ref{prop-ass-max} $q_{F}=1/2(F_{r,i}+\Id_{H})$. 
For $\eps$ in $(0,1/4)$ and $r$ positive, 
 let $l_{\eps,F,i}$ be the smallest integer such that
  $\displaystyle\sum_{l=l_{\eps,F,i}+1}^{+\infty}(\alpha_{r,i}+2))^l/l!<\eps/(3\alpha_{r,i}+6)$.
Let us define  $$u_{F,\eps}=\sum_{l=0}^{l_{\eps,F,i}}(2\imath\pi
q_{F})^l/l!-q_F\sum_{l=1}^{l_{\eps,F,i}}(2\imath\pi)^l/l!.$$
It is straightforward to check that
\begin{itemize}
\item $u_{F,\eps}-\Id_{H}$ is indeed an element of
    $C(\Ga/\Ga_i,\K(H))\rtm\Ga$ with propagation less than
    $l_{\eps,F,i}h_r(s)$.
\item $\|u_{F,\eps}-e^{2\imath\pi q_{F}}\|\leq \eps/3$.
\end{itemize}
In view of remark \ref{rem-closed}, $u_{F,\eps}$ is a $\eps$-unitary.
Moreover,
if  $x$  in $K^\text{top}_1(\Ga,C(\Ga/\Ga_i))$ comes  from
$x_F$ in $KK^\Ga_1(\pr,C(\Ga/\Ga_i))$ for $F$ in some
$\Psi^{\Ga_i}(\pr)$, then $\mu_{\Ga,C(\Ga/\Ga_i),\text{max},1}(x)$ is the class  of $u_{F,\eps}$ in $K_1(C(\Ga/\Ga_i)\rtm\Ga)$.
\begin{remark}
According to remark  \ref{rem-norm-equ}, for all $\eps$ in $(0,1/4)$,
$i$ positive integer and $s$ positive real, then the sets
$\{j_{\eps,F,i};\,F\text{ K-cycle of }\Psi_1^{\Ga_i}(\pr)\text{ of
  propagation less than }s\}$ and $\{l_{\eps,F,i};\,F\text{ K-cycle of }\Psi_1^{\Ga_i}(\pr)\text{ of
  propagation less than }s\}$ are bounded. Thereby,  if $J_{\eps,i,s}$ and
$L_{\eps,i,s}$ are respectively their suppremium, then for all
K-cycle $F$ of $\Psi_1^{\Ga_i}(\pr)$ with propagation less than $s$,
we get that $\widetilde{p_{F,\eps}}$ and ${p_{F,\eps}}$ have
propagation less than $(2J_{\eps,i,s}+1)h_r(s)$ and ${u_{F,\eps}}$ has
propagation less than $L_{\eps,i,s}h_r(s)$
\end{remark}
With notations of lemma \ref{sub-roe-resfin},
let $x$ be an element of $K_*^\text{top}(\Ga,B_\Ga)$. Under the
identification 
$KK^\Ga_*(\pr,B_\Ga)\cong\prod_{i\in\N}KK^\Ga_*(\pr,C(\Ga/\Ga_i))$
of proposition  \ref{prop-prod}, we can assume that $x$ comes from an
element $(x_{F_i})_{i\in\N}$, where
\begin{itemize}
\item $F_i$ is   K-cycle of $\Psi_1^{\Ga_i}(\pr)$ for
  every positive integer $i$;
\item there exists a real $s$ such that $F_i$ has propagation less
  than $s$ for every positive integer $i$. 
\end{itemize}
By viewing
$B_\Ga\rtm\Ga=\ell^\infty(\cup_{i\in\N}\Ga/\Ga_i,\K(H))\rtm\Ga$ as an
algebra of  multipliers  of
$\oplus_{i\in\N}C(\Ga/\Ga_i,\K(H))\rtm\Ga$, we see that $B_\Ga\rtm\Ga$
is indeed  a closed $*$-subalgebra of
$\prod_{i\in\N}\left(C(\Ga/\Ga_i,\K(H)\rtm\Ga\right)$. In particular, with above
notations, if $x$ is even and  since
$\|\lambda(S_j)\|\leq\|\lambda_\Ga(S_i)_{i\in\N}\|$ for all
uniformally bounded family
$S_i$ in $\prod_{i\in\N}\Psi^{\Ga_i}(\Ga_i)$ with propagation
uniformally bounded and all integer
$j$, then  the idempotent $(e_{F_i})_{i\in\N}$ and hence  the
projector $(e'_{F_i})_{i\in\N}$ belong  to $M_2(\widetilde{B_\Ga\rtm\Ga})$
and moreover,
$\mu_{\Ga,B_\Ga,\text{max},0}(x)=[(e'_{F_i})_{i\in\N}]-\left[\left(\begin{smallmatrix}
    \Id_{H}&0\\0&0\end{smallmatrix}\right)\right]$. Furthermore,  the family of integers
$(j_{\eps,F_i,i})_{i\in\N}$ is bounded and hence  $(\widetilde{p_{F_i,\eps}})_{i\in\N}$ and
  $(p_{F_i,\eps})_{i\in\N}$ are $\eps$-projector in
    $M_2(\widetilde{B_\Ga\rtm\Ga})$. Since
    $\|\widetilde{p_{F_i,\eps}}-e'_{F_i}\|<\eps$ for all integer $i$,
      we finally get that
      $(\phi_\eps({\widetilde{p_{F_i,\eps}}}))_{i\in\N}$ is a
      projector of $M_2(\widetilde{B_\Ga\rtm\Ga})$ homotopic to
      $(e'_{F_i})_{i\in\N}$ and hence 
\begin{eqnarray*}
\mu_{\Ga,B_\Ga,\text{max},0}(x)&=& [(\phi_\eps({\widetilde{p_{F_i,\eps}}}))_{i\in\N}]-[(\Id_H)_{i\in\N}]\\
&=&
[(\phi_\eps({p_{F_i,\eps}}))_{i\in\N}]-[(I_{n_{F_i,\eps}})_{i\in\N}].\end{eqnarray*}
In the same way, in the odd case we get that $(u_{F_i})_{i\in\N}$ is a
$\eps$-unitary of $\widetilde{B_\Ga\rtm\Ga}$ and
$$\mu_{\Ga,B_\Ga,\text{max},1}(x)=[(u_{F_i})_{i\in\N}].$$
\subsection{Asymptotic statements}\label{sec-asymp}
For any integer $i$ and any positive real $r,r',s,s'$ and any $\eps$
in $(0,1/72)$, let us consider the following statements
\begin{description}
\item[$\mathbf{QI_0(i,r,r',s,{\eps})}$] for any  (even) K-cycle $F$ of
  $\Psi_1^{\Ga_i}(\pr)$, then
  $(p_{F,\eps},n_{F,\eps})\sim_{18\eps,s}(0,0)$ implies that $x_{F}$
  lies in the kernel of the homomorphism 
$$KK_0^\Ga(\pr,C(\Ga/\Ga_i))\longrightarrow
KK_0^\Ga(P_{r'}(\Ga),C(\Ga/\Ga_i))$$ induced by the inclusion
$\pr\hookrightarrow P_{r'}(\Ga)$.
\item[$\mathbf{QI_1(i,r,r',s,{\eps})}$] for any (odd)  K-cycle $F$ of
  $\Psi_1^{\Ga_i}(\pr)$, then
  $u_{F,\eps}\sim_{\eps,s}\Id_{H} $ implies that $x_{F}$
  lies in the kernel of the homomorphism 
$$KK_1^\Ga(\pr,C(\Ga/\Ga_i))\longrightarrow
KK_1^\Ga(P_{r'}(\Ga),C(\Ga/\Ga_i))$$ induced by the inclusion
$\pr\hookrightarrow P_{r'}(\Ga)$.
\item[$\mathbf{QS_0(i,r,s,s',{\eps})}$] For any $\eps$-projector $p$ in
  some $M_k(C(\Ga/\Ga_i)\rtm\Ga)$ with propagation less than $s$, and
  any integer $n$, there exists a (even) K-cycle $F$ of
  $\Psi_1^{\Ga_i}(\pr)$ such that $(p_{F,\eps},n_{F,\eps})\sim_{18\eps,s}(p,n)$.
\item[$\mathbf{QS_1(i,r,s,s',{\eps})}$] For any $\eps$-unitary  $u$ in
  $C(\Ga/\Ga_i,\K(H)\widetilde{\rtm}\Ga$  with propagation less than
  $s$,  there exists a (odd) K-cycle $F$ of
  $\Psi_1^{\Ga_i}(\pr)$ such that $u_{F,\eps}\sim_{\eps,s} u$.
\end{description}
\begin{remark}\label{rem-asympy-st} It is straightforward to check that if two $\eps$-projectors
  are $\eps$-closed, then they are homotopic as $18\eps$-projectors and
  hence conditions $QI_0$ and $QS_0$ do not depend on a particular  choice of
    sequences of polynomial functions $(P_n)_{n\in\N}$ and
    $(Q_n)_{n\in\N}$ used in the definition of $p_{F,\eps}$. Moreover,
    replacing $n_{F,\eps}$ by any integer $n$ with $n\geq n_{F,\eps}$
    and $p_{F,\eps}$ by $q_n\widetilde{p_{F,\eps/2}}q_n$ does not either affect
    conditions $QI_0$ and $QS_0$.
\end{remark}
\begin{theorem}\label{theo-QI} Let $\Ga$ be   a  finitely generated group  residually finite  with
respect to a family 
$\Ga_0\supset\Ga_1\supset\ldots\Ga_n\supset\ldots$ of  normal finite
index subgroups and let 
 $l$ be in $\{0,1\}$. The following statements are equivalent:
\begin{enumerate}
\item  For any positive real $r$ the
following condition holds : there is an $\eps$ in $(0,1/72)$  such
that for any positive real $s$, there exists  an
integer $j$ and a 
positive real $r'$ for which  $QI_l(i,r,r',s,\eps)$ is true for all
$i\geq j$.
\item The maximal coarse assembly map $\mu_{\x,\text{max},l}$ is
  injective.
\end{enumerate}
\end{theorem}
\begin{proof}  Let us give the proof in the
even case, the odd one been quite similar (even simpler). In view of theorem \ref{thm-iso-khom},  condition (ii)
  is equivalent to injectivity of $\mu_{\Ga,A_\Ga,\mx,0}$. Let us
  prove that condition (i) implies injectivity of $\mu_{\Ga,A_\Ga,\mx,0}$.
According to remark \ref{rem-intertwin}, this amounts to prove that
for any  $x$ in $K^\text{top}_0(\Ga,B_\Ga)$, then the condition
$\mu_{\Ga,B_\Ga,\text{max},0}(x)\in K_0(B_{\Ga,0}\rtm\Ga)$ implies
that $x$ belongs indeed to  $K^\text{top}_0(\Ga,B_{\Ga,0})$. Up to
replace $\cup_{i\in\N}\Ga/\Ga_i$ by $\cup_{i\geq i_0}\Ga/\Ga_i$
for some integer $i_0$, we can actually assume that
$\mu_{\Ga,B_\Ga,\text{max},0}(x)=0$.
Suppose that $x$ comes from an element $(x_{F_i})_{i\in\N}$ in some
$KK^\Ga_0(\pr,B_\Ga)\cong \prod_{i\in\N}KK^\Ga_0(\pr,C(\Ga/\Ga_i))$
for  $r$  positive real, where $(F_i)_{i\in\N}$ is a family of K-cycles
in $\prod_{i\in\N}\Psi_1^{\Ga_i}(\pr)$ with propagation uniformally bounded. Then
there exist integers $k$ and $n$ and a projector homotopy  in
$M_{n+k+2}(\widetilde{B_\Ga\rtm\Ga})$ between
$\left(\begin{smallmatrix}
    (\phi_\eps(p_{F_i,\eps}))_{i\in\N}&0\\0&p_{n,k}\end{smallmatrix}\right)$ and$\left(\begin{smallmatrix}(I_{n_{F_i,\eps}})_{i\in\N}&0\\0&p_{n,k+1}\end{smallmatrix}\right)$,
where 
$p_{n,k}$ is the projector $\left(\begin{smallmatrix}
    I_n&0\\0&0\end{smallmatrix}\right)$ of 
  $M_{n+k}(\widetilde{B_\Ga\rtm\Ga})$. Hence we can find a homotopy of
  $18\eps$-projector $P:[0,1]\to M_{n+k+2}(\widetilde{B_\Ga\rtm\Ga})$
  between $\left(\begin{smallmatrix}
    (p_{F_i,\eps})_{i\in\N}&0\\0&p_{n,k}\end{smallmatrix}\right)$
and$\left(\begin{smallmatrix}(I_{n_{F_i,\eps}})_{i\in\N}&0\\0&p_{n,k+1}\end{smallmatrix}\right)$
  such that for some $s$   real, $P(t)$ has propagation less than $s$
  for every $t$ in $[0,1]$. From this,   by
  using for every
  integer $j$ the projection
  $$B_\Ga\rtm\Ga=\ell^\infty(\cup_{i\in\N},\K(H))\rtm\Ga\to
  C(\Ga/\Ga_j,\K(H))\rtm\Ga$$ and proceeding as we did in section
  \ref{subsec-prop} to obtain  ${p_{F,\eps}}$ from
  $\widetilde{p_{F,\eps/2}}$, we get that
  $(p_{F_j,\eps},n_{F_j,\eps})\sim_{18\eps,s}(0,0)$. If  $\eps$ is in
  $(0,1/72)$ and $j$ is an  integer satisfy  the assumptions  of the
  theorem for  $s$  as above, then there exists a $r'$ such that
  $x_{F_i}$ lies in the kernel of
 $KK^\Ga_0(\pr,C(\Ga/\Ga_i))\to KK^\Ga_0(P_{r'}(\Ga),C(\Ga/\Ga_i))$
 for all integer $i\geq j$. This implies that   $x$ comes indeed from an element in
 $\bigoplus_{i=0}^{j-1} KK^\Ga_0(P_{r'}(\Ga),C(\Ga/\Ga_i))$ and  hence
 belongs to $\bigoplus_{i\in\N}K^\top_0(\Ga,C(\Ga/\Ga_i))$.

Conversely, assume that for some positive real $r$, then for any $\eps$
in $(0,1/72)$ there exists a positive real $s$ such that for every
integer $j$ and positive real $r'$, there exists an integer $i$ with $i\geq j$ for which
$QI_0(i,r,r',s,\eps)$ does not hold. Let us prove that
$\mu_{\Ga,\a,\mx,0}$ is not injective. If $r$ is as above, let us fix
$\eps$ in $(0,1/72)$ and $(r'_n)_{n\in\N}$ an increasing and unbounded
sequence of positive reals. Then we can find an increasing sequence
$(j_i)_{i\in\N}$ of integers, and for each integer $i$ a K-cycle
$F_{j_i}$ in $\Psi^{\Ga_{j_i}}_1(\pr)$ such that
$(p_{F_{j_i},\eps},n_{F_{j_i},\eps})\sim_{18\eps,s}(0,0)$ and $x_{F_{j_i}}$ does
not belong to the kernel of $KK^\Ga_0(\pr,C(\Ga/\Ga_{j_i}))\to
KK^\Ga_0(P_{r'_i}(\Ga),C(\Ga/\Ga_{j_i}))$.
 By using a cut-off
function for the action of $\Ga$ on $\pr$  and in view of remark
\ref{rem-asympy-st},  we can actually assume that the family
$(F_{j_i})_{i\in\N}$ as propagation uniformally bounded.
Define for any integer $k$ the K-cycle $F_k$ of $\Psi_1^{\Ga_k}(\pr)$
to be $F_{j_i}$ if $k=j_i$ for some integer $i$ and $\Id_{H_\pr}$
otherwise. Let $x$ be the element of $K^\top_0(\Ga,B_\Ga)$ arising
from $(x_{F_i})_{i\in\N}$. We clearly have
$\mu_{\Ga,B_\Ga,\mx,0}(x)=0$ and  $x$ does not sit in
$\bigoplus_{i\in\N}K^\top_0(\Ga,C(\Ga/\Ga_i))$. Hence, the image of $x$  under the
epimorphism 
$K^\top_0(\Ga,B_\Ga)\to K^\top_0(\Ga,A_\Ga)$ is a non vanishing
element of the kernel of $\mu_{\Ga,\a,\mx,0}$.
\end{proof}
\begin{corollary}\label{cor-QI}
If $\Ga$ is residually finite, finitely generated and uniformally
embeddable into a Hilbert space, then $\Ga$ satisfies condition (i) of
theorem \ref{theo-QI}.
\end{corollary}
\begin{remark} As already mentionned, under the assumption of
  corollary \ref{cor-QI}, the reduced assembly map
  $\mu_{\Ga,\a,\red,*}$ is injective. Moreover, the group $\Ga$ is
  K-exact and  hence, in view of the   proof
  of theorem \ref{theo-QI} and if  we consider conditions $QI_0$ and $QI_1$ with 
 reduced assembly maps instead of maximal ones,  we get in this
 setting  an analogue of corollary \ref{cor-QI} for $\Ga$.
\end{remark}
\begin{theorem}\label{thm-QS} Let $\Ga$ be   a  finitely generated group  residually finite  with
respect to a family 
$\Ga_0\supset\Ga_1\supset\ldots\Ga_n\supset\ldots$ of  normal finite
index subgroups and let 
 $l$ be in $\{0,1\}$.  Assume that there exists $\eps$ in $(0,1/72)$
such that the following condition holds:  for every positive
real $s$, there exist positive real $r$ and $s'$ and an integer $j$
such that $QS_l(i,r,s,s',\eps)$ is true for all integer $i\geq j$. Then
$\mu_{\pr,\mx,l}$ is surjective.
\end{theorem}
\begin{proof}
As for injectivity, it is enought in view of theorem \ref{thm-iso-khom}   to prove that
$\mu_{\Ga,A_\Ga,\mx,l}$ is surjective and 
according to remark \ref{rem-intertwin}, this amounts to prove that
for any  $z$ in $K_l(B_\Ga\rtm\Ga)$, there exists an element  $x$ in
$K^\text{top}_l(\Ga,B_\Ga)$ such that 
$\mu_{\Ga,B_\Ga,\text{max},l}(x)-z$ belongs  to
$K_l(B_{\Ga,0}\rtm\Ga)$. As before, we give the proof in the even
case, the odd case being quite similar.
Recall that we have fixed an identification
$\K(H)\cong\overline{\cup_{n\in\N}M_n(\C)}$. It is then
straightforward to check that every element in $K_0(B_\Ga\rtm\Ga)$ can
be written down as the difference of the classes of projector that
belongs to $\left(\prod_{i\in\N}M_{n_i}(C(\Ga/\Ga_i))\right)\rtm\Ga$
for some sequence of integers $(n_i)_{i\in\N}$.
Let $p=(p_i)_{i\in\N}$ be a such projector viewed as an element of
$\prod_{i\in\N}\left(M_{n_i}(C(\Ga/\Ga_i))\rtm\Ga\right)$. Let $\eps$
be as in the assumption of the theorem. We can indeed assume that there
exists a positive real $s$ and $(q_i)_{i\in\N}$ an $\eps$-projector of
$\left(\prod_{i\in\N}M_{n_i}(C(\Ga/\Ga_i)\right)\rtm\Ga\subset\prod_{i\in\N}\left(M_{n_i}(C(\Ga/\Ga_i))\rtm\Ga\right)$
with propagation less than  $s$ and such that
$p=(\phi_\eps(q_i))_{i\in\N}$. Let  $r$ and $s'$ be  positive reals and
let  $j$ be  a
positive integer such that $QS_0(i,r,s,s',\eps)$ is true for every
integer $i\geq j$, i.e there exists a
 K-cycle $F_i$ in $\Psi_1^{\Ga_i}(\pr)$ such that
$(p_{F_i,\eps},n_{F_i,\eps})\sim_{18\eps,s}(q_i,0)$. By using a cut-off
function for the action of $\Ga$ on $\pr$  and in view of remark
\ref{rem-closed},  we can actually assume that the family
$(F_i)_{i\in\N}$ as propagation uniformally bounded.    If we set
$F_i=\Id_{H_\pr}$ for every positive integer $i$ with  $i\leq j-1$ and
then consider the element $x$ of $K^\top_0(\Ga,B_\Ga)$ coming from
$(x_{F_i})_{i\in\N}\in KK^\Ga_0(\pr,B_\Ga)$, we get that  
$\mu_{\Ga,A_\Ga,\mx,0}(x)-[p]$ belongs to
  $\bigoplus_{i\in\N}K_0(C(\Ga/\Ga_i)\rtm\Ga)$.
\end{proof} 
As we shall see, up to a slight modification in the sequence of
finite index  normal
subgroups in $\Ga$, we get a converse result for theorem
\ref{thm-QS}. This  allows in particular
to deal at least with group that satisfies the strong Baum-Connes
conjecture.
Let us set for any $i$ integer $X_i(\Ga)=\coprod_{j\geq i}\Ga/\Ga_j$
and $X^\infty(\Ga)=\coprod_{i\in\N}X_i(\Ga)$ provided with the action
of $\Ga$ inherited by the action on $\x$. Let us equip
$X^\infty(\Ga)$ with a $\Ga$-invariant metric $d$ such that the
restriction of $d$ to each $X_i(\Ga)$ coincides with the metric on
$\x$ and $d(X_i(\Ga),X_j(\Ga))\geq i+j$ for every integer $i$ and
$j$. Let us set $A^\infty_\Ga=\ell^\infty(\xin,\K(H))/C_0(\xin,\K(H))$.
The space $\xin$ is indeed construct in the same way as $\x$ by
considering the sequence of finite index normal subgroups
$\Ga_0\supset\Ga_1\supset\Ga_1\supset\Ga_2\supset\Ga_2\supset\Ga_2\supset\Ga_3\ldots$
and hence, according to theorem \ref{thm-equ}, we get 
\begin{proposition}\label{prop-equiv-inf} The following assertions are equivalent
\begin{enumerate}
\item   
The maximal coarse assembly map $$\mu_{X^\infty(\Ga),\text{max},*}:\lim_r
K_*(P_r(\xin),\C)\to  K_*(C_\mx(\xin))$$ is surjective.
\item the maximal assembly map   
$$\mu_{\Ga,\ai,\text{max},*}:K^\text{top}_*(\Ga,A^\infty_\Ga)\to K_*(\ai\rtm\Ga)$$
is surjective.
\end{enumerate}
\end{proposition}
We have of course  analogous statements for injectivity and
isomorphism.
We are now in position to give a weak converse result for theorem
\ref{thm-QS}.
\begin{theorem}\label{thm-QS-converse}  Let $\Ga$ be   a  finitely generated group,  residually finite  with
respect to a family 
$\Ga_0\supset\Ga_1\supset\ldots\Ga_n\supset\ldots$ of  normal finite
index subgroups and let 
 $l$ be in $\{0,1\}$.  Assume that the maximal coarse assembly map $$\mu_{X^\infty(\Ga),\text{max},l}:\lim_r
K_l(P_r(\xin),\C)\to  K_l(C_\mx(\xin))$$  is onto.
Then  there exists $\eps$ in $(0,1/72)$
such that the following condition is satisfied:  for every positive
reals $s$, there exist  positive reals $r$ and $s'$ and an integer $j$
such that $QS_l(i,r,s,s',\eps)$ is true for all integer $i\geq
j$.\end{theorem}
\begin{proof}As before we give the prove for the even case.
Assume that for all $\eps$ in $(0,1/72)$ there exists a positive real
$s$ such that for all positive reals $r$ and $s'$ and integer $j$
there exists an integer $i$ with $i\geq j$ for which $QS_0(i,r,s,s',\eps)$
does not holds. In view of proposition \ref{prop-equiv-inf}, let us show that
$\mu_{\Ga,\ai,\text{max},0}$ is not
onto.
Let us fix $\eps$ in $(0,1/72)$ and $(s'_i)_{i\in\N}$ and
$(r_i)_{i\in\N}$ increasing and unbounded sequences of positive reals.
Then for each integer $k$, there exist  an increasing sequence of
integers $(j^k_i)_{i\in\N}$, an $\eps$-projector $q_{j^k_i,k}$ with
propagation less than $s$ in some
$M_{n_{j^k_i,k}}(C(\Ga/\Ga_{j^k_i})\rtm\Ga)$ and an integer $m_{j^k_i,k}$ such that 
there is no K-cycle $F$ in $\Psi^{\Ga_{j^k_i}}(P_{r_k}(\Ga))$ for which
$(x_{F,\eps},n_{F,\eps})\sim_{18\eps,s'_i}(q_{j^k_i,k},m_{j^k_i,k})$.
For $j$ and $k$ integers such that $j\geq k$, define 
\begin{itemize}
\item $q_{j,k}$ to be $q_{j^k_i,k}$ if $j=j^k_i$ for some integer $i$
    and  $q_{j,k}=0$ otherwise.
\item $m_{j,k}$ to be $m_{j^k_i,k}$ if $j=j^k_i$ for some integer $i$
    and  $m_{j,k}=0$ otherwise.
\end{itemize}
Let us set $\bi=\ell^\infty(\xin,\K(H))$ and
$\bio=C_0(\xin,\K(H))$. As in the proof of theorem \ref{thm-QS}
surjectivity will fail if there is  no $x$ in $K^\top_0(\Ga,\bi)$ such
that 
$\mu_{\Ga,\bi,\mx,0}(x)-[(\phi_\eps(q_{j,k})_{j\in\N,\,k\leq j}]+[(I_{m_{j,k}})_{j\in\N,\, k\leq j}]$ lies in $K_0(\bio\rtm\Ga)$.
Suppose that such an  $x$ exists, coming from an element  $y$ in 
$KK^\Ga_0(\pr,\bi)$ and let us fix $k$ an integer such that
$r_{k}\geq r$. Define then  $y_k$ as the image of $y$
under the composition $$KK^\Ga(\pr,\bi)\to KK^\Ga(P_{r_k}(\Ga),\bi)\to
KK^\Ga(P_{r_k}(\Ga),\ai),$$ where the first map is induced by the
inclusion $\pr\hookrightarrow P_{r_k}(\Ga)$ and the second by the
projection homomorphism $\ai\to \ell^\infty(X_k(\Ga),\K(H))$. We can
assume that  $y_k=(x_{F_{j,k}})_{j\in\N,\,k\leq j}$ where 
 $(F_{j,k})_{j\in\N,\,k\leq j}$ is a family of K-cycles in
 $\prod_{j\in\N,\,k\leq j}\Psi^{\Ga^j}_1(P_{r_k}(\Ga))$ with
 propagation uniformally bounded. Then
 $(\phi_\eps(p_{F_{j,k},\eps}))_{j\in\N,\,k\leq j}$ is a projector in
 $K_0(\ell^\infty(X_k(\Ga))\rtm\Ga)$ and by naturality of the assembly
 map, we get that $[(\phi_\eps(p_{F_{j,k},\eps}))_{j\in\N,\,k\leq
   j}]-[(I_{n_{F_{j,k},\eps}})_{j\in\N,\,k\leq
   j}]-[(\phi_\eps(q_{j,k}))_{j\in\N,\,k\leq j}]+[(I_{m_{j,k}})_{j\in\N,\,
   k\leq j}]$ lies in $K_0(C_0(X_k(\Ga),\K(H))\rtm\Ga)$. By taking $k$ big
 enought, we can indeed assume that $$[(\phi_\eps(p_{F_{j,k},\eps}))_{j\in\N,\, k\leq
   j}]-[(I_{n_{F_{j,k},\eps}})_{j\in\N,\,k\leq
   j}]=[(\phi_\eps(q_{j,k})_{{j\in\N},\,{k\leq j}}]-[(I_{m_{j,k}})_{j\in\N,\,
   k\leq j}].$$
Thus, up to stabilisation, there is a homotopy 
in some
$M_n\left(\l^\infty(X_k(\Ga),\K(H))\widetilde{\rtm}\Ga\right)$ of
$18\eps$-projectors  $\left(\begin{smallmatrix}p_{F_{j,k},\eps}&0\\0& I_{m_{j,k}}\end{smallmatrix}\right)_{{j\in\N},\, k\leq
   j}$ and 
 $\left(\begin{smallmatrix}q_{j,k}&0\\0&I_{n_{F_{j,k},\eps}}\end{smallmatrix}\right)_{{j\in\N},\,{ k\leq j}}$ with finite
propagation 
between.
Proceeding as we did in section \ref{subsec-prop} to obtain
${p}_{F,\eps}$ from  $\tilde{p}_{F,\eps/2}$, we actually get that
there exists 
\begin{itemize}
\item a positive real $s'$; 
\item two  sequences of integers $(i_j)_{j\in\N,\,{ k\leq j}}$ and $(i'_j)_{j\in\N,\,{ k\leq j}}$; 
\item  a homotopy of
  $18\eps$-projector $P:[0,1]\to \left(\prod_{j\in\N,\,{ k\leq j}} M_{i'_j}(C(\Ga/\Ga_j))\right)\rtm\Ga$
  between $\left(\begin{smallmatrix}
    p_{F_{j,k},\eps}&0\\0&I_{i_j+m_{j,k}}\end{smallmatrix}\right)_{j\in\N,\,k\leq j}$
and$\left(\begin{smallmatrix}q_{j,k}&0\\0&I_{i_j+n_{F_{j,k},\eps}}\end{smallmatrix}\right)_{j\in\N,\,k\leq j}$
  such that  $P(t)$ has propagation less than $s'$
  for every $t$ in $[0,1]$.
\end{itemize}
If   $i$ is an integer such that $s'_i\geq s'$ and $j^k_i\geq k$ , then 
$(p_{F_{j^i_k,k},\eps},n_{F_{j^i_k,k},\eps})\sim_{18\eps,s'_i}(q_{j^i_k,k},m_{j^k_i,k})$,
which is in contradiction with the way we have chosen 
$(q_{j^k_i,k})_{(i,k)\in\N^2}$ and $(m_{j^k_i,k})_{(i,k)\in\N^2}$.
\end{proof}

\begin{corollary} Let $\Ga$ be   a  finitely generated group  residually finite  with
respect to a family 
$\Ga_0\supset\Ga_1\supset\ldots\Ga_n\supset\ldots$ of  normal finite
index subgroups and let 
 $l$ be in $\{0,1\}$.
If $\Ga$ satisfies the strong Baum-Connes conjecture (for example if
$\Ga$ has the Haagerup property) then  there exists $\eps$ in $(0,1/72)$
such that the following condition is satisfied:  for every positive
reals $s$, there exist  positive real $r$ and $s'$ and an integer $j$
such that $QS_l(i,r,s,s',\eps)$ is true for all integer $i\geq
j$. 
\end{corollary}
Since group satisfying the strong Baum-Connes conjecture are
K-amenable \cite{tumoy}, the same result holds if we replace in the definition of
conditions $QS_0$ and $QS_1$ maximal assembly maps by the reduced
one. More generally, in this setting, the  analogue of the hypothesis
of theorem \ref{thm-QS-converse} implies the surjectivity of the reduced Baum-Connes assembly map
$\mu_{\Ga,\bi,\red,l}:K^\top_l(\Ga,\bi)\to K_l(\bi\rtm\Ga)$ for $l$ in
$\{0,1\}$ for K-exact groups (in particular for groups that embed
uniformly in a Hilbert space). 

\section*{Acknowledgements}
The first author wants to thank Guoliang Yu and the Department of
Mathematics at Vanderbilt University, Nashville,  for their warm hospitality during
a visit in Nashville, in October 2007 and February 2009, where much of
this paper was written.
He  is also indebted to  Vincent Lafforgue and Georges Skandalis
for some  helpful discussions.
The second author would like to thank Herv\'e Oyono-Oyono and the
Department of Mathematics of Clermont-Ferrand for their great
hospitality during  a visit in June 2006 where this work was
initiated. The second author also acknowledges partial financial
support from NSF and CNSF.

{\small
\bibliographystyle{plain}
\bibliography{expanders}
}

\end{document}